\pdfoutput=1
\documentclass[a4paper,11pt,pdf,reqno,french,english]{amsart}
\usepackage{amsmath, amsfonts, amssymb, amsthm, thmtools, wasysym, graphics, graphicx, xcolor, frcursive,comment,bbm}
\usepackage[shortlabels]{enumitem}
\usepackage{caption}
\usepackage{subcaption} 
\usepackage{lscape}
\usepackage{tikz-cd}
\usepackage{comment}
\usepackage[utf8]{inputenc}
\usetikzlibrary{braids} 

\makeatletter
\newcommand{\thickhline}{%
    \noalign {\ifnum 0=`}\fi \hrule height 1pt
    \futurelet \reserved@a \@xhline
}

\definecolor{darkblue}{rgb}{0.0,0,0.7} 
 
\definecolor{darkred}{rgb}{0.7,0,0}
\usepackage[
backend=biber,
style=numeric,
isbn=false,
maxnames=50,
minalphanames=6,
maxalphanames=10
]{biblatex}
\AtEveryBibitem{
	\clearfield{urlyear}
	\clearfield{urlmonth}
}
\usepackage[all]{xy}
\usepackage{mathtools}
\usepackage{utfsym}
\usepackage[T1]{fontenc}
\usepackage{stmaryrd}
\usepackage{float}
\usepackage{adjustbox}

\usepackage{vmargin}            
\setmarginsrb{3.5cm}{2.2cm}{3.5cm}{2.2cm}{0cm}{0.6cm}{0cm}{0.3cm}

\usepackage{caption,lipsum}
\usepackage{graphicx}                  
\usepackage{pstricks,pst-plot,pst-text,pst-tree,pst-eps,pst-fill,pst-node,pst-math}
\usepackage{setspace}

\usepackage{multicol}

\usepackage{yhmath}
\usepackage[colorlinks=true,citecolor=cyan,backref=none, hypertexnames=false]{hyperref}
\hypersetup{
    citecolor=magenta,
    colorlinks=true,
    linkcolor=blue,
    filecolor=cyan,      
    urlcolor=magenta,
}

\numberwithin{equation}{section}

\newcommand{\Z}{\mathbb{Z}}

\newcommand{\R}{\mathbb{R}}
\newcommand{\llangle}{\langle\langle}
\newcommand{\rrangle}{\rangle\rangle}
\newcommand{\N}{\mathbb{N}}
\newcommand{\C}{\mathbb{C}}
\newcommand{\B}{\mathcal{B}}

\newcommand{\K}{\mathbb{K}}

\newcommand{\s}{\sigma}

\newcommand{\xto}{\xrightarrow}

\makeatletter
\newcommand*\bigcdot{\mathpalette\bigcdot@{.5}}
\newcommand*\bigcdot@[2]{\mathbin{\vcenter{\hbox{\scalebox{#2}{$\m@th#1\bullet$}}}}}
\makeatother
\theoremstyle{plain}

\newcommand{\bbar}[1]{\overline{#1}}

\newtheorem{theorem}{Theorem}[section]
\newtheorem{lemma}[theorem]{Lemma}
\newtheorem{proposition}[theorem]{Proposition}
\newtheorem{corollary}[theorem]{Corollary}

\theoremstyle{definition}

\newtheorem{definition}[theorem]{Definition}

\newtheorem{remark}[theorem]{Remark}

\newtheorem{notation}[theorem]{Notation}

\newcommand{\intv}[1]{[\![#1]\!]}
\newcommand{\jg}[1]{\big( \begin{smallmatrix} #1\end{smallmatrix}\big)}
\newcommand{\Hh}{\mathcal{H}}
\newcommand{\Ll}{\mathcal{L}}
\newcommand{\muls}[1]{ \left\{\kern-0.6em \left\{ #1\right\}\kern-0.6em \right\} }

\title[Generalized $J$-groups, $J$-braid groups and Seifert link groups]{Generalized $J$-groups, $J$-braid groups and Seifert link groups}

\author{Owen Garnier}
\address{O.~Garnier: Departamento de \'Algebra e Instituto de Matem\'aticas de la Universidad de Sevilla, Spain}
\author{Igor Haladjian}
\address{I.~Haladjian: Institut Denis Poisson, CNRS UMR 7019, Faculté des Sciences et Techniques, Université de Tours, Parc de Grandmont, 
37200 TOURS, France}
\DeclareRobustCommand{\SkipTocEntry}[5]{}

\begin{document}

\thispagestyle{empty}
\begin{abstract}
The family of $J$-groups was introduced by Achar and Aubert with the goal of providing Coxeter-like combinatorial tools for studying rank 2 complex reflection groups. However, $J$-groups lack an explicit presentation with abstract reflections as generators. This gap was filled by Gobet, and later by the second author, for the subfamily of so-called $J$-reflection groups. The obtained presentations then gave rise to a concept of $J$-braid group, which coincides with the link groups of torus necklaces.

In this paper we study a generalization of $J$-groups. We determine which of these groups are finitely generated. We show that, as for classical $J$-groups, the family of finite generalized $J$-groups coincides with the family of rank 2 complex reflection groups. We also show that finitely generated generalized $J$-groups coincide with what we call the torsion quotients of $J$-braid groups. We deduce explicit presentations for all finitely generated generalized $J$-groups, where the generators are abstract reflections. We also complete the classification of these groups up to reflection isomorphism.

As a byproduct of these results, we obtain that a quotient of a Seifert link group obtained by adding torsion to meridians somehow determines the link up to isotopy. Moreover, such a quotient is finite if and only if it is isomorphic to a complex reflection group of rank two.
\end{abstract}

\maketitle

\tableofcontents

\section{Introduction}
It is well-known since the work of Coxeter that the family of (finite) real reflection group coincides with the family of finite Coxeter groups (see \cite{CoxeterClassification},\cite{CoxeterClassification2}). This allows for the use of the rich combinatorics of Coxeter groups in the study of real reflection groups. Unfortunately, no analogue of the result of Coxeter is known for complex reflection groups. The search for such an analogue, at least in the case of rank 2 complex reflection groups, is the motivation for Achar \& Aubert's paper \cite{AA}, where they defined the family of $J$-groups:

Let $k,n,m$ be positive integers. The group $J(k,n,m)$ is defined by the group presentation
\[\langle s,t,u~|~s^k=t^n=u^m=1,stu=tus=ust\rangle,\]
and it is called a \emph{parent $J$-group}. Now, if $k',n',m'$ are positive and pairwise coprime integers which respectively divide $k,n,m$, the \emph{$J$-group} $J\jg{k&n&m\\k'&n'&m'}$ is defined as the normal closure of $\{s^{k'},t^{n'},u^{m'}\}$ in $J(k,n,m)$.

We know that every rank two complex reflection group is the normal closure of some reflections in one of the groups $G(2c,2,2)$, $G_7$, $G_{11}$ or $G_{19}$ in Shephard-Todd's notation. Inspecting \cite[Tables 1-4]{BMR}, one sees that the aforementioned groups are parent $J$-groups and so every rank two complex reflection group is a (finite) $J$-group. The main result of \cite{AA} is the converse statement:

\begin{theorem}\cite[Theorem 1.2]{AA}\label{Classification1}
A group is a rank two complex reflection group if and only if it is a finite $J$-group.
\end{theorem}

In \cite{VCRG}, generalizing the work of Gobet in \cite{GobetToric}, the second author studied the specific case of $J$-groups where (at least) one of $k',n'$ or $m'$ is equal to 1. In this case we have a uniform presentation for these groups \cite[Theorem 2.29]{VCRG}. Moreover, these groups are classified up to reflection isomorphism \cite[Theorem 1.6]{VCRG} and their centers are known \cite[Theorem 1.5]{VCRG}. Later, in \cite{BraidsJIgor}, the second author also defined and studied $J$-braid groups attached to $J$-reflection groups, thus giving a Coxeter-like generalization of complex reflection groups of rank two and their braid groups. Moreover, connections with the triangle Von Dyck groups on the one hand and Seifert link groups on the other hand were obtained in \cite{Necklace}.

\vspace{1em} The purpose of this article is to define and study a family of groups we call \textit{generalized $J$-groups}, and to extend to this case previous results on $J$-reflection groups. Since our generalization contains the whole family of classical $J$-groups, our results also applies in this case. We also complete the correspondence between Seifert links and $J$-braid groups. As a byproduct, we give a new perspective on rank two complex reflection groups as the family of finite torsion quotients of Seifert link groups.

\begin{definition}[\textbf{Generalized $J$-groups}]\label{def:generalized_j_groups_intro}
Let $k_1,\dots,k_p\in\N_{\geqslant 2}\cup\{\infty\}$ and write $K:=(k_1,\dots,k_p)$. The group $J(K)$ is defined by the group presentation 
\begin{equation}\label{eq:pres_generalized_j_groups_intro}
    \langle s_1,\dots,s_p\,|\, s_1s_2\cdots s_p=s_2\cdots s_ps_1=\cdots=s_ps_1\cdots s_{p-1},\, s_i^{k_i} \text{ if $k_i<\infty$}\rangle.
\end{equation}
and it is called a \emph{parent generalized $J$-group}. Now, if $k_1',\dots,k_p'\in\N_{\geqslant 1}$ are respective divisors of $k_1,\dots,k_p$, and writing $K'=(k_1',\dots,k_p')$, the \emph{generalized $J$-group} $J\binom K{K'}$ is defined as the normal closure of $\{s_i^{k_i'}~|~i\in \intv{1,p}\}$ in $J(K)$.
\end{definition}

Notice that for $p\leqslant 3$, $k_1,k_2,k_3<\infty$, and $k_1',k_2',k_3'$ pairwise coprime, we recover the notion of classical $J$-groups. 

Following the definition given by Gobet in \cite{GobetToric}, we call conjugates in $J(K)$ of non-trivial powers of elements in $\{s_i^{k_i'}~|~i\in \intv{1,p}\}$ the \textit{reflections of $J\binom K{K'}$}. Still following \cite{GobetToric}, we say that an isomorphism between two generalized $J$-groups is a \textit{reflection isomorphism} if it induces a bijection on the set of reflections.

\vspace{1em}In Section \ref{sec:generalized_j-groups}, we study several group theoretic properties of generalized $J$-groups. Our approach takes advantage of a relationship between generalized $J$-groups and some alternating subgroups of particular Coxeter groups. We call these Coxeter groups polygonal Coxeter groups (as a generalization of the triangle Coxeter groups) and we study them in Section \ref{sec:prelim_coxeter}.

We can separate three mutually exclusive cases in the study of generalized $J$-groups:
\begin{enumerate}[(1)]
\item The generalized parent $J$-group is a finite group.
\item The generalized parent $J$-group is infinite, but the generalized $J$-group has finite index in its generalized parent $J$-group.
\item The generalized $J$-group is nontrivial and has infinite index in its generalized parent $J$-group.
\end{enumerate}
These three cases of course cover all generalized $J$-groups, and they are recognizable from looking at the generalized $J$-group alone:

\begin{theorem}[Proposition \ref{prop:finite-type_=_finiment_engendre} and Corollary \ref{cor:finite_generalised_j_groups}]\label{prop:finite_type_=_finiment_engedre_intro}\hfill
\begin{itemize}
\item A generalized $J$-group is finite if and only if we are in case $(1)$.
\item A generalized $J$-group is infinite and finitely generated if and only if we are in case $(2)$.
\item A generalized $J$-group is infinitely generated if and only if we are in case $(3)$. 
\end{itemize}
\end{theorem}

Since we easily determine which generalized parent $J$-groups are finite (see Lemma \ref{lem:finite_index_J_or_coxeter}), we obtain the list of all finite generalized $J$-groups. In particular no new finite groups arise as generalized $J$-groups:

\begin{corollary}[Corollary \ref{cor:finite_generalized_j_groups_are_crgs}]\label{cor:finite_generalized_j_groups_intro}
A group is isomorphic to a finite generalized $J$-group if and only if it is isomorphic to a complex reflection group of rank two.
\end{corollary}

Specializing to the case of classical $J$-groups, Theorem \ref{prop:finite_type_=_finiment_engedre_intro} answers \cite[Question 1.1.2]{TheseIgor}. 

Then we compute the center of generalized $J$-groups. We obtain that the center of a generalized $J$-groups allows us to recognize cases $(1)$, $(2)$ and $(3)$:

\begin{proposition}[Proposition \ref{order_center_generalized_j_groups}]\hfill
\begin{itemize}
    \item In case $(1)$, the center is either a finite cyclic group or a  direct product of two finite cyclic groups.
    \item In case $(2)$, the center is infinite cyclic.
    \item In case $(3)$, the center is trivial.
\end{itemize}
\end{proposition}
This result again generalizes \cite[Theorem 1.5]{VCRG} for $J$-reflection groups. Lastly, we get a description of torsion elements in infinite generalized $J$-groups (i.e. cases $(2)$ and $(3)$).

\begin{corollary}[Corollary \ref{cor:torsion_generalized_j_groups}]
If a generalized $J$-group is infinite, its torsion elements are precisely its finite order reflections.
\end{corollary}

\vspace{1em} In Section \ref{sec:j-braid_groups_and_torsion_quotients}, we give the definition of $J$-braid groups and their torsion quotients. The family of $J$-braid groups was originally defined by the second author in \cite{BraidsJIgor} as a combinatorial generalization of rank 2 complex braid groups. The family of $J$-braid groups actually consists of 4 related subfamilies depending on two positive integer parameters $n,m$. 

The group $\B_*^*(n,m)$ has generators $\{x_1,\ldots,x_n,y,z\}$ and a presentation depending on $n$ and $m$ (see Definition \ref{BraidPresentation1} for the actual presentation). The conjugates of the generators are called the (braid) reflections of $\B_*^*(n,m)$. The groups $\B^*(n,m)$ (resp. $\B_*(n,m)$, $\B(n,m)$) are defined as the quotient of $\B_*^*(n,m)$ by the normal closure of $y$ (resp. of $z$, of $\{y,z\}$).

After some reminders on $J$-braid groups, we introduce the other main concept of the article, which is that of \emph{torsion quotient of $J$-braid group}. Concretely, a torsion quotient of a given $J$-braid group $\B$ is the quotient of $\B$ by the normal closure of some (non-trivial) powers of its braid reflections. As in the case of generalized $J$-groups, we define reflections as conjugates of nontrivial powers of the generators, and we extend to this context the definition of reflection isomorphism.

The main result of Section \ref{sec:j-braid_groups_and_torsion_quotients} is that torsion quotients of $J$-braid groups are reflection isomorphic to finite-type generalized $J$-groups. In order to be able to give a precise statement, we need some notation. A torsion quotient of a $J$-braid group $\B$ is determined by the exponents of the generators of which we consider the normal closure. The exponents of the generators of the form $x_i$ form a tuple $K$, while the integers $b,c$ are the respective exponents for the generators $y$ and $z$ (if $\B$ admits such generators). With this notation we have the following result:

\begin{theorem}[Corollary \ref{cor:embedding_torsion_quotients}]\label{theo:embedding_torsion_quotients_intro}
Let $n,m$ be positive integer with gcd $d$, and let $n':=\frac{n}{d}$, $m':=\frac{m}{d}$. Let also $K=(k_1,\dots,k_d)\in(\N_{\geqslant 2}\cup \{\infty\})^d$ and $b,c\in(\N_{\geqslant 2}\cup\{\infty\})$. We have the following reflection isomorphisms:
\begin{itemize}[itemsep=3pt]
\item $\B^*_*(n,m; K,b,c)\cong_{\mathrm{Ref}}J\jg{K&bn'&cm'\\d\cdot 1&n'&m'}$.
\item $\B^*(n,m;K,c) \cong_{\mathrm{Ref}} J\jg{K&n'&cm'\\d\cdot 1&n'&m'}$.
\item $\B_*(n,m;K,b) \cong_{\mathrm{Ref}} J\jg{K&bn'&m'\\d\cdot 1&n'&m'}$.
\item $\B(n,m;K) \cong_{\mathrm{Ref}} J\jg{K&n'&m'\\d\cdot 1&n'&m'}$.
\end{itemize}
where $d\cdot 1$ denotes the $d$-tuple $(1,\ldots,1)$.
\end{theorem}

In particular, every torsion quotient of a $J$-braid group is reflection isomorphic to some finitely generated generalized $J$-group. Section \ref{sec:gen_j-groups_as_torsion_quotients} is devoted to the proof of the converse statement. This proof relies on a preliminary classification of finitely generated generalized $J$-groups (see Lemma \ref{lem:finite_type_generalized_list}). We then proceed in a case-by-case approach, the last case being handled partially by computer. 

\begin{theorem}[Theorem \ref{theo:generalized_j_groups_torsion_quotients}]\label{theo:generalized_j_groups_torsion_quotients_intro}
Every finitely generated generalized $J$-group is reflection isomorphic to a torsion quotient of $J$-braid group.
\end{theorem}

Since torsion quotients of $J$-braid groups are defined by explicit presentations with abstract reflections as generators, this result also provides explicit presentations of generalized $J$-groups with abstract reflections as generators.

One of the preliminary motivations for this paper was to determine the finite torsion quotients of $J$-braid groups, and more specifically of $J$-braid groups of the form $\B(n,m)$. The defining presentation of the group $\B(n,m)$ coincides with the defining presentation of the \emph{circular group} $G(n,m)$ studied by the first author in \cite{GarnierHoso}. This question had already been investigated in several cases by different authors. 
\begin{itemize}
    \item In \cite{CoxDihedralShephard}, Coxeter studied torsion quotients of dihedral Artin groups, that is, of groups of the form $\B(2,m;K)$. He classifies all finite such quotients (see \cite[Section 6]{CoxDihedralShephard}).
    \item The results of Achar and Aubert in \cite{AA} give in particular a complete list of the finite torsion quotients of $\B(3,3)$ (namely, finite parent $J$-groups).
    \item In \cite{GobetToric}, Gobet studied groups of the form $\B(n,m;K)$ with $n\wedge m=1$. They are called \emph{toric reflection groups} in \emph{loc. cit.} and they are denoted by $W(k,n,m)$. Gobet realizes these groups as $J$-groups, showcasing the first connection between $J$-groups and link groups. He moreover completes their classification up to reflection isomorphism, yielding in particular a classification of the finite groups $\B(n,m;K)$ with $n,m$ coprime. 
\end{itemize}

Theorem \ref{theo:embedding_torsion_quotients_intro} combined with the classification of finite generalized $J$-groups then yields the following result:

\begin{proposition}
There are no other finite torsion quotients of $J$-braid groups than those identified in \cite{CoxDihedralShephard},\cite{AA} and \cite{GobetToric}.
\end{proposition}
Since circular groups are particular cases of $J$-braid groups, this proposition also applies to circular groups.

\vspace{1em} Notice that the generalized $J$-groups appearing in Theorem \ref{theo:embedding_torsion_quotients_intro} share the property that the bottom set of parameters contains at most two entries different from $1$, in which case they are coprime. We call generalized $J$-groups sharing this property \emph{reduced $J$-groups}. Note that if $p=3$, the family of reduced $J$-groups coincides with that of $J$-reflection groups as defined by the second author in \cite{VCRG}. Combining Theorem \ref{theo:embedding_torsion_quotients_intro} and Theorem \ref{theo:generalized_j_groups_torsion_quotients_intro}, we obtain:

\begin{proposition}
Every finitely generated generalized $J$-group is reflection isomorphic to a reduced $J$-group.
\end{proposition}

In Section \ref{sec:classif}, we complete the classification of reduced $J$-groups up to reflection isomorphism:

\begin{proposition}[Proposition \ref{prop:classif_red_j-group}]\label{prop:classif_red_j-group_intro}
Let $W_1:=J\binom{K}{K'}$ and $W_2:=J\binom{L}{L'}$ be two reduced $J$-groups. 
\begin{enumerate}[(a)]
\item The group $W_1$ is abelian if and only if $p\leqslant 2$, in which case $W_1\cong_{\mathrm{Ref}} J\jg{k_1/k_1'\\1}$ or $W_1\cong_{\mathrm{Ref}}J\jg{k_1/k_1'&k_2/k_2'\\1&1}$ depending on the length of $K$.
\item If $W_1$ and $W_2$ are nonabelian, then we have $W_1\cong_{\mathrm{Ref}}W_2$ if and only if $p=q$ and if there is $\sigma\in S_p$ such that $\sigma(K)=L$ and $\sigma(K')=L'$.
\end{enumerate}
\end{proposition}

Combining this with Theorem \ref{theo:generalized_j_groups_torsion_quotients_intro} and Theorem \ref{theo:embedding_torsion_quotients_intro}, we may see reduced $J$-groups as canonical representatives of reflection isomorphism classes of generalized $J$-groups. Moreover, using Theorem \ref{theo:embedding_torsion_quotients_intro}, we can complete the classification of torsion quotients of $J$-braid groups up to reflection isomorphism. We skip the complete statement here (see Theorem \ref{theo:classif_tors_quot_j-braid} for a complete statement). We give here a consequence on the classification of $J$-braid groups

\begin{corollary}[Corollary \ref{cor:iso_quot_tors_implies_iso_j-braid}]\label{cor:iso_quot_j-braid_intro}
Let $\B,\B'$ be two $J$-braid groups, and let $W,W'$ be respective torsion quotients of $\B$ and $\B'$. If $W$ and $W'$ are reflection isomorphic, then so are $\B$ and $\B'$.
\end{corollary}

Specializing to the case of $J$-reflection groups, Corollary \ref{cor:iso_quot_j-braid_intro} solves \cite[Conjecture 3.2.3]{TheseIgor}.

\vspace{1em} Lastly, we detail the connection with Seifert links in Section \ref{sec:seifert_links}. In \cite{Necklace}, the second author showed that the family of $J$-braid groups coincides with that of torus necklace groups, an important family of Seifert link groups. Under this correspondence, braid reflections of the $J$-braid group correspond to meridians in the link group. In this spirit, given a link $L$, we can define a \emph{torsion quotient} of $L$ as quotient of its link group by the normal closure of some powers of meridians. Combining the results of \cite{Necklace} with Theorem \ref{theo:embedding_torsion_quotients_intro} and Theorem \ref{theo:generalized_j_groups_torsion_quotients_intro}, we obtain one of the main results of this article:

\begin{theorem}\label{theo:everyone_friend_intro}
The following families of groups coincide, up to isomorphism preserving the set of generators
\begin{itemize}
\item Finitely generated generalized $J$-groups (generators: reflections)
\item Reduced $J$-groups (generators: reflections)
\item Torsion quotients of $J$-braid groups (generators: reflections)
\item Torsion quotients of torus necklaces (generators: nontrivial powers of images of meridians).
\end{itemize}
\end{theorem}

Combining Corollary \ref{cor:finite_generalized_j_groups_intro} and Theorem \ref{theo:everyone_friend_intro}, with very little additional work we are able to give a new point of view for complex reflection groups of rank two:

\begin{corollary}\label{finite_crg_seifert_quotient_intro}
The complex reflection groups of rank two are precisely the finite torsion quotients of Seifert links.
\end{corollary}

The last main result of this article is the classification of torsion quotients of Seifert links up to isomorphisms which preserves the generators. The following result is slightly more general than Corollary \ref{cor:iso_quot_j-braid_intro}, although the additional cases are easily dealt with.
\begin{theorem}[Theorem \ref{theo:class_torsion_quotient_seifert}]
If two torsion quotients of Seifert links are isomorphic in a way which preserves the nontrivial powers of images of meridians, then the underlying links are isotopic.\end{theorem} 

\subsection*{Acknowledgement}
We thank Axel Péneau for enlightening discussions about \\Fuchsian groups. We also thank Jon McCammond for emphasizing the geometric connection between $J$-groups and triangle groups.\\
The first author is supported by the projects PID2022-138719NA-I00 funded by MCIN/AEI/10.13039/501100011033 and by FEDER, UE.\\
The second author is supported by the ANR
project CORTIPOM (ANR-21-CE40-0019).


\section{Preliminaries on Coxeter groups}\label{sec:prelim_coxeter}
\subsection{Elementary definitions and results}
Recall that a \emph{Coxeter graph} is an edge-weighted simple graph $\Gamma=(\Gamma,m)$, where $m$ takes values in $\N_{\geqslant 2}$. Given a Coxeter graph $\Gamma$, we will write $S$ for the set of vertices of $\Gamma$.  The \emph{Coxeter group}\index{Coxeter groups} $W=W[\Gamma]$ attached to a Coxeter graph $\Gamma$ is defined by the following group presentation:

\begin{equation}\label{CoxPresGraph}
    \left\langle S\,|\, s^2=1 \, \forall s\in S ; \, (st)^{m_{\{s,t\}}}=1\, \mathrm{for \, all} \,  \{s,t\}\in E(\Gamma)\right\rangle.
\end{equation}
The cardinality of $S$ is the \emph{rank} of $W[\Gamma]$. If $S$ cannot be decomposed as a disjoint union $S_1\sqcup S_2$ such that every vertex of $S_1$ is connected to every vertex of $S_2$ by an edge of length two, we say that $W[\Gamma]$ is \emph{irreducible}.

For a subset $S'\subset S$, we denote by $W_{S'}$ the subgroup of $W$ generated by $S'$. It is classical (see \cite[Section IV.1.8]{BourbakiLie4-6}) that $W_{S'}$ is the Coxeter group attached to the full subgraph of $\Gamma$ whose vertex set is $S'$. The subgroup $W_{S'}$ is called a \emph{standard parabolic subgroup of $W$}. A \emph{parabolic subgroup of $W$} is the conjugate of a standard parabolic subgroup

Any Coxeter group $W$ admits a natural morphism $W\to \Z/2\Z$, sending every generator $s\in S$ to $1$. The kernel of this morphism is the \emph{alternating subgroup} $W^+$ of the Coxeter group $W$.

Let $W$ be a Coxeter group. We give a short list of results we will use in the sequel:
\begin{itemize}
    \item The order of an element $w\in W$ is finite if and only if $w$ belongs to a finite parabolic subgroup of $W$ \cite[Exercise 2, Section 4, Chapter 5]{BourbakiLie4-6}.
    \item If the rank of $W$ is at least $3$, then $Z(W^+)\subset Z(W)$. In particular, if $W$ is infinite and irreducible, then $Z(W^+)$ is trivial \cite[Proposition 3.1]{GobetToric}. \item If $W$ is irreducible, infinite and non-affine, each finite index subgroup of $W$ has trivial center \cite[Proposition 6.4]{ParisIrrCox}.
\end{itemize}

\subsection{Polygonal Coxeter groups}\label{sec:polycox}
One part of the approach of Gobet in \cite{GobetToric} for studying torsion quotients of circular groups with coprime parameters consists in relating such a quotient to the alternating subgroup of a triangle Coxeter Group. We will use similar techniques when dealing with arbitrary circular groups, but the Coxeter groups involved form a larger family, which we call \emph{polygonal Coxeter groups}.

In this section, we fix an integer $d\geqslant 3$, along with a $d$-tuple $K:=(k_1,\ldots,k_d)$ in $(\N_{\geqslant 2}\cup \{\infty\})^d$. 

\begin{definition}[\textbf{Polygonal Coxeter group}]
We define a Coxeter graph $\Gamma_K$ as follows:
\begin{itemize}
    \item The set of vertices is a set $\{v_1,\ldots,v_d\}$ of cardinality $d$ (the indices are seen modulo $d$).
    \item All edges have form $\{v_i,v_{i+1}\}$. Moreover, the pair $\{v_i,v_{i+1}\}$ is actually an edge if and only if $k_i<\infty$, in which case its weight is $k_i$.
\end{itemize}
The group $W_K$ is the Coxeter group $W[\Gamma_K]$, which we call the \emph{polygonal Coxeter group} attached to $K$.
\end{definition}
By construction, the Coxeter presentation of $W_K$ is
\begin{equation}\label{PresCoxeterW}
\left\langle
\begin{array}{c}
   s_1,s_2,\dots,s_d
\end{array}
\;\middle|\;
\begin{array}{c}
  s_i^2=1 \, \text{for all } \, i=1,\dots,d\\
  (s_is_{i+1})^{k_i}=1\, \text{for all } \, i=1,\dots,d \,\, \text{such that } \, k_i<\infty
\end{array}\right\rangle,
\end{equation}
where the indices are seen modulo $d$.

\begin{remark}\label{AffinePolygonalCox}
Using the classification of finite Coxeter groups (see \cite{CoxeterClassification2}),  the group $W_K$ is finite if and only if $K\in\{(2,3,3),(2,3,4),(2,3,5)\}\cup\{(2,2,l)\}_{l\geqslant 2}$ up to permutation. The only cases in which $W_K$ can be infinite and not irreducible are $K=(\infty,2,2)$ $K=(2,\infty,2)$, $K=(2,2,\infty)$ and $K=(2,2,2,2)$. The only cases were $W_K$ is affine are $K=(3,3,3),(2,3,6),(2,4,4)$ or $(2,2,2,2)$ (again up to permutation).
\end{remark}

In \cite{BourbakiLie4-6}, a presentation for the alternating group of a given Coxeter group is given. Specialising this presentation to $W_K$ yields the following result:

\begin{proposition}\cite[Exercise 9, Section 1, Chapter 4]{BourbakiLie4-6}\label{PresBourbaki}
The group $W_K^+$ admits the following group presentation:
\begin{equation}\label{PresBourbakiW}
    \left\langle
\begin{array}{c}
   b_1,\dots,b_{d-1}
\end{array}
\;\middle|\;
\begin{array}{c}
  b_1^{k_d}=b_{d-1}^{k_{d-1}}=1 \\
  (b_ib_{i+1}^{-1})^{k_i}=1\, \text{for all } \, i\in \intv{1,d-2} 
\end{array}\right\rangle,
\end{equation}
where the relation corresponding to $i\in \intv{1,d-2}$ is empty if $k_i=\infty$.
\end{proposition}

Using this, we obtain another presentation of $W_K^+$:

\begin{corollary}\label{AlternatingK}
The group $W_K^+$ admits the group presentation 
\begin{equation}\label{PresAlternatingW}
\left\langle
\begin{array}{c}
   a_1,a_2,\dots,a_d
\end{array}
\;\middle|\;
\begin{array}{c}
  a_i^{k_i}=1 \,  \text{for } \, i\in \intv{1,d} \text{ with } \, k_i<\infty\\
  a_1a_2\cdots a_d=1
\end{array}\right\rangle.
\end{equation}
\end{corollary}

\begin{proof}
Let $H$ be the group defined by Presentation \eqref{PresAlternatingW}. 
We define an isomorphism between the groups $H$ and $W_K^+$.\\
First, we define a morphism $\varphi$ from $H$ to $W_K^+$ by setting
\[\forall i\in \intv{1,d},a_i\mapsto\begin{cases}b_ib_{i+1}^{-1} & \text{if }i\in \intv{1,d-2},\\ b_{d-1}&\text{if }i=d-1,\\b_1^{-1}&\text{if }i=d.\end{cases}\]
Indeed, we have $\varphi(a_i)^{k_i}=1$ for all $i\in \intv{1,d}$ and
\[1=(b_1b_2^{-1})(b_2b_3^{-1})\cdots(b_{d-2}b_{d-1}^{-1})b_{d-1}b_1^{-1}=\varphi(a_1)\varphi(a_2)\cdots\varphi(a_d).\]
The inverse of this morphism is given by setting $b_1\mapsto a_d^{-1}$ and $b_i\mapsto (a_da_1\cdots a_{i-1})^{-1}$ for $i\in \intv{2,d}$. This concludes the proof.
\end{proof}

We finish this section by establishing some group theoretic results on polygonal Coxeter groups. Namely, we describe the centers and torsion elements of infinite polygonal Coxeter groups, and, we give the classification of polygonal Coxeter group up to group isomorphism.

\begin{corollary}\label{CenterAlternatingW}
The center of $W_K^+$ is trivial. 
\end{corollary}
\begin{proof}
If $W_K$ is finite, this is proven in \cite[Proof of Theorem 3.3]{GobetToric}. If $W_K$ is infinite and irreducible, then the result is a direct application of \cite[Proposition 3.1]{GobetToric}.\\
The only cases where $W_K$ is infinite and reducible are $K=(\infty,2,2)$ and $K=(2,2,2,2)$. For $K=(\infty,2,2)$, Presentation $\eqref{PresAlternatingW}$ reads 
\begin{equation*}
    \left \langle a_1,a_2,a_3\,|\, a_2^2=a_3^2=a_1a_2a_3=1\right\rangle,
\end{equation*}
thus $W_K^+$ is isomorphic to $\Z/2\Z*\Z/2\Z$, whose center is trivial.\\
For $K=(2,2,2,2)$, the group $W_K$ is $(\Z/2\Z*\Z/2\Z)\times (\Z/2\Z*\Z/2\Z)$, whose center is trivial. Applying \cite[Proposition 3.1]{GobetToric} again yields that $Z(W_K^+)$ is trivial in this case. This concludes the proof.
\end{proof}

\begin{proposition}\label{PermutePolygonalCoxeter}
Let $K\in (\N_{\geqslant 2}\cup\{\infty\})^d$ and $\sigma\in S_d$. Write $\sigma(K)$ for the tuple $(k_{\sigma(1)},\dots,k_{\sigma(d)})$. We have $W_{K}^+\cong W_{\sigma(K)}^+$.
\end{proposition}

\begin{proof}
In this proof, we use Presentation \eqref{PresAlternatingW}. It is enough to show the result for the permutations $(1,2)$ and $(1,2,\dots,d)$, since these two permutations generate $S_d$. For $(1,2,\dots,d)$, the result is immediate. Let then $\sigma=(1,2)$.\\
For the remainder of the proof, we write $a_i$ for the generators of $W_K^+$ and $b_i$ for the generators of $W_{\sigma(K)}^+$.
We obtain a well-defined morphism from $W_K^+$ to $W_{\sigma(K)}^+$ by setting 
\[a_i\mapsto \begin{cases} b_2^{-1}&\text{if }i=1,\\b_1^{-1}&\text{if }i=2,\\ (b_3\cdots b_{i-1})b_i^{-1}(b_3\cdots b_{i-1})^{-1}&\text{if }i\in\intv{3,d}.\end{cases}\]
Exchanging the role of $b_i$ and $a_i$ yields a morphism from $W_{\sigma(K)}$ to $W_{\sigma^2(K)}=W_K$, which is the inverse of the first morphism. This concludes the proof. 
\end{proof}

\begin{proposition}\label{TorsionPolygonalCox}
If $W_K$ is infinite, then an element $w$ in $W_K^+$ has finite order if and only it is a conjugate to a power of some generator $a_i$ in Presentation \eqref{PresAlternatingW} such that $k_i$ is finite.
\end{proposition}

\begin{proof}
We mentioned in Section \ref{sec:prelim_coxeter} that the torsion elements of $W_K$ are exactly the conjugates of elements lying in finite standard parabolic subgroups. The finite standard parabolic subgroups of $W_K$ are those generated by two adjacent vertices in $\Gamma_K$ with finite length. 
Now, recall that given a Coxeter group $W$ and a parabolic subgroup $W'<W$, the alternating subgroup of $W'$ can be described as $W'\cap W^+$. Since the elements $\{a_1,\dots,a_d\}$ are the generators of the alternating subgroups of the standard parabolic subgroups generated by two adjacent vertices, this concludes the proof. 
\end{proof}

\begin{corollary}\label{MaximalFinitePolygonalCox}
Assume that $W_K$ is infinite, and let $p$ be the number of coordinates of $K$ equal to $\infty$. There are $d-p$ conjugacy classes of maximal finite subgroups of $W_K^+$. The finite groups in these $d-p$ classes are isomorphic to the $\Z/k_i\Z$ such that $k_i$ is finite, in which case they are conjugates to $\langle a_i\rangle$.
\end{corollary}

\begin{proof}
The proof of \cite[Proposition 4.6]{GobetToric} only uses that every finite standard parabolic subgroup of the considered Coxeter groups has rank at most two, which is true in our case. Thus, the aforementioned proof carries out in our setting. 
\end{proof}

\begin{lemma}\label{lem:no_roots}
Assume that $W_K^+$ is infinite. For $i\in \intv{1,d}$, the group $\langle a_i \rangle$ is a maximal cyclic group of $W_K^+$.\end{lemma}
\begin{proof}
By Proposition \ref{PermutePolygonalCoxeter}, we can assume that $i=1$.

First, assume that $k_1$ is finite. In this case, we know by Corollary \ref{MaximalFinitePolygonalCox} that $\langle a_1\rangle$ is a maximal finite subgroup of $W_K^+$. Since $a_1$ has finite order, any cyclic group containing $a_i$ must be finite and we have the result.

Now, assume that $k_1$ is infinite. In this case, we can replace $a_1$ by $(a_2\cdots a_d)^{-1}$ in Presentation \eqref{PresAlternatingW} to obtain 
$W_K^+\simeq \Z/k_2\Z*\cdots*\Z/k_d\Z$. The existence of a normal form in free products of cyclic group then implies that $(a_2\cdots a_d)$ does not admit (nontrivial) roots. Again, this gives the desired result.
\end{proof}

\begin{proposition}\label{prop:centralizer_rot_polyg_coxeter_group}
Assume that $W_K^+$ is infinite. For $i\in \intv{1,d}$ and $n\in \Z^*$, if $a_i^n$ is nontrivial, then $C_{W_K^+}(a_i^n)=\langle a_i\rangle$.
\end{proposition}
\begin{proof}
If $W_K^+$ is affine, then $a_i^n$ is an affine rotation, with exactly one fixed point $P$. If $g$ is a positive isometry of the Euclidean plane which centralizes $a_i^n$, then $g.P=P$ and $g$ is a rotation of center $P$. The centralizer of $a_i^n$ in the group of positive isometries is then isomorphic to $\R$. Since $W_K^+$ is discrete, the centralizer of $a_i^n$ in $W_K^+$ is a discrete subgroup of $\R$: it is a cyclic group. Since $\langle a_i\rangle \subset C_{W_K^+}(a_i^n)$, we have the result by Lemma \ref{lem:no_roots}.

If $W_K^+$ is not affine, hence a Fuchsian group, we can apply \cite[Theorem 2.3.2 and 2.3.5]{katok}, which gives that the centralizer of $a_i^k$ in $W_K^+$ is a cyclic group. Again, since $\langle a_i\rangle \subset C_{W_K^+}(a_i^n)$, we have the result by Lemma \ref{lem:no_roots}. 
\end{proof}

\begin{proposition}\label{ClassificationPolygonalCox}
Let $K_1\in (\N_{\geqslant 2}\cup\{\infty\})^{d_1}$, $K_2\in (\N_{\geqslant 2}\cup\{\infty\})^{d_2}$. We have $W_{K_1}^+\cong W_{K_2}^+$ if and only if $d_1=d_2$ and $K_1=K_2$ up to permutation.
\end{proposition}

\begin{proof}
The if statement is exactly Proposition \ref{PermutePolygonalCoxeter} and we only have to show the only if part.

If the groups are finite, this can be seen for example in \cite[Table 2]{GobetToric}. Assume then that the groups are infinite.\\
If $W_{K_1}^+\cong W_{K_2}^+$, the quotients of these groups by the normal closure of their torsion elements are isomorphic as well. Using Proposition \ref{TorsionPolygonalCox}, these quotients are isomorphic to a free group of rank one less than the number of infinite indices in $K_1$ and $K_2$, hence these numbers are equal, say $p$.\\
Moreover, Corollary \ref{MaximalFinitePolygonalCox} identifies the multisets $\{\!\{k_{1,1},\dots,k_{1,d_1-p}\}\!\}$ and \\ $\{\!\{k_{2,1},\dots,k_{2,d_2-p}\}\!\}$. Thus, $d_1-p=d_2-p$ so that $d_1=d_2$. We conclude that up to permutation, the finite parameters of $K_1$ and $K_2$ are the same, which concludes the proof.
\end{proof}

\subsection{Subgroups of alternating polygonal Coxeter groups}

In the following sections, we will be able to study torsion quotients by relating them to particular subgroups of some group of the form $W_K^+$.

In this section, we fix a positive integer $d$, along with a $d$-tuple $K:=(k_1,\ldots,k_d)$ in $(\N_{\geqslant 2}\cup \{\infty\})^d$. We also fix another tuple $K'=(k_1',\ldots,k_d')$ in $(\N_{\geqslant 1}\cup \{\infty\})^d$ such that $k_i'$ divides $k_i$ if $k_i<\infty$, and such that $k_i'=\infty$ implies $k_i=\infty$.

We consider the presentation \eqref{PresAlternatingW} of the group $W_K^+$, and we define $W_K^+(K')$ to be the normal closure of the elements $a_i^{k_i'}$ in $W_K^+$. By construction, the quotient $W_K^+/W_K^+(K')$ is isomorphic to $W_{K'}^+$, or more precisely to $W_{K''}^+$, where $K''$ is obtained from $K'$ by removing the entries equal to $1$.

The main purpose of this section is to prove the following result:

\begin{proposition} \label{prop:infinitely_generated_in_coxeter}
If the group $W_K^+(K')$ is nontrivial, then it is finitely generated if and only if it has finite index in $W_K^+$.

\end{proposition}
An easy consequence of Schreier's lemma is that a finite index subgroup in a finitely generated group is again finitely generated. In particular, it is sufficient to show that if $W_K^+(K')$ is nontrivial and has infinite index in $W_K^+$, then it is not finitely generated. In order for $W_K^+(K')$ to have infinite index in $W_K^+$, it is necessary for $W_K^+$ to be infinite.

Let $X$ be either the Euclidean plane or the hyperbolic plane. We can consider a convex polygon $P$ in $X$ whose dihedral angles are $\pi/k_i$ for $i\in \intv{1,d}$ (or $0$ when $k_i$ is infinite). By \cite[Theorem 7.1.3 and 7.1.4]{ratcliffe}, the group generated by the reflections of $X$ relative to the sides of $P$ is isomorphic to $W_K$. Moreover, this group admits $P$ as a fundamental domain for its action on $X$.

If $P$ is an affine polygon, then $W_K$ is an affine Coxeter group. By Remark \ref{AffinePolygonalCox}, we have $K=(3,3,3),(2,3,6),(2,4,4)$ or $(2,2,2,2)$. We will deal with these cases at the end of the proof. Otherwise, $P$ is an hyperbolic polygon, and $W_K$ is a discrete group of isometries of the hyperbolic plane. The subgroup $W_K^+$ is then a discrete group of orientation preserving isometries of the hyperbolic plane, i.e. a \emph{Fuchsian group}. 

The theory of Fuchsian group is well known, we use results from \cite{ratcliffe} and \cite{katok} regarding these groups. 

The first thing we need to do in this case is to prove that if $W_K^+(K')$ is nontrivial, then it is infinite. If $W_K^+(K')$ is nontrivial, then it contains some nontrivial element $w:=a_i^{k_i'}$ with $k_i'<k_i$. If this element has infinite order, then of course $W_K^+(K')$ is infinite, so we can assume that $w$ has finite order. By \cite[Theorem 2.3.3]{katok}, the centralizer $H$ of $w$ in $W_K^+$ consists of the elements of $W_K^+$ whose actions on the hyperbolic plane fix the same points $w$. Since $H$ is discrete as a subgroup of $W_K^+$, it is a then a cyclic group by \cite[Theorem 2.3.5]{katok}. Moreover, $H$ is finite since it contains the finite order element $w$. Since $W_K^+$ is infinite, the index of $H$ in $W_K^+$ is infinite since $H$ is finite. We conclude by the orbit-stabilizer theorem that the conjugacy class of $w$ in $W_K^+$, which is included in $W_K^+(K')$, is infinite so that $W_K^+(K')$ is infinite.

Now, let us show that $W_K^+(K')$ is finitely generated if and only if it has finite index in $W_K^+$. By \cite[Theorem 3.5.4 and 4.6.1]{katok}, a Fuchsian group is finitely generated if and only if it is \emph{geometrically finite}. By \cite[Theorem 4.1.1 and 4.5.1 and 4.5.2]{katok}, the area of the fundamental domain of a Fuchsian group is finite if and only if the group is geometrically finite and of the \emph{first kind}.

Now, the polygon $P$ has finite hyperbolic area by \cite[Exercise 3.5.13]{ratcliffe}, thus $W_K^+$ is a geometrically finite Fuchsian group of the first kind. By \cite[Theorem 12.2.1]{ratcliffe}, $W_K^+$ is nonelementary since it is of the first kind, and thus $W_K^+(K')$ is also a Fuchsian group of the first kind by \cite[Theorem 12.2.14]{ratcliffe} since it is infinite and normal in $W_K^+$.

Since $W_K^+(K')$ is of the first kind, we obtain that it is finitely generated if and only if the area of a fundamental fomain for $W_K^+(K')$ is finite. However, by \cite[Theorem 6.7.3]{ratcliffe}, the area of such a domain is the product of the area of $P$ with the index of $W_K^+(K')$ in $W_K^+$, which terminates the proof in this case. 

It remains to prove the result for the case where $W_K$ is affine. First assume that $K=(2,2,2,2)$. If $W_K^+(K')$ is nontrivial, then $W_K^+/W_K^+(K')=W_{K'}^+$ is a quotient of $W_{(2,2,2)}^+$ which is finite. The group $W_K^+(K')$ then has finite index in $W_K^+$ and there is nothing to show.

Otherwise, we have $d=3$ and $K=(k_1,k_2,k_3)$. If $k_i=1$ for some $i$, then $W_K^+/W_K^+(K')=W_K'^+$ is a quotient of a group of the form $W_{(p,q)}^+$, which is always finite. We can then assume that no $k_i'$ is equal to one. In particular if $k_i$ is a prime, then we can assume that $k_i'=k_i$. We can also assume that $K'\neq K$ since otherwise $W_K^+(K')$ is trivial and there is nothing to show.
\begin{itemize}
\item If $K=(3,3,3)$, then there is no remaining possibility for $K'$ and there is nothing to show.
\item If $K=(2,3,6)$, then $K'=(2,3,2)$ or $(2,3,3)$. In both cases $W_K^+/W_K^+(K')$ is finite by Remark \ref{AffinePolygonalCox} and there is nothing to show.
\item If $K=(2,4,4)$, then $K'=(2,2,2),(2,2,4)$ or $(2,4,2)$. In all cases the quotient $W_K^+/W_K^+(K')$ is finite by Remark \ref{AffinePolygonalCox} and there is nothing to show.
\end{itemize}

Now that the proof of Proposition \ref{prop:infinitely_generated_in_coxeter} is completed, we finish this section with a corollary regarding the center of the group $W_K^+(K')$.

\begin{corollary}\label{cor:center_wkk'}
Assume that $W_K^+$ is infinite. If $W_K^+(K')$ is nontrivial, then its center is trivial.
\end{corollary}

\begin{proof}
If $W_K^+(K')$ is nontrivial, then it contains some nontrivial element $a_i^{k_i'}$. Up to permuting the $a_i$ (using Proposition \ref{ClassificationPolygonalCox}), we can assume that $a_1^{k_1'}$ is nontrivial.

Now, the center of $W_K^+(K')$ is included in the centralizer of $a_1^{k_1'}$ in $W_K^+$, which is equal to $\langle a_1\rangle$ by Proposition \ref{prop:centralizer_rot_polyg_coxeter_group}. Assume that $Z(W_K^+(K'))$ contains a nontrivial element $z$. We have $z= a_i^n$ for some integer $n$. Any $x\in W_K^+(K')$ lies in the centralizer of $z$ in $W_K^+$, and thus $x\in \langle a_i \rangle$. We then have $W_K^+(K')\subset \langle a_i\rangle$.

Now, as during the proof of Proposition \ref{prop:infinitely_generated_in_coxeter}, the group $W_K^+$ is a discrete isometry group for $X$ either a euclidean or hyperbolic plane. Moreover, a fundamental domain for the action of $W_K^+$ on $X$ is a convex polygon $P$.

The group $W_K$ is generated by the reflections along the sides of $P$, and the elements $a_1,\ldots,a_d$ are rotations around the vertices of $P$. Let $v_1$ (resp. $v_2$) be the vertex of $P$ which is a center for the rotation $a_1$ (resp. $a_2)$. The element $a_1^{k_1'}$ is also a rotation centered at $v_1$. Since $v_1\neq v_2$, the point $v_1$ is not a fixed point of $a_2$, and thus $a_2a_1^{k_1'}a_2^{-1}$ is a rotation around a point different from $v_1$. In particular, $a_2a_1^{k_1'}a_2^{-1}$ is not a power of $a_1$. Since $a_2a_1^{k_1'}a_2^{-1}\in W_K^+(K')$, this is a contradiction.\end{proof}

\section{Generalized $J$-groups}\label{sec:generalized_j-groups}
The family of $J$-groups was first defined in \cite{AA}, where the authors proved in particular that the finite groups in this family are precisely the complex reflection groups of rank 2. In this section, we give a generalization of the definition of $J$-groups and we prove some group-theoretic results in this new setting.

For the remainder of the section, fix a positive integer $p$, along with a $p$-tuple $K=(k_1,\ldots,k_p)$ in $(\N_{\geqslant 2}\cup \{\infty\})^p$. We also fix another tuple $K'=(k_1',\ldots,k_p')$ in $(\N_{\geqslant 1})^p$ such that $k_i'$ divides $k_i$ for each $i$ (by convention, we say that every nonzero integer divides $\infty$).

\subsection{Definition, reflections}

\begin{definition}[\textbf{Generalized $J$-groups}]\label{def:generalized_j-groups}
The group $J(K)$ is defined by the group presentation 
\begin{equation*}
    \left\langle s_1,\dots,s_p\, |\, s_1s_2\cdots s_p=s_2\cdots s_ps_1=\dots=s_ps_1\cdots s_{p-1},\, s_i^{k_i}=1 \, \text{if $k_i<\infty$}\right\rangle.
\end{equation*}
It is called a \emph{generalized parent $J$-group}.\\ The \emph{generalized $J$-group} $J\binom{K}{K'}$ is defined as the normal closure of $\{s_1^{k_1'},\ldots,s_p^{k_p'}\}$ in $J(K)$.
\end{definition}

By extension, we will often consider groups of the form $J(K)$ with elements of $K$ possibly equal to $1$. Considering the above presentation, this is equivalent to considering $J(\hat{K})$, where $\hat{K}$ is obtained from $K$ by removing the entries equal to $1$.

Since $k_i'$ divides $k_i$ for all $i$, we have a natural surjective morphism $J(K)\to J(K')$. The kernel of this morphism is precisely the generalized $J$-group $J\binom{K}{K'}$. 

\begin{remark}[The case $k_i'=\infty$]
The definition of generalized $J$-group also makes sense if we allow $K'$ to contain infinity elements. Assume for instance that $k_p'=k_p=\infty$, and let $z:=s_1\cdots s_p$. By looking at the presentation of $J(K)$, we can show that $J(K)$ decomposes as a direct product $\langle s_1,\ldots,s_{p-1}\rangle \times \langle z\rangle$. Moreover, the group $\langle s_1,\ldots,s_{p-1}\rangle$ is a free product of the cyclic groups $\langle s_i\rangle$. The group $J\binom{K}{K'}$ is then the normal closure of $s_1^{k_1'},\ldots,s_{p-1}^{k_{p-1}'}$ in $\langle s_1,\ldots,s_{p-1}\rangle$. We obtain that $J\binom{K}{K'}$ is also a (possibly infinitely generated) free product of cyclic groups. We chose not to consider this case here, since on the one hand free products of cyclic groups are well-studied, and on the other hand it would require to specify particular cases for several theorems below.
\end{remark}

Just as in the case of $J$-groups, we have a natural notion of abstract reflections in generalized $J$-groups. This definition imitates the definition originally given in \cite{GobetToric} in the case of toric reflection groups, which are particular (classical) $J$-groups.

\begin{definition}[\textbf{Reflections}]
The set of \emph{reflections} of the generalized parent $J$-group $J(K)$ is defined as the set of all conjugates of nontrivial powers of $s_1,\ldots,s_p$ in $J(K)$. It is denoted by $R(J(K))$. \\ The set of \emph{reflections} of the generalized $J$-group $J\binom{K}{K'}$ is defined as the set of all conjugates of nontrivial powers of $s_1^{k_1'},\ldots,s_p^{k_p'}$ in $J(K)$. It is denoted by $R(J\binom{K}{K'})$.
\end{definition}
Notice that $R(J\binom{K}{K'})$ is not the set of conjugates of the $s_i^{k_i'}$ in $J\binom{K}{K'}$, but rather in the generalized parent $J$-group $J(K)$. In particular, while computing the conjugacy classes of reflections in a generalized parent $J$-group is rather easy, it is more difficult in an arbitrary generalized $J$-group (see Proposition \ref{prop:conj_reflections_j-groupe}).

Another reasonable definition for the set of reflections of a generalized $J$-group would be the intersection with said $J$-group of the set of reflections of its parent generalized $J$-group. The following two lemmas show that the two definitions coincide.

\begin{lemma}\label{lem:reflexion_belonging_to_j_group}
Let $i\in \intv{1,p}$, and let $n\in \N$. The reflection $s_i^n$ belongs to $J\binom{K}{K'}$ if and only if $n$ is a multiple of $k_i'$.
\end{lemma}
\begin{proof}
The if part is immediate by definition of $J\binom{K}{K'}$. Conversely,  let $\bbar{s_i}$ be the image of $s_i$ in $J(K')$. Since $J\binom{K}{K'}$ is the kernel of the projection map $J(K)\to J(K')$, we have $s_i^n\in J\binom{K}{K'}$ if and only if $\bbar{s_i}^n$ is trivial. It is then sufficient to show that the order of $\bbar{s_i}$ is $k_i'$. By construction, we know that this order divides $k_i'$. Conversely, we can compute the abelianization of $J(K')$ using its defining presentation. We obtain that this abelianization is a direct product of $p$ cyclic groups $\Z/k_i'\Z$. In particular, the order of the image of $\bbar{s_i}$ in this abelianization in $k_i'$, whence the order of $\bbar{s_i}$ is precisely $k_i'$, which finishes the proof.
\end{proof}

\begin{lemma}
We have $R(J\binom{K}{K'})=R(J(K))\cap J\binom{K}{K'}$.
\end{lemma}
\begin{proof}
Since reflections in $J\binom{K}{K'}$ are in particular reflections in $J(K)$, we have the direct inclusion. Conversely, let $r\in R(J(K))$ be a reflection which belongs to $J\binom{K}{K'}$. By definition, there is an element $g\in J(K)$ which conjugates $r$ to some nontrivial reflection $s_i^m$ of $J(K)$. Since $J\binom{K}{K'}$ is normal in $J(K)$, we have $s_i^m\in J\binom{K}{K'}$. By Lemma \ref{lem:reflexion_belonging_to_j_group}, we have $s_i^m\in R(J\binom{K}{K'})$ and thus $r\in R(J\binom{K}{K'})$ as we wanted to show.
\end{proof}

Since a generalized $J$-group is defined as the normal closure of a set of reflections in its generalized parent $J$-group, it is generated by reflections. However, it is not obvious at this stage that it can be generated by a finite number of reflections. In fact this will only be true for a particular family of generalized $J$-groups which is introduced below (Definition \ref{def:finite-type}).

We now naturally extend the notion of reflection (iso)morphisms introduced in \cite{GobetToric}:

\begin{definition}[\textbf{Reflection isomorphism}]
Let $J_1$, $J_2$ be two generalized $J$-groups. A group morphism $J_1\xto\varphi J_2$ is a \textit{reflection morphism} if $\varphi(R(J_1))\subset R(J_2)\cup\{1_{J_2}\}$. Moreover, the groups $J_1$ and $J_2$ are said to be \textit{reflection isomorphic} if there exists a group isomorphism $J_1\to J_2$ which restricts to a bijection between $R(J_1)$ and $R(J_2)$. In this case, we write $$H_1\cong_{\mathrm{Ref}} H_2.$$
\end{definition}

An important fact is that the ordering of the tuples $K$ and $K'$ does not impact the reflection isomorphism type of the generalized $J$-group $J\binom{K}{K'}$:

\begin{proposition}[\textbf{Permutation of parameters}]\label{prop:permutation_of_parameters_generalized_j_groups}
Let $\sigma \in S_p$ and write $\sigma(K)=(k_{\sigma(1)},k_{\sigma(2)},\dots,k_{\sigma(p)})$ (and similarly for $\sigma(K')$). The groups $J\binom K{K'}$ and $J\binom{\sigma(K)}{\sigma(K')}$ are reflection isomorphic.
\end{proposition}
\begin{proof}
Since generalized $J$-groups are normal closures in generalized parent $J$-groups, it is enough to prove the result for generalized parent $J$-groups.\\
If $p=1$, there is nothing to show. If $p=2$, the result is obtained by considering the automorphism swapping $s_1$ and $s_2$.\\
We assume from now on that $p\geqslant 3$. Since $S_p$ is generated by the transposition $s:=(1~2)$ and the cycle $c:=(1~2~\cdots~p)$, it is enough to show the result for these two permutations. Writing $t_1,\dots,t_p$ for the generators of $J(K)$, we obtain a well-defined morphism from $J(K)$ to $J(s(K))$ by setting
\[s_i\mapsto \begin{cases} t_2^{-1}&\text{if }i=1,\\t_1^{-1}&\text{if }i=2,\\(t_3\cdots t_{i-1})t_i^{-1}(t_3\cdots t_{i-1})^{-1}&\text{if }i\in \intv{3,p}.\end{cases}\]
Exchanging the roles of $s_i$ and $t_i$ yields a morphism from $J(s(K))$ to $J(s^2(K))=J(K)$, which is the inverse of the first morphism.
Moreover, $\alpha(s_1)$ is conjugate to $t_2^{-1}$, $\alpha(s_2)$ is conjugate to $t_1^{-1}$ and $\alpha(s_i)$ is conjugate to $t_i^{-1}$ for $i\in \intv{3,p}$, so that $\alpha$ is a reflection isomorphism.

Now, by definition of $J(c(K))$, a reflection isomorphism $\gamma:J(K)\to J(c(K))$ is obtained by sending $s_1$ to $t_{i+1}$ for $i\in \intv{1,p-1}$ and by sending $s_p$ to $t_1$. This concludes the proof.
\end{proof}

\subsection{Determination of finite generalized $J$-groups}
One of the main results of \cite{AA} is the determination of finite $J$-groups. More specifically, the authors show that the family of finite $J$-groups coincide with the family of complex reflection groups of rank $2$. In this section, we generalize this result to arbitrary generalized $J$-groups.

If $p\leqslant 2$, then $J(K)$ is finite if and only if the elements of $K$ are all finite. In particular, since the elements of $K'$ are finite, the group $J(K')$ is finite and $J\binom{K}{K'}$ has finite index in $J(K)$. In particular, $J\binom{K}{K'}$ is finite if and only if its generalized parent $J$-group is also finite.

Obtaining a similar result in the case $p\geqslant 3$ is possible but more intricate. Just as $J$-reflection groups are related to triangle Coxeter groups (see for instance \cite[Theorem 1.5]{VCRG}), generalized $J$-groups are related to polygonal Coxeter groups. More precisely, we have the following result:

\begin{proposition}[\textbf{Center of parent generalized $J$-groups}]\label{prop:center_parent_generalized_j_groups}
Assume $p\geqslant 3$. The center of $J(K)$ is generated by the product $s_1s_2\cdots s_p$. Moreover, the quotient $J(K)/Z(J(K))$ is isomorphic to the alternating polygonal Coxeter group $W_K^+$.
\end{proposition}
\begin{proof}
First, $s_1s_2\cdots s_p$ is central in $J(K)$ (it is invariant by conjugation by any of the $s_i$'s). Using Corollary \ref{AlternatingK}, the correspondence $s_i\mapsto a_i$ for $i\in \intv{1,p}$ induces an isomorphism between $J(K)/\langle s_1\cdots s_p \rangle$ and $W_K^+$. Since the latter group has trivial center by Corollary \ref{CenterAlternatingW}. We deduce that $Z(J(K))\subset \langle s_1s_2\cdots s_p\rangle \subset Z(J(K))$, which concludes the proof. 
\end{proof}

This result allows us to determine which generalized parent $J$-groups are finite, by relating them to alternating polygonal Coxeter groups.

\begin{lemma}\label{lem:finite_index_J_or_coxeter}
Assume $p\geqslant 3$. The group $J(K)$ is finite if and only if $W^+_K$ is finite. This is equivalent to having (up to permutation) $K=(2,2,l)$ for $l\geqslant 2$ or $K\in \{(2,3,3),(2,3,4),(2,3,5)\}$.\end{lemma}
\begin{proof}
By Corollary \ref{AlternatingK}, there is a natural quotient $J(K)\to W_K^+$ sending $s_i$ to $a_i$. Therefore, if $W_K^+$ is infinite, then so is $J(K)$. As said in Remark \ref{AffinePolygonalCox}, the group $W_K$ (hence $W_K^+$) is finite if and only if $K$ belongs up to permutation to the set $\{(2,2,l),(2,3,3),(2,3,4),(2,3,5)\}$ (with $l\geqslant 2$. The corresponding groups $J(K)$ are known to be finite \cite[Theorem 1.2]{AA}, this finishes the proof.
\end{proof}

In order to obtain a complete description of finite generalized $J$-groups, we introduce the definition of finite-type generalized $J$-group.

\begin{definition}[\textbf{Finite-type}]\label{def:finite-type}
A nontrivial generalized $J$-group $J\binom{K}{K'}$ is said to have \emph{finite type} if 
the group $J(K')$ is finite. Otherwise, it is said to have infinite type.
\end{definition}

The main interest of this definition is that it coincides with the family of finitely generated generalized $J$-groups:

\begin{proposition}\label{prop:finite-type_=_finiment_engendre}
A nontrivial generalized $J$-group is finitely generated if and only if it has finite-type. Moreover, a finite-type generalized $J$-group is finitely presented.
\end{proposition}
\begin{proof}
Let $J\binom{K}{K'}$ be a generalized $J$-group. 
If $J\binom{K}{K'}$ has finite index in $J(K)$, then it is finitely presented since $J(K)$ is finitely presented (this is a consequence of the Reidemeister Schreier method).

Conversely, assume that $J\binom{K}{K'}$ has infinite type. If $p\leqslant 2$, then since the elements of $K'$ are all finite, the group $J(K')$ is always finite and $J\binom{K}{K'}$ has finite-type. We can then assume that $p\geqslant 3$.

Consider the morphism $J(K)\to W_K^+$ given by Proposition \ref{prop:center_parent_generalized_j_groups}. The image of $J\binom{K}{K'}$ under this morphism is the group $W_K^+(K')$. If $K'$ contains 2 or less entries different from $1$, then $J(K')$ is a finite abelian group and $J\binom{K}{K'}$ has finite type. We can then assume that $K'$ contains at least three entries different from $1$, and we can apply Lemma \ref{lem:finite_index_J_or_coxeter}. We obtain that the group $J(K)/J\binom{K}{K'}\simeq J(K')$ is finite if and only if the group $W_K^+/W_K^+(K')\simeq W_{K'}^+$ is finite. Thus $J\binom{K}{K'}$ has finite index in $J(K)$ if and only if $W_K^+(K')$ has finite index in $W_K^+$.

If $J\binom{K}{K'}$ is nontrivial and has infinite index in $J(K)$, then $W_K^+(K')$ is nontrivial and has infinite index in $W_K^+$. The group $W_K^+(K')$ is thus infinitely generated by Proposition \ref{prop:infinitely_generated_in_coxeter}. Since $J\binom{K}{K'}$ surjects onto $W_K^+(K')$, we have the result.
\end{proof}

\begin{remark}
Notice that the above proposition fails if we allow the elements of $K'$ to be infinite. For instance, the group $J\jg{\infty&\infty\\1&\infty}$ is isomorphic to $\Z$ but has infinite index in its parent generalized $J$-group $\Z^2$.
\end{remark}

\begin{remark}
The above proposition applies in particular to the case $p=3$. It answers \cite[Question 1.1.2]{TheseIgor} about finite generatedness of (classical) $J$-groups. However, even if we showed that a finite-type generalized $J$-group is finitely generated, we have not yet showed that it can be generated by a finite number of reflections. This will be seen in Section \ref{sec:gen_j-groups_as_torsion_quotients}, where we realize all finite-type generalized $J$-groups as torsion quotients of $J$-braid groups (see Remark \ref{rem:finite-type_j-groups_generated_by_reflections})
\end{remark}

\begin{corollary}[\textbf{Finite generalized $J$-groups}]\label{cor:finite_generalised_j_groups}
A nontrivial generalized $J$-group is finite if and only if its generalized parent $J$-group is finite.\end{corollary}
\begin{proof}
If the generalized parent $J$-group of a generalized $J$-group is finite, then the generalized $J$-group is finite as a subgroup of its generalized parent $J$-group.

Conversely, assume that $J\binom{K}{K'}$ is a nontrivial finite generalized $J$-group. In particular, $J\binom{K}{K'}$ is finitely generated. By Proposition \ref{prop:finite-type_=_finiment_engendre}, the index of $J\binom{K}{K'}$ in its generalized parent $J$-group is finite. Since $J\binom{K}{K'}$ is finite, this implies that its generalized parent $J$-group is also finite.
\end{proof}

Combining this corollary with Lemma \ref{lem:finite_index_J_or_coxeter}, we obtain a complete classification of finite generalized $J$-groups.

\begin{corollary}\label{cor:finite_generalized_j_groups_are_crgs}
A nontrivial generalized $J$-group is finite if and only if it is reflection isomorphic to a complex reflection group of rank two. Conversely, every complex reflection group of rank 2 is realized as a finite generalized $J$-group.
\end{corollary}

\begin{proof}
Using Corollary \ref{cor:finite_generalised_j_groups} and Lemma \ref{lem:finite_index_J_or_coxeter}, the generalized $J$-group $J\binom{K}{K'}$ can only be finite if $p\leqslant 3$. If $p\leqslant 2$, the finite generalized $J$-groups are direct product of finite cyclic groups, which are (reducible) rank two complex reflection groups. If $p=3$, then we know that the parent $J$-group $J(K)$ is a complex reflection group (see \cite[Theorem 1.2]{AA} or simply \cite[Tables 1-3]{BMR}). The considered $J$-group is then a complex reflection group as the normal closure in a finite reflection group of a set of reflections. The converse statement that every complex reflection group of rank 2 is realized as a finite $J$-group is shown in \cite[Theorem 1.2]{AA}.
\end{proof}

\subsection{Conjugacy of reflections}
We now study the conjugacy of reflections in generalized $J$-groups. First, using the description of centralizers obtained in alternating polygonal Coxeter groups, we can describe the centralizer of reflections in parent generalized $J$-groups.

\begin{proposition}[\textbf{Centralizer of a reflection}]\label{prop:centralizer_reflections}
Let $i\in \intv{1,p}$, and let $n\in \Z^*$ be such that $a_i^n$ is nontrivial. The centralizer of $a_i^n$ in $J(K)$ is abelian and generated by $a_i$ and $Z(J(K))$.
\end{proposition}
\begin{proof}
If $J(K)$ is abelian, then $J(K)=Z(J(K))$ is the centralizer of any element and the result is immediate.
If $p\geqslant 3$ and $J(K)$ is finite, then $W:=J(K)$ is an irreducible complex reflection group of rank $2$ and $a_i$ is a reflection of $W$. The center of $J(K)$ is cyclic and generated by $z:=s_1\cdots s_p$ after Proposition \ref{prop:center_parent_generalized_j_groups}.  The group $H:=\langle a_i,z\rangle$ is an abelian group with order $k_i|Z(W)|$. Moreover, $H$ is included in the centralizer of $a_i^{n}$, and it is sufficient to show that this centralizer has cardinality $k_i|Z(W)|$.

In the Shephard-Todd notation, we have either $W=G(2l,2,2)$ for $l\geqslant 2$, of $W\in \{G_7,G_{11},G_{19}\}$. In each case, we have $|W|=|Z(W)|^2$. If $K\in (2,3,3),(2,3,4)$ or $(2,3,5)$, the result is obtained by computer (using for instance the CHEVIE package of {\tt GAP3} \cite{gap3}). If $K=(l,2,2)$ for some $l\geqslant 2$, then $J(K)\simeq G(2l,2,2)$. More precisely, the identification is given by
\[a_1\mapsto \begin{pmatrix} e^{2i\pi/l}&0\\0&1\end{pmatrix},~a_2\mapsto \begin{pmatrix} 0&1\\1&0\end{pmatrix},~a_3\mapsto \begin{pmatrix} 0&e^{i\pi/l}\\e^{-i\pi/l}&0\end{pmatrix}.\]
The centralizers of these elements and their nontrivial powers are easily computed using this representation, which gives the desired result.

Lastly, assume that $p\geqslant 3$ and that $J(K)$ is infinite. Let $g\in C_{J(K)}(a_i^n)$. writing $\pi:J(K)\to W_K^+$ to denote the canonical quotient, we then have $\pi(g)\in C_{W_K^+}(\pi(a_i)^n)$. Now, using Lemma \ref{lem:finite_index_J_or_coxeter}, the group $W_K^+$ is infinite since $J(K)$ is as well, hence Proposition \ref{prop:centralizer_rot_polyg_coxeter_group} applies. This shows that $\pi(g)\in\langle \pi(a_i)\rangle$. Since $\ker(\pi)=\langle a_1a_2\cdots a_p\rangle$, in turn equal to $Z(J(K))$ by Proposition \ref{prop:center_parent_generalized_j_groups},  we obtain that $g\in \langle a_i,a_1a_2\cdots a_p\rangle$. This concludes the proof.
\end{proof}

Using this result, we can describe the conjugacy classes of reflections in a generalized $J$-group. As usual, we begin with the case of generalized parent $J$-groups. 

\begin{lemma}\label{lem:conj_refl_parent}
Two reflections $r,r'\in J(K)$ are conjugate in $J(K)$ if and only if we can write $r=gs_i^ng^{-1}$ and $r'=hs_i^nh^{-1}$ with $g,h\in J(K)$ and $n\in \intv{1,k_i-1}$.
\end{lemma}
\begin{proof}
By definition of $R(J(K))$, we can write $r=gs_i^ng^{-1}$ and $r'=hs_j^mh^{-1}$ with $i,j\in \intv{1,p}$ and $n\in \intv{1,k_i-1}$, $m\in \intv{1,k_j-1}$. The two reflections $r$ and $r'$ are then conjugate in $J(K)$ if and only if the reflections $s_i^n$ and $s_j^m$ are conjugate. If $j=i$ and $m=n$, then $s_i^n=s_j^m$ and $r,r'$ are conjugate in $J(K)$. Conversely if $s_i^n$ and $s_j^m$ are conjugate, then their images in the abelianization $A$ of $J(K)$ are equal. By Definition \ref{def:generalized_j-groups}, $A$ decomposes as a direct product $\Z/k_1\Z\times \cdots \times \Z/k_i\Z$, generated by the images of the $s_i$'s. If the images of $s_i^n$ and $s_j^m$ are equal, we then obtain $i=j$ and $m=n$ modulo $k_i$. Since $m,n\in \intv{1,k_i}$ by assumption, we have $m=n$ as we wanted to show.
\end{proof}

In particular, we see that a reflection is never conjugate to any of its powers other than itself. For linear reflection groups this is easily shown by looking at eigenvalues, but here we have the result without using a linear representation of generalized $J$-groups.

\begin{proposition}[\textbf{Conjugacy classes of reflections}]\label{prop:conj_reflections_j-groupe}
Let $Q$ denote the quotient of $J(K')$ by the image of $Z(J(K))$ under the projection $J(K)\to J(K')$.

For $i\in \intv{1,p}$, let $\bbar{s_i}$ denote the image of $s_i$ in $Q$. Let $n\in \N$ be such that $s_i^n\in J\binom{K}{K'}$ is nontrivial. 
\begin{enumerate}[(a)]
\item For $g,h\in J(K)$, the elements $gs_i^{n}g^{-1}$ and $hs_i^{n}h^{-1}$ are conjugate in $J\binom{K}{K'}$ if and only if the images of $g$ and $h$ in $Q$ lie in the same left $\langle \bbar{s_i}\rangle$-coset.
\item There is a natural bijection between the set of conjugacy classes in $J\binom{K}{K'}$ of reflections conjugate to $s_i^{n}$ in $J(K)$ and the left $\langle \bbar{s_i}\rangle$-cosets in $Q$.
\end{enumerate}
\end{proposition}
\begin{proof}
Point $(b)$ is a direct consequence of point $(a)$, we thus only have to prove the latter. Since $J\binom{K}{K'}$ is normal in $J(K)$, it is sufficient to consider the case where $h=1$. Let $r:=s_i^{n}$. 

By definition, $r$ and $grg^{-1}$ are conjugate in $J\binom{K}{K'}$ if and only if there is an element of $J\binom{K}{K'}$ which conjugates $r$ to $grg^{-1}$. The set of elements of $J(K)$ which conjugate $r$ to $grg^{-1}$ is the left coset $gC_{J(K)}(r)$. The intersection $gC_{J(K)}(r)\cap J\binom{K}{K'}$ is nontrivial if and only if there is some $x\in C_{J(K)}(r)$ whose image in $J(K')$ is equal to that of $g$. By Proposition \ref{prop:centralizer_reflections}, the centralizer $C_{J(K)}(r)$ is generated by $r$ and $Z(J(K))$. Thus $r$ and $grg^{-1}$ are conjugate in $J\binom{K}{K'}$ if and only if there is some $m\in \Z^*$, $z\in Z(J(K))$ such that the images of $g$ and of $s_i^mz$ in $J(K')$ are equal. Taking the quotient by the image of $Z(J(K))$ in $J(K)$. We obtain that $r$ and $grg^{-1}$ are conjugate in $J\binom{K}{K'}$ if and only if there is some $m\in \Z^+$ such that the image of $g$ in $Q$ in $\bbar{s_i}^m$. This is the desired result.
\end{proof}

The fact that Proposition \ref{prop:conj_reflections_j-groupe} gives a complete description of the conjugacy classes of reflections in $J\binom{K}{K'}$ comes from Lemma \ref{lem:conj_refl_parent}: if $r,r'\in R(J\binom{K}{K'})$ are conjugate in $J\binom{K}{K'}$, then they are in particular conjugate in $J(K)$. By Lemma \ref{lem:conj_refl_parent}, this implies that we are in the situation of Proposition \ref{prop:conj_reflections_j-groupe}.

\begin{corollary}[\textbf{Counting conjugacy classes of reflections}]\label{cor:counting_conjugacy_classes_of_reflections}
Let $Q$ denote the quotient of $J(K')$ by the image of $Z(J(K))$ under the projection $J(K)\to J(K')$. For $i\in \intv{1,p}$, let $\bbar{s_i}$ denote the image of $s_i$ in $Q$.
\begin{enumerate}[(a)]
\item The group $Q$ is finite if and only if $J\binom{K}{K'}$ has finite-type.
\item In this case, the number of conjugacy classes of reflections in $J\binom{K}{K'}$ is given by
\[\sum_{i=1}^p \sum_{j=1}^{k_i/k_i'-1} [Q:\langle \bbar{s_i}\rangle].\]
\end{enumerate}
\end{corollary}
\begin{proof}
$(a)$ If $J\binom{K}{K'}$ has finite-type, then $Q$ is finite as a quotient of the finite group $J(K')$. Conversely, assume that $J\binom{K}{K'}$ has infinite type. We can assume that $p\geqslant 3$, since $p\leqslant 2$ forces $J\binom{K}{K'}$ to have finite-type. Moreover, we can assume that $K'$ contains 3 or more entries different than $1$, since $J(K')$ is finite and abelian if this is not the case. The center of $J(K)$ is then generated by $s_1\cdots s_p$ by Proposition \ref{prop:center_parent_generalized_j_groups}, and the group $Q$ is $W_{K'}^+$. We then conclude by Lemma \ref{lem:finite_index_J_or_coxeter} that if $J\binom{K}{K'}$ has infinite type, then $J(K')$ and $W_{K'}^+$ are infinite. 

$(b)$ After Lemma \ref{lem:conj_refl_parent}, every reflection in $J\binom{K}{K'}$ is conjugate in $J(K)$ to exactly one reflection of the form $s_i^n$ for $n\in \intv{1,k_i-1}$. Moreover, by Lemma \ref{lem:reflexion_belonging_to_j_group}, $n$ must be a multiple of $k_i'$. The number of conjugacy classes of reflections in $J\binom{K}{K'}$ which are conjugate to $s_i^n$ is given by Proposition \ref{prop:conj_reflections_j-groupe} $(b)$.
\end{proof}

\subsection{Center and torsion elements}
As a consequence of the determination of centralizers of reflections, we can determine the center of arbitrary generalized $J$-groups. This was already proven by the second author in \cite[Theorem 2.2.29]{TheseIgor} for the family of $J$-reflection groups (which are particular (classical) $J$-group).

\begin{corollary}[\textbf{Center of generalized $J$-groups}]\label{cor:center_generalized_j_groups}
The center of $J\binom{K}{K'}$ is given by $Z(J\binom K{K'})=J\binom K{K'}\cap Z(J(K))$. Moreover, if $J\binom K{K'}$ is nontrivial, it is abelian if and only if $J(K)$ is also abelian. 
\end{corollary}

\begin{proof}
If $p\leqslant 2$, then $J(K)$ is abelian and the result is immediate. Moreover, if $J(K)$ is finite, then we are in the well-known case of irreducible rank two complex reflection groups.

Assume thus that $p\geqslant 3$ and that $J(K)$ is infinite. Consider the natural morphism $J(K)\to W_K^+$. This morphism restricts to a surjective morphism $\pi:J\binom{K}{K'}\to W_K^+(K')$. Since $W_K^+(K')$ has trivial center by Corollary \ref{cor:center_wkk'}, the center of $J\binom{K}{K'}$ is included in $\ker \pi$. Since the kernel of the morphism $J(K)\to W_K^+$ is $Z(J(K))$, the kernel $\ker\pi$ is the intersection of $J\binom{K}{K'}$ with $Z(J(K))$. We then have $Z(J(\binom{K}{K'}))\subset J\binom{K}{K'}\cap Z(J(K))$. The other inclusion is immediate.

Now, if $J(K)$ is abelian, then $J\binom K{K'}$ is abelian as a subgroup of $J(K)$. Conversely, if $J(K)$ is nonabelian, then $p\geqslant 3$ and $Z(J(K))$ is cyclic and generated by $a_1\cdots a_p$. However, $J\binom K{K'}$ contains a nontrivial reflection, which is noncentral by Proposition \ref{prop:conj_reflections_j-groupe}. We then have $Z(J\binom K{K'})=Z(J(K))\cap J\binom K{K'}\subsetneq J\binom K{K'}$ and $J\binom K{K'}$ is not abelian. Moreover, if $p\geqslant 3$, the group $J(K)$ is not abelian since its inner automorphism group $W_K^+$ is non-trivial. This concludes the proof.
\end{proof}

By Corollary \ref{cor:finite_generalised_j_groups}, finite generalized $J$-groups have finite-type. Thus a generalized $J$-group is either finite, or infinite and finite-type, or infinite-type. Moreover, theses cases are mutually exclusive. We will show that these cases can be recovered by only considering the centers of generalized $J$-groups.

First, we show that finiteness of a generalized parent $J$-group can be characterized by only considering its center. This result was originally showed by the second author \cite[Corollary 2.2.37]{TheseIgor} for the case $p=3$.

\begin{proposition}[\textbf{Order of center of parent generalized $J$-groups}]\label{prop:order_center_parent_generalized_j_groups}~\\
A generalized parent $J$-group is finite if and only if its center is finite.
\end{proposition}

\begin{proof}
If $p\leqslant 2$, the group $J(K)$ is abelian so that the result is trivial.\\
If $p=3$, the statement is that of \cite[Corollary 2.2.37]{TheseIgor}.\\
Assume now that $p\geqslant 4$ (in which case $J(K)$ is always infinite), we have to show that $Z(J(K))$ is infinite. If $L$ is a tuple of length $p'\geqslant 3$ obtained from $K$ by deleting entries, then there is a natural quotient $J(K)\to J(L)$ which induces a quotient $Z(J(K))\to Z(J(L))$. It is then sufficient to show that $Z(J(L))$ is infinite. \\ Assume that $p=4$. If $K$ contains a subset $L$ of cardinal $3$ such that $J(L)$ is infinite, then case $p=3$ applies and $Z(J(L))$ is infinite. All remaining cases for $p=4$ are of the form $(2,2,n,m)$ with $n,m\geqslant 2$. By \cite[Theorem 2.2.10 and Remark 3.1.7]{TheseIgor}, we have $J(2,2,n,m)\cong J\jg{2 & 2n & 2m\\ 1 & 2 & 2}$, in which case the result follows from \cite[Corollary 2.2.37]{TheseIgor}.\\ If $p\geqslant 5$, then $K$ admits a subtuple $L$ of length $4$. The group $J(L)$ is infinite in this case by Proposition \ref{prop:center_parent_generalized_j_groups}, and $Z(J(L))$ is infinite by the case $p=4$. This concludes the proof.
\end{proof}

\begin{proposition}[\textbf{Order of center of generalized $J$-groups}]\label{order_center_generalized_j_groups}~\\Let $J\binom{K}{K'}$ be a generalized $J$-group.
\begin{itemize}
    \item If $J\binom{K}{K'}$ is finite, then its center is either a finite cyclic group or a direct product of two finite cyclic groups.
    \item If $J\binom{K}{K'}$ is infinite and finite-type, then its center is infinite cyclic.
    \item If $J\binom{K}{K'}$ is infinite-type, then its center is trivial (and $J\binom{K}{K'}\simeq W_K^+(K')$). 
\end{itemize}
\end{proposition}

\begin{proof}
By Corollary \ref{cor:finite_generalised_j_groups}, a generalized $J$-group is finite if and only if its generalized parent $J$-group is finite. If this is the case, then we are in the well-known case of rank two complex reflection groups. This covers the first case.\\
If $J\binom{K}{K'}$ is infinite and finite type, then $J(K)$ is infinite and $J(K')$ is finite. We know by Corollary \ref{cor:center_generalized_j_groups} that $Z(J\binom K{K'})=J\binom K{K'}\cap Z(J(K))$ and by Proposition \ref{prop:order_center_parent_generalized_j_groups} that $Z(J(K))$ is infinite cyclic. Thus, since $J(K')$ is finite, the group $J\binom K{K'}$ contains at least one non-trivial element of $Z(J(K))$ (otherwise the quotient $J(K)/J\binom K{K'}\cong J(K')$ would contain an element of infinite order, namely the image of the generator of $Z(J(K))$ ). This covers the second case.\\
If $J\binom{K}{K'}$ is infinite-type, then $J(K')$ is infinite. The group $J\binom K{K'}$ cannot contain a central element of $J(K)$. Indeed, if it did, the image of the generator of $Z(J(K))$ in $J(K)/J\binom K{K'}\cong J(K')$ would have finite order, which is false by Proposition \ref{prop:order_center_parent_generalized_j_groups}. This shows that $Z(J\binom K{K'})=1$. We conclude in particular that if $J(K')$ is infinite, the restriction of the quotient $J(K)\to W_K^+$ to $J\binom K{K'}$ is injective. This concludes the proof.
\end{proof}

We finish this section by giving a complete description of the torsion elements in a infinite generalized $J$-group.

\begin{corollary}[\textbf{Torsion in generalized $J$-groups}]\label{cor:torsion_generalized_j_groups}
The torsion elements of an infinite generalized $J$-group are precisely its finite order reflections. 
\end{corollary}

\begin{proof}
Since $J\binom K{K'}\cap R(J(K))=R(J\binom K{K'})$, it is enough to show the result for $J(K)$. If $p\leqslant 2$, then $J\binom{K}{K'}$ is a direct product of cyclic group, which is infinite if and only if at least one of the factors is infinite cyclic. The result is immediate in this case. Now, assume that $p\geqslant 3$, and write $\pi:J(K)\to W_K^+$ for the natural quotient. If $x\in J(K)$ is a torsion element, $\pi(x)$ is a torsion element in $W_K^+$. Moreover, by Proposition \ref{prop:order_center_parent_generalized_j_groups}, $x$ is not central, thus $\pi(x)$ is nontrivial. By Proposition \ref{TorsionPolygonalCox}, $\pi(x)$ is conjugate to a nontrivial power of some $\pi(a_i)$ with $k_i<\infty$ for $i\in \intv{1,d}$. Up to conjugating $x$ we can assume that $\pi(x)=\pi(a_i)^q$ for $q<k_i$.  We then have $x=a_i^q(a_1a_2\cdots a_p)^r$ for some integer $r$. By Proposition \ref{prop:order_center_parent_generalized_j_groups}, the order of $a_1a_2\cdots a_p$ is infinite, which forces $r=0$ since $x$ is a torsion element. We obtain $x=a_i^q$, as we wanted to show.
\end{proof}

\section{$J$-braid groups and torsion quotients}\label{sec:j-braid_groups_and_torsion_quotients}
\subsection{Reminders on $J$-braid groups}
In \cite{VCRG}, the second author studies a particular family of (classical) $J$-groups called \emph{$J$-reflection groups}, which contains in particular all finite $J$-groups (i.e. all complex reflection groups of rank $2$). Later in \cite{BraidsJIgor}, the second author introduced a family of so-called \emph{$J$-braid group}, naturally attached $J$-reflection groups. This construction generalizes the definition of braid group of complex reflection groups of rank 2. The definition of $J$-braid group given in \cite{BraidsJIgor} is via a presentation by generators and relations involving two positive integers $n,m$. In \cite{BraidsJIgor}, the integers $n,m$ were assumed to be coprime. This assumption was later removed in \cite{Necklace}, giving rise to \emph{generalized} $J$-braid groups, attached to \emph{generalized} $J$-reflection groups. For readability, we will call these groups $J$-braid groups if needed.

For the remainder of this section we fix $n,m\in \N_{\geqslant 1}$ and we write $m=qn+r$ with $0\leqslant q$ and $0\leqslant r\leqslant n-1$. We also define $d=n\wedge m$.

\begin{definition}[\textbf{$J$-braid groups}]\label{BraidPresentation1}\hfill
\begin{itemize}
\item The group $\B_*^*(n,m)$ is defined by the following presentation:
\begin{subequations}
\label{classicalBraidPres1} 
\begin{align}
&(1) \,\, \mathrm{Generators}\!:\,  \{x_1,\dots,x_n,y,z\};\notag\\
&(2) \,\, \mathrm{Relations}\!: \,\notag\\
  &  x_1\cdots x_nyz=zx_1\cdots x_ny,\label{classicalBraidPres1:1}\\
       & x_{i+1}\cdots x_nyz\delta^{q-1}x_1\cdots x_{i+r}=x_i\cdots x_nyz\delta^{q-1}x_1\cdots x_{i+r-1}, \, \forall 1\leqslant i \leqslant n-r,\label{classicalBraidPres1:2}\\
      &  x_{i+1}\cdots x_nyz\delta^qx_1\cdots x_{i+r-n}=x_i\cdots x_nyz\delta^qx_1\cdots x_{i+r-n-1},\, \forall n-r+1\leqslant i \leqslant n,\label{classicalBraidPres1:3}
\end{align}
\end{subequations}
where indices are taken modulo $n$ and where $\delta$ denotes $x_1\cdots x_ny$.
\item The group $\B^*(n,m)$ is the quotient of $\B_*^*(n,m)$ by the normal closure of $y$.
\item The group $\B_*(n,m)$ is the quotient of $\B_*^*(n,m)$ by the normal closure of $z$.
\item The group $\B(n,m)$ is the quotient of $\B_*^*(n,m)$ by the normal closure of $\{y,z\}$.
\end{itemize}
\end{definition}

\begin{remark}\label{rem:circular_are_j-braid}
By definition, in the presentation of $\B(n,m)$, only relations of the form \eqref{classicalBraidPres1:2} and \eqref{classicalBraidPres1:3} remain, and they state that the product $x_i\cdots x_{i+m}$ (with indices seen mod $n$) does not depend on $i$. In other words, the presentation given for the group $\B(n,m)$ coincides with the presentation given in \cite{GarnierHoso} of the \emph{circular group} denoted there by $G(n,m)$. Circular groups then appear as particular examples of $J$-braid groups. 
\end{remark}

This combinatorial definition of $J$-braid groups was justified by the following theorem:

\begin{theorem}{\cite[Theorem 2.2.10 and Remark 3.1.7]{TheseIgor}}
Let $b,c\in \N_{\geqslant 1}$. The group $J\jg{ k & bn & cm \\
1 & n & m}$ is isomorphic to $\B_*^*(n,m)/\llangle x_1^k,\dots,x_n^k,y^b,z^c\rrangle$.
\end{theorem}

In the following sections, we are going to consider all quotients of $J$-braid groups obtained by adding torsion to the generators. For readability, it is convenient to name the conjugates of the generators in a $J$-braid group. We follow the notation of \cite{Necklace}:

\begin{definition}[\textbf{Braid reflections}]
We call \emph{braid reflections} the conjugates of the generators of $\B_*^*(n,m)$, $\B^*(n,m)$, $\B_*(n,m)$ and $\B(n,m)$.
\end{definition}

Notice that we did not give a topological meaning to $J$-braid groups (yet). The terminology of braid reflections relies (for now) entirely on the analogy with the case of $J$-braid groups which are complex braid groups. We thus chose to use the term \emph{braid reflections} instead of \emph{braided reflection}.

By construction of $J$-braid groups, we have a commutative square of groups
\begin{equation}\label{eq:square_j-braid_groups}
\begin{tikzcd}
	\B^*_*(n,m) & \B^*(n,m) \\
	\B_*(n,m) & \B(n,m)
	\arrow["y=1",from=1-1, to=1-2]
	\arrow["z=1",from=1-2, to=2-2]
	\arrow["z=1"',from=1-1, to=2-1]
	\arrow["y=1"',from=2-1, to=2-2]
\end{tikzcd}
\end{equation}
Contrary to \cite{BraidsJIgor} and \cite{Necklace}, we do not make the assumption that $m\geqslant 2$ (resp. $n\geqslant 2$) in order to define $\B_*(n,m)$ (resp. $\B^*(n,m)$). In particular, the above square is always defined, but we need to be precise on the results involving groups $\B^*(1,m)$ and $\B_*(n,1)$.

Presentation \eqref{classicalBraidPres1} reads \[\B_*^*(1,m)=\langle x,y,z~|~xyz=zxy,~ yz(xy)^{m-1}x=xyz(xy)^{m-1}\rangle,\] so that $\B^*(1,m)$ is isomorphic to $\Z^2$. Now, there exists an isomorphism between $\B^*(1,m)$ and $\B_*(m,1)$ which sends $y$ to $z$ (see Corollary \ref{cor:swap_parameters_torsion_quotient}), hence $\B_*(m,1)$ is again isomorphic to $\Z^2$. Notice that if $1\in\{n,m\}$ we also have $\B(n,m)\cong \Z$ by \cite[Corollary 2.11]{GarnierHoso}.

We mentioned above that circular groups are particular cases of $J$-braid groups. In fact, the family of $J$-braid groups coincides (up to abstract group isomorphism) with the family of circular groups: 
\begin{theorem}[\textbf{Isomorphism type of $J$-braid groups}]\cite[Theorem 3.1.23 and 3.1.31]{TheseIgor}\label{TheoremBraidIsoIntroCircularFR}
\begin{itemize}
    \item The group $\B_*^*(n,m)$ is isomorphic to $\B(d+2,d+2)$.
    \item The group $\B_*(n,m)$ is isomorphic to  $\B(d+1,(d+1)m')$.
    \item The group $\B^*(n,m)$ is isomorphic to $\B((d+1)n',d+1)$.
\end{itemize}
Note however that these various isomorphisms do not send braid reflections to braid reflections in general.
\end{theorem}
\begin{proof}
The only cases which are not covered in \cite[Theorem 3.1.23 and 3.1.31]{TheseIgor} are the isomorphisms $\B_*(n,1)\simeq \B(2,2)\simeq \B^*(1,m)$, which we already showed above.
\end{proof}

Using this result, it is rather easy to determine the center of $J$-braid groups:

\begin{corollary}[\textbf{Center of $J$-braid groups}]\cite[Corollary 3.1.35]{TheseIgor}\\
Let $m':=\frac{m}{d}$ and $n':=\frac{n}{d}$.
\begin{itemize}
\item The center of $\B_*^*(n,m)$ is infinite cyclic and generated by $\delta^{m'}z^{n'}$.
\item The center of $\B^*(n,m)$ is infinite cyclic and generated by $\delta^{m'}z^{n'}$, except for $\B^*(1,m)\simeq \Z^2$.
\item The center of $\B_*(n,m)$ is infinite cyclic and generated by $\delta^{m'}$, except for $\B_*(n,1)\simeq \Z^2$. 
\item The center of $\B(n,m)$ is infinite cyclic and generated by $\delta^{m'}$, except if $n=1$ or $m=1$ or $n=m=2$, in which case $\B(n,m)$ is abelian.
\end{itemize}
\end{corollary}
\begin{proof}
The fourth case was obtained in \cite[Corollary 2.11]{GarnierHoso}, thus it only remains to observe the result for $\B^*(1,m)$ and for $\B_*(n,1)$, which is immediate.
\end{proof}

\begin{lemma}[\textbf{Conjugacy classes of braid reflections}]\label{lem:conjugacy_classes_of_braid_reflections}
The pairs of conjugate generators of $\B_*^*(n,m)$ (resp. $\B^*(n,m),\B_*(n,m),\B(n,m)$) are exactly the pairs of the form $\{x_i,x_{i+kd}\}$ for $k\in \N$ (where we see the indices mod $n$). \\In particular, a complete set of representatives of conjugacy classes of braid reflections is given by
\begin{itemize}
\item $\{x_1,\ldots,x_d,y,z\}$ for $\B^*_*(n,m)$.
\item $\{x_1,\ldots,x_d,z\}$ for $\B^*(n,m)$.
\item $\{x_1,\ldots,x_d,y\}$ for $\B_*(n,m)$.
\item $\{x_1,\ldots,x_d\}$ for $\B(n,m)$.
\end{itemize}
\end{lemma}
\begin{proof}
Consider the presentation of $\B_*^*(n,m)$. For $i\in \intv{1,n-r}$, Equation \eqref{classicalBraidPres1:2} implies that $x_i$ is conjugate to $x_{i+r}$ by $x_{i+1}\cdots x_ny\delta^{q-1}x_1\cdots x_{i+r-1}$. Similarly for $i\in \intv{n-r+1,n}$, Equation \eqref{classicalBraidPres1:3} implies again that $x_i$ is conjugate to $x_{i+r-n}$. Seeing the indices mod $n$, we obtain that $x_i$ is conjugate to $x_{i+r}$ for all $i\in \intv{1,n}$. Since the gcd of $n$ and $r$ is equal to that of $n$ and $m$ (i.e. to $d$), we obtain that $x_{i+kd}$ is always conjugate to $x_i$ for $i\in \intv{1,n}$ and $k\in \N$. Since this holds in $\B_*^*(n,m)$, it also holds in its quotients $\B^*(n,m),\B_*(n,m)$ and $\B(n,m)$. 

Conversely, the presentation of $\B^*_*(n,m)$ induces a presentation of its abelianization $A$. We obtain that $A$ is generated by $\bar{x_1},\ldots,\bar{x_n},\bar{y},\bar{z}$ with the only relations being $\bar{x_{i}}=\bar{x_{i+r}}$ for $i\in \intv{1,n-r}$ and $\bar{x_i}=\bar{x_{i+r-n}}$ for $i\in \intv{n-r+1,n}$. The group $A$ is then free abelian generated by $\bar{x_1},\ldots,\bar{x_d},\bar{y},\bar{z}$. Since two conjugate elements in $\B^*_*(n,m)$ must be identified in $A$, we deduce that no two elements of $\{x_1,\ldots,x_d,y,z\}$ are conjugate in $\B_*^*(n,m)$, which finishes the proof in this case. Similar computations give the abelianization of $\B^*(n,m),\B_*(n,m)$ and of $\B(n,m)$, which gives the result in these cases.
\end{proof}

\subsection{Torsion quotients of $J$-braid groups}
In this section we fix two positive integers $n,m$, and we let $d:=n\wedge m$ denote the gcd of $n$ and $m$. We also fix a $d$-tuple $K:=(k_1,\ldots,k_d)$ of elements in $\N_{\geqslant 2}\cup \{\infty\}$, along with two elements $b,c\in \N_{\geqslant 2}\cup\{\infty\}$. 

\begin{definition}[\textbf{Torsion quotient}]\label{def:torsion_quotient_j-braid}\hfill\begin{itemize}
\item The group $\B_*^*(n,m;K,b,c)$, defined as the quotient of $\B_*^*(n,m)$ by the normal closure of $\{x_1^{k_1},\ldots,x_d^{k_d},y^b,z^c\}$ is called a \emph{torsion quotient} of $\B_*^*(n,m)$. 
\item The group $\B^*(n,m;K,c)$, defined as the quotient of $\B^*(n,m)$ by the normal closure of $\{x_1^{k_1},\ldots,x_d^{k_d},z^c\}$ is called a \emph{torsion quotient} of $\B^*(n,m)$. 
\item The group $\B_*(n,m;K,b)$, defined as the quotient of $\B_*(n,m)$ by the normal closure of $\{x_1^{k_1},\ldots,x_d^{k_d},y^b\}$ is called a \emph{torsion quotient} of $\B_*(n,m)$. 
\item The group $\B(n,m;K)$, defined as the quotient of $\B(n,m)$ by the normal closure of $\{x_1^{k_1},\ldots,x_d^{k_d}\}$ is called a \emph{torsion quotient} of $\B(n,m)$. 
\end{itemize}
\end{definition}

By construction of torsion quotients, square \eqref{eq:square_j-braid_groups} induces a square between torsion quotients:
\begin{equation}\label{eq:square_torsion_quotients}
\begin{tikzcd}
	\B^*_*(n,m;K,b,c) & \B^*(n,m;K,c) \\
	\B_*(n,m;K,b) & \B(n,m;K)
	\arrow["y=1",from=1-1, to=1-2]
	\arrow["z=1",from=1-2, to=2-2]
	\arrow["z=1"',from=1-1, to=2-1]
	\arrow["y=1"',from=2-1, to=2-2]
\end{tikzcd}
\end{equation}
where we abusively denote by $y$ and $z$ the respective images of $y$ and $z$ in the torsion quotients.

At first glance, we only imposed torsion relations on the first $d$ generators $x_i$ of the presentation of a $J$-braid group in order to define torsion quotients. However, by Lemma \ref{lem:conjugacy_classes_of_braid_reflections}, two generators $x_i,x_j$ of a $J$-braid group are conjugate if and only if $i$ and $j$ are equivalent modulo $d$. In particular, the image of $x_{d+i}$ in the associated torsion quotient has order $k_{i}$ for all $i\in \intv{1,d}$.  

Conversely, for $(k_1,k_2,\ldots,k_n)\in \N_{\geqslant 2}\cup \{\infty\}$, the quotient of a $J$-braid group by the relations $x_i^{k_i}=1$ for all $i\in \intv{1,n}$ (along with $y^b=z^c=1$ if needed) is easily shown to be a torsion quotient in the sense of Definition \ref{def:torsion_quotient_j-braid}. 

\begin{remark}\label{rem:b=1_c=1}
We could theoretically allow for $b=1$ or $c=1$ in the definition of torsion quotient. In this case every torsion quotient could be described as a torsion quotient of $\B_*^*(n,m)$, with for instance $\B_*(n,m;K,c)=\B_*^*(n,m;K,1,c)$. 
\end{remark}

\begin{remark}
Assume that $d=1$. In this case, the tuple $K$ is actually a single element $k$. The presentation of $\B_*^*(n,m;k,b,c)$ coincides with the presentation of the \emph{$J$-reflection group} $W_b^c(k,bn,cm)$ given in \cite[Theorem 2.29]{VCRG}. More generally, if $K=(k,\ldots,k)$, then the presentation of $\B_*^*(n,m;K,b,c)$ coincides with the presentation of the \emph{generalized $J$-reflection group} $W_b^c(k,bn,cm)$. This also holds if $b=1$ or if $c=1$. Torsion quotients of $J$-braid groups then generalize generalised $J$-reflection groups.
\end{remark}

\begin{definition}[\textbf{Reflections}]
Let $\B$ be a $J$-braid group and let $W$ be a torsion quotient of $\B$. The nontrivial powers of conjugates of the images in $W$ of the braid reflections of $\B$ are called the \emph{reflections} of $W$. We denote the set of reflections of $W$ by $R(W)$. We say that two torsion quotients $W$ and $W'$ are \emph{reflection isomorphic} if there exists a group isomorphism $\varphi:W\to W'$ such that $\varphi(R(W))=R(W')$.
\end{definition}

\begin{remark}\label{rem:conjugacy_classes_torsion_quotients}
The description of the conjugacy classes given in Lemma \ref{lem:conjugacy_classes_of_braid_reflections} extends to all torsion quotients of $J$-braid groups, as the proof carries out word by word.
\end{remark}

\subsection{Embedding results for torsion quotients of $J$-braid groups}
In this section we fix two positive integers $n,m$, and we let $d:=n\wedge m$ denote the gcd of $n$ and $m$. We write $m=qn+r$ the euclidean division of $m$ by $n$. We also fix a $d$-tuple $K:=(k_1,\ldots,k_d)$ of elements in $\N_{\geqslant 2}\cup \{\infty\}$, along with two elements $b,c\in \N_{\geqslant 2}\cup\{\infty\}$. 

\begin{notation} We will write $d\cdot K$ for the $dp$-tuple $(k_1,\ldots,k_p,k_1,k_2\ldots,k_p)$ obtained by repeating $K$ $d$-times. Moreover, for $k\in \N\cup \{\infty\}$, $d\cdot k$ will denote the $d$-tuple $d\cdot (k)$.
\end{notation}

Our main tool for studying torsion quotients of $J$-braid groups is the following theorem, which provides an embedding from $\B_*^*(pn,pm;d\cdot K,b,c)$ to $J(K,bn,cm)$ which maps reflections to powers of reflections.

\begin{theorem}[\textbf{Embedding of $J$-braid group}]\label{theo:embedding_j-braid}
The correspondance
\[\begin{cases} x_{l}\mapsto s_{p+1}^js_is_{p+1}^{-j}&\text{for }l\in \intv{1,pn}\text{ and }l=jp+i,~j\in \intv{0,n-1},i\in \intv{1,p},\\y\mapsto s_{p+1}^{n},\\z\mapsto s_{p+2}^{m}.\end{cases}\]
induces an injective morphism from $\B^*_*(pn,pm,d\cdot K,b,c)$ to $J(K,bn,cm)$, whose image is the normal closure of $\{s_1,\ldots,s_{p},s_{p+1}^{n},s_{p+2}^{m}\}$ in $J(K,bn,cm)$.
\end{theorem}

The proof of this theorem is rather intricate and we postpone it until the next section. For the end of this section, we give some corollaries, either for torsion quotients or for the groups $\B^*(n,m),\B_*(n,m)$ and $\B(n,m)$.

\begin{corollary}[\textbf{Embedding of torsion quotients}]\label{cor:embedding_torsion_quotients}
Let $\varphi$ be the embedding of Theorem \ref{theo:embedding_j-braid}.
\begin{itemize}[itemsep=3pt]
\item The morphism $\varphi$ induces an embedding $\B^*_*(pn,pm;d\cdot K,b,c)\to J(K,bn,cm)$, which exhibits $\B_*^*(pn,pm;d\cdot K,b,c)$ as $J\jg{K&bn&cm\\p\cdot 1&n&m}$.
\item The morphism $\varphi$ induces an embedding $\B^*(pn,pm;d\cdot K,c)\to J(K,n,cm)$, which exhibits $\B^*(pn,pm;d\cdot K,c)$ as $J\jg{K&n&cm\\p\cdot 1&n&m}$.
\item The morphism $\varphi$ induces an embedding $\B_*(pn,pm;d\cdot K,b)\to J(K,bn,m)$, which exhibits $\B_*(pn,pm;d\cdot K,b)$ as $J\jg{K&bn&m\\p\cdot 1&n&m}$.
\item The morphism $\varphi$ induces an embedding $\B(pn,pm;d\cdot K)\to J(K,n,m)$, which exhibits $\B(pn,pm;d\cdot K)$ as $J\jg{K&n&m\\p\cdot 1&n&m}$.
\end{itemize}
In all these cases, the notions of reflections for the quotient torsion of $J$-braid groups on the one hand and for the generalized $J$-groups on the other hand coincide.
\end{corollary}
\begin{proof}
The first point is simply a rephrasing of the definitions. The other points follow from considering the square of natural quotients 
\begin{equation*}
\begin{tikzcd}
	J(K,bn,cm) & J(K,n,cm) \\
	J(K,bn,m) & J(K,n,m)
	\arrow[from=1-1, to=1-2]
	\arrow[from=1-2, to=2-2]
	\arrow[from=1-1, to=2-1]
	\arrow[from=2-1, to=2-2]
\end{tikzcd}
\end{equation*}

\end{proof}

\begin{corollary}[\textbf{Permuting torsion coefficients}]\label{cor:permutation_parameters_torsion_quotient}
Let $\sigma\in S_d$ be a permutation, and let $\sigma(K)$ be the $d$-tuple $(k_{\sigma(1)},\ldots,k_{\sigma(d)})$. We have
\begin{itemize}
\item $\B_*^*(n,m;K,b,c)$ and $\B_*^*(n,m;\sigma(K),b,c)$ are reflection isomorphic.
\item $\B^*(n,m;K,c)$ and $\B^*(n,m;\sigma(K),c)$ are reflection isomorphic.
\item $\B_*(n,m;K,b)$ and $\B_*(n,m;\sigma(K),b)$ are reflection isomorphic.
\item $\B(n,m;K)$ and $\B(n,m;\sigma(K))$ are reflection isomorphic.
\end{itemize}
\end{corollary}

\begin{proof}
This is a direct application of Proposition \ref{prop:permutation_of_parameters_generalized_j_groups} (permutation of parameters), since torsion quotients of $J$-braid groups are generalized $J$-groups by Corollary \ref{cor:embedding_torsion_quotients}.
\end{proof}

\begin{corollary}[\textbf{Swap of parameters in torsion quotients}]\label{cor:swap_parameters_torsion_quotient}
The following torsion quotients are reflection isomorphic:
\begin{itemize}
\item $\B_*^*(n,m;K,b,c)$ and $\B_*^*(m,n;K,c,b)$
\item $\B^*(n,m;K,c)$ and $\B_*(m,n;K,c)$. 
\item $\B(n,m;K)$ and $\B(m,n;K)$.
\end{itemize}
\end{corollary}

\begin{proof}
The argument is precisely the same as that of Corollary \ref{cor:permutation_parameters_torsion_quotient}, using the reflection isomorphism between $J(K,bn',cm')$ and $J(K,cm',bn')$ given by Proposition \ref{prop:permutation_of_parameters_generalized_j_groups}.
\end{proof}

Since we can take $K=(\infty,\ldots,\infty)$, and $b=c=\infty$ if needed, Corollary \ref{cor:embedding_torsion_quotients} implies in turn the following result.

\begin{corollary}\label{cor:embedding_other_j-braid}
Let $\varphi:\B_*^*(pn,pm)\to J((p+2)\cdot \infty)$ be the embedding of Theorem \ref{theo:embedding_j-braid}.
\begin{enumerate}
\item The morphism $\varphi$ induces an embedding $\B^*(pn,pm)\to J(p\cdot\infty,n,\infty)$.
\item The morphism $\varphi$ induces an embedding $\B_*(pn,pm)\to J((p+1)\cdot \infty,m)$.
\item The morphism $\varphi$ induces an embedding $\B(pn,pm)\to J(p\cdot \infty,n,m)$.
\end{enumerate}
Moreover, for each of these morphisms, the image is a finite index normal subgroup.
\end{corollary}

\subsection{Proof of the embedding theorem}
In this section we fix two positive integers $n,m$, and we let $d:=n\wedge m$ denote the gcd of $n$ and $m$. We write $m=qn+r$ the euclidean division of $m$ by $n$. We also fix a $d$-tuple $K:=(k_1,\ldots,k_d)$ of elements in $\N_{\geqslant 2}\cup \{\infty\}$, along with two elements $b,c\in \N_{\geqslant 2}\cup\{\infty\}$. 

This section is devoted to the proof of Theorem \ref{theo:embedding_j-braid}. The core of the proof is to compute a presentation of the normal closure of $\{s_1,\ldots,s_p,s_{p+1}^{n},s_{p+2}^{m}\}$ in $J(K,bn,cm)$ using the Reidemeister-Schreier method. The proof is separated in several intermediate results, and the first reduction is made possible by the following elementary group theoretic result. 

\begin{lemma}\label{lem:poufpouf}
Let $G$ be a group, and let $H$ be a normal subgroup of $G$. Assume that the natural map $Z(G)\to G/H$ is surjective. Then for every $x\in H$, the conjugacy classes of $x$ in $H$ and in $G$ coincide.
\end{lemma}
\begin{proof}
Let $x\in H$, and let $\mathrm{Cl}_H(x)$ (resp. $\mathrm{Cl}_G(x)$) denote the conjugacy class of $x$ in $H$ (resp. in $G$). Since $H$ is a subgroup of $G$, it is immediate that $\mathrm{Cl}_H(x)\subset \mathrm{Cl}_G(x)$. Conversely, let $gxg^{-1}\in \mathrm{Cl}_G(x)$. By assumption, we can write $g=zh$ with $z\in Z(G)$ and $h\in H$. We then have $gxg^{-1}=hxh^{-1}\in \mathrm{Cl}_H(x)$, which terminates the proof.
\end{proof}

\begin{lemma}\label{lem:normal_closure_is_normal_closure}
Let $N$ be the normal closure of $\{s_1,\ldots,s_{p+1},s_{p+2}^{m}\}$ in $J(K,bn,cm)$. The normal closures of $\{s_1,\ldots,s_p,s_{p+1}^{n},s_{p+2}^{m}\}$ in $J(K,bn,cm)$ and in $N$ are equal.
\end{lemma}

\begin{proof}
The element $\alpha:=s_1\cdots s_{p+2}$ of $J(K,bn,cm)$ is central. Its image in the quotient $J(K,bn,cm)/N\simeq \Z/m\Z$ is $1$, which is a generator. By Lemma \ref{lem:poufpouf}, the conjugacy classes of an element of $N$ in $N$ or in $J(K,bn,cm)$ are equal.

Now, the normal closure of a finite subset is generated by the union of the conjugacy classes of its elements. Since the union of the conjugacy classes of $\{s_1,\ldots,s_p,$ $s_{p+1}^{n},s_{p+2}^{m}\}$ in $J(K,bn,cm)$ or in $N$ are equal, this finishes the proof.
\end{proof}

Using this lemma, we can first compute a presentation of the group $N$ using the Reidemeister-Schreier method, and then compute a presentation for the normal closure of $\{s_1,\ldots,s_p,s_{p+1}^{n},s_{p+2}^{m}\}$ in $N$, again using the Reidemeister-Schreier method.

\begin{proposition}\label{NormalClosureBraid1}
The correspondance
\[\begin{cases} x_i\mapsto s_i &\text{for }i\in \intv{1,p+1},\\z\mapsto s_{p+2}^{m}.\end{cases}\]
induces an injective morphism from $\B^*((p+1),(p+1)m;(K,bn),c)$ to $J(K,bn,cm)$, whose image is the normal closure of $\{s_1,\ldots,s_{p+1},s_{p+2}^{m}\}$ in $J(K,bn,cm)$. 
\end{proposition}
\begin{proof}
To simplify the notations in this proof, we denote $bn$ by $k_{p+1}$. Let $N$ denote the normal closure of $\{s_1,\ldots,s_{p+1},s_{p+2}^{m}\}$ in $J(K,k_{p+1},cm)$. The quotient group \\$J(K,bn,cm)/N$ is isomorphic to $\Z/m\Z$ and a Schreier transversal for $N$ is given by $\{s_{p+2}^j\}_{j\in \intv{0,m-1}}$. By the Reidemeister-Schreier algorithm, the group $N$ admits the following group presentation: 
\begin{subequations}
\label{NormalClosureEQ1} 
\begin{align}
&(1) \,\, \mathrm{Generators}\!:\,  \{x_{i,j}\}_{(i,j)\in \intv{1,p+2}\times\llbracket 0,m-1\rrbracket};\notag\\
&(2) \,\, \mathrm{Relations}\!: \,\notag\\
&x_{i,j}^{k_i} \, \text{ for all $i$ such that $k_i<\infty$ and $0\leqslant j\leqslant m-1$}, \label{NormalClosureEQ1:-1}\\
& x_{p+2,j}x_{p+2,j+1}\cdots x_{p+2,j+cm-1}=1 \text{ if $c<\infty$,} \label{NormalClosureEQ1:0}\\
  &  x_{p+2,j}=1 \text{ for all $0\leqslant j \leqslant m-2$},\label{NormalClosureEQ1:1}\\
       &  x_{i,j}x_{i+1,j}\cdots x_{p+2,j}x_{1,j+1}\dots x_{i-1,j+1}=x_{i+1,j}x_{i+2,j}\cdots x_{p+2,j}x_{1,j+1}\cdots x_{i,j+1}\,(*),\label{NormalClosureEQ1:2}\\
      &  (*) \text{ for all $1\leqslant i< p+2$, \, $0\leqslant j\leqslant m-1$.}\notag
\end{align}
\end{subequations}
where the index $i$ is taken modulo $p+2$, and the index $j$ is taken modulo $m-1$. Moreover, for all $i\in \intv{1,p+1}$ $j\in \intv{0,m-1}$, we have $x_{i,j}=s_{p+2}^js_is_{p+2}^{-j}$, while $x_{p+2,m-1}=s_{p+2}^{m}$.
Now, for $j\in \intv{0,m-2}$, we have $x_{p+2,j}=1$ so that Equation \eqref{NormalClosureEQ1:0} reads
\begin{equation}\label{NormalClosureEQ1:3}
x_{p+2,m-1}^c=1 \text{ if $c<\infty$}
\end{equation}
and 
Equation \eqref{NormalClosureEQ1:2} reads 
\begin{equation}\label{NormalClosureEQ1:4}
x_{i,j}x_{i+1,j}\cdots x_{p+1,j}x_{1,j+1}\dots x_{i-1,j+1}=x_{i+1,j}x_{i+2,j}\cdots x_{p+1,j}x_{1,j+1}\cdots x_{i,j+1}.
\end{equation}

For any $l\in \intv{1,m(p+1)}$, there is a unique way to write $l$ as $j(p+1)+i$ with $j\in\intv{0,m-1}$ and $i\in \intv{1,p+1}$. Using this, we relabel the generators by setting $y_{j(p+1)+i}:=x_{i,j}$ $i\in \intv{0,m-1}$ and $i\in \intv{1,p+1}$. We also set $z:=x_{p+2,m-1}$. The group $N$ then admits the presentation 

\begin{subequations}
\label{NormalClosureEQ2} 
\begin{align}
&(1) \,\, \mathrm{Generators}\!:\,  \{y_1,\ldots,y_{m(p+1)},z\};\notag\\
&(2) \,\, \mathrm{Relations}\!: \,\notag\\
& y_l^{k_i}=1 \text{ if $l=j(p+1)+i$ and $k_i<\infty$,} \label{NormalClosureEQ2:-1} \\
& z^c=1 \text{ if $c<\infty$}, \label{NormalClosureEQ2:0}\\
       &  y_iy_{i+1}\cdots y_{i+p}=y_{i+1}y_{i+2}\cdots y_{i+p+1}\text{ for all $1\leqslant i< (m-1)(p+1)+1$},\label{NormalClosureEQ2:1}\\
      &  y_{(m-1)(p+1)+i}\cdots y_{m(p+1)}zy_1\cdots y_{i-1}=y_{(m-1)(p+1)+i+1}\cdots y_{m(p+1)}zy_1\cdots y_{i}\ (*).\label{NormalClosureEQ2:2}\\
      &(*) \text{ for all $1\leqslant i<p+2$}\notag
\end{align}
\end{subequations}
\noindent
Rephrasing Equation \eqref{NormalClosureEQ2:1} yields that, for $i\in \intv{1,(m-1)(p+1)}$, we have 
\[y_{i+p+1}=y_{i}^{y_{i+1}\cdots y_{i+p}}.\]
Since $y_i$ commutes with itself, we deduce that
\[y_{i+p+1}=y_i^{y_iy_{i+1}\cdots y_{i+p}}=y_i^{y_1\cdots y_{p+1}}.\]
Let us define $\delta:=y_1\cdots y_{p+1}$. By an immediate induction on the above formula, we obtain that
\begin{equation}\label{NormalClosureEQ2:3}
\forall j\in \intv{0,m-1},i\in\intv{1,p+1},~y_{j(p+1)+i}=y_i^{(\delta^j)}.
\end{equation}
Using this, we can delete the generators $y_l$ for $l>p+1$ from the presentations, along with Equation \eqref{NormalClosureEQ2:1}.

Now, for $i\in \intv{1,p+1}$, we have
\[y_{(m-1)(p+1)+i}\cdots y_{m(p+1)}=(y_i\cdots y_{p+1})^{(\delta^{m-1})}.\]
In particular, Equation \eqref{NormalClosureEQ2:2} is equivalent to 
\[\delta^{1-m}y_i\cdots y_{p+1}\delta^{m-1}zy_1\cdots y_{i-1}=\delta^{1-m}y_{i+1}\cdots y_{p+1}\delta^{m-1}zy_1\cdots y_i.\]
By simplifying this equation by $\delta^{1-m}$, we obtain that
\begin{equation}\label{NormalClosureEQ2:4}\forall i\in \intv{1,p+1}, y_i\cdots y_{p+1}\delta^{m-1}zy_1\cdots y_{i-1}=y_{i+1}\cdots y_{p+1}\delta^{m-1}zy_1\cdots y_i.\end{equation}
By an immediate induction, we obtain that
\[y_1\cdots y_{p+1}\delta^{m-1}z=\delta^{m-1}zy_1\cdots y_{p+1}\Leftrightarrow z\delta=\delta z\]
By adding the relation $z\delta=\delta z$ to our presentation of $N$, we can ask \eqref{NormalClosureEQ2:4} to be verified only for $i\in \intv{1,p}$, and not for $i\in \intv{1,p+1}$.

We finally obtained that the group $N$ admits the following presentation
\begin{subequations}
\begin{align}
&(1) \,\, \mathrm{Generators}\!:\,  \{y_1,\ldots,y_{p+1},z\};\notag\\
&(2) \,\, \mathrm{Relations}\!: \,\notag\\
& y_i^{k_i} \text{ for all $i$ such that $k_i<\infty$}\\
& z^c=1 \text{ if $c<\infty$,}\\
		&  y_1\cdots y_{p+1}z=zy_1\cdots y_{p+1},\\
       &  y_iy_{i+1}\cdots y_{p+1}z\delta^{m-1}y_1\cdots y_{i-1}=y_{i+1}\cdots y_{p+1}z\delta^{m-1}zy_1\cdots y_i\text{ for all $1\leqslant i\leqslant p$},
\end{align}
\end{subequations}
where $y_i=x_{i,0}=s_i$ and where $z=x_{p+2,m-1}=s_{p+2}^{m}$. As this presentation is precisely the defining presentation of $\B^*((p+1),(p+1)m;(K,bn),c)$, we have the desired result.
\end{proof}

\begin{remark}Since $\B^*(p+1,(p+1)m,(K,bn),c)$ is generated by $x_1,\ldots,x_{p+1}$ and $z$, the above proposition also proves that the subgroup $\langle s_1\cdots,s_{p+1},s_{p+2}^{m}\rangle$ in $J(K,bn,cm)$ is normal. \end{remark}

\begin{remark}\label{SelfIso1}
By Theorem \ref{TheoremBraidIsoIntroCircularFR}, the normal closure of $\{s_1,\dots,s_{p+1},s_{p+2}^{m}\}$ in $J((p+2)\cdot \infty)$ is in turn isomorphic to $J((p+2)\cdot\infty)$.
\end{remark}

\begin{proof}[Proof of Theorem \ref{theo:embedding_j-braid}]
We aim to compute a presentation of the normal closure of $\{s_1,\ldots,s_{p},s_{p+1}^{n},s_{p+2}^{m}\}$ in $J(K,bn,cm)$. By Lemma \ref{lem:normal_closure_is_normal_closure}, it is sufficient to compute a presentation of the normal closure of $\{s_1,\ldots,s_{p+1}^{n},s_{p+2}^{m}\}$ in the subgroup $N$ of $J(K,bn,cm)$ which is the normal closure of $\{s_1,\ldots,s_p,s_{p+1},s_{p+2}^{m}\}$. By Proposition \ref{NormalClosureBraid1}, we can identify $N$ with $\B^*(p+1,(p+1)m;(K,bn),c)$.
The euclidean division of $(p+1)m$ by $p+1$ is $(p+1)m=(p+1)\times  m+0$. The group $\B^*(p+1,(p+1)m;(K,bn),c)$ is then defined by the following presentation:
\begin{subequations}
\begin{align}
&(1) \,\, \mathrm{Generators}\!:\,  \{x_1,\dots,x_{p+1},z\};\notag\\
&(2) \,\, \mathrm{Relations}\!: \,\notag\\
 & x_i^{k_i}=1 \text{ if $k_i<\infty$,} \\
 & z^c=1 \text{ if $c<\infty$,}\\
  &  x_1\cdots x_{p+1}z=zx_1\cdots x_{p+1},\\
       & x_{i+1}\cdots x_{p+1}z\delta^{m-1}x_1\cdots x_{i}=x_i\cdots x_{p+1}z\delta^{m-1}x_1\cdots x_{i-1}, \, \forall 1\leqslant i \leqslant p+1,
\end{align}
\end{subequations}
where $k_{p+1}$ denotes $bn$.

Under the identification $N\simeq \B^*(p+1,(p+1)m;(K,bn),c)$, the normal closure we have to compute is the normal closure of $\{x_1,\ldots,x_p,x_{p+1}^{n},z\}$ in $\B^*(p+1,(p+1)m;(K,bn),c)$. Let us denote by $H$ this normal closure.
We have $\B^*(p+1,(p+1)m;(K,bn),c)/H\simeq \Z/n\Z$ and a Schreier transversal for $H$ is $\{x_{p+1}^j\}_{j\in \intv{0,n-1}}$. By the Reidemeister-Schreier algorithm, the group $H$ is generated by the elements $\{x_{i,j}\}_{i\in \intv{1,p+1},j\in \intv{0,n-1}}$ along with the elements $\{z_j\}_{j\in \intv{0,n-1}}$, where 
\begin{itemize}
\item $x_{i,j}=x_{p+1}^{j}x_ix_{p+1}^{-j}$ for $i<p+1$, 
\item $x_{p+1,j}=1$ for $j<n-1$,
\item $x_{p+1,n-1}=x_{p+1}^{n}$,
\item $z_j=x_{p+1}^j z x_{p+1}^{-j}$.
\end{itemize}
In order for the relations of the presentation obtained by the Reidemeister-Schreier method to be readable, we introduce some intermediate elements.
\begin{itemize}
\item We set $y:=x_{p+1,n-1}$.
\item For $j\in \intv{0,n-1}$, we set $\delta_j=x_{1,j}\cdots x_{p,j}$ and $\delta'_j=x_{1,j}\cdots x_{p+1,j}$. Notice that $\delta_j=\delta'_j$ for $j<n-1$, and that $\delta'_{n-1}=\delta_{n-1}y$.
\item We set $\Delta=\delta_0\cdots \delta_{n-1}y$. Notice that, taking indices modulo $n$ and setting $D_j:=\delta'_{j+1}\cdots\delta'_{j+m-1}$, we have
\[D_j=\begin{cases} 
\Delta^q\delta_0\cdots \delta_{r-2}&\text{if }j=n-1,\\
\delta_{j+1}\cdots \delta_{n-1}y\Delta^{q}\delta_0\cdots \delta_{r+j-n-1}&\text{if }j\in \intv{n-r,n-2},\\
\delta_{j+1}\cdots \delta_{n-1}y\Delta^{q-1}\delta_0\cdots \delta_{r+j-1}&\text{if }j\in \intv{0, n-r-1},
\end{cases}\]
where $m=qn+r$ is the euclidean division of $m$ by $n$.
\item For $i\in \intv{1,p+1}$ and $j\in \intv{0,n-1}$, we define $P_{i,j}$ as
\[\begin{cases}
x_{i,n-1}\cdots x_{p,n-1}x_{p+1,n-1}z_0D_{n-1}x_{1,r-1}\cdots x_{i-1,r-1}&\text{if }j=n-1,\\
x_{i,j}\cdots x_{p+1,j}z_{j+1}D_jx_{1,r+j-n}\cdots x_{i-1,r+j-n} &\text{if }j\in \intv{n-r,n-2},\\
x_{i,j}\cdots x_{p+1,j}z_{j+1}D_jx_{1,r+j}\cdots x_{i-1,r+j}&\text{if }j\in \intv{0,n-r-1}.
\end{cases}\]
\end{itemize}
With these definitions, the relations for of $H$ obtained via the Reidemeister-Schreier method are
\begin{subequations}\label{reid-schrei_penible1}
\begin{align}
 & x_{i,j}^{k_i}=1 \text{ if $k_i<\infty$},\, \forall 1\leqslant i\leqslant p,\, 0\leqslant j\leqslant n-1 \label{reid-schrei_penible1:-2}\\
 & y^b=1 \text{ if $b<\infty$}, \label{reid-schrei_penible1:-1}\\
 & z_j^c=1 \text{ if $c<\infty$}, \, \forall 0\leqslant j \leqslant n-1 \label{reid-schrei_penible1:0}\\
  &  \delta_j z_{j+1}=z_j\delta_j, \, \forall 1\leqslant j\leqslant n-2\label{reid-schrei_penible1:1}\\
   &  \delta_{n-1} yz_{0}=z_{n-1}\delta_{n-1}y, \label{reid-schrei_penible1:2}\\
       & P_{i+1,j}=P_{i,j}, \, \forall 1\leqslant i \leqslant p, 0\leqslant j\leqslant n-1 \label{reid-schrei_penible1:3}
\end{align}
\end{subequations}
Equation \eqref{reid-schrei_penible1:1} implies that, for $j\in \intv{1,n-2}$, we have $z_{j+1}=z_j^{\delta_j}=z_0^{\delta_0\cdots \delta_j}$. Replacing $z_{n-1}$ with $z_0^{\delta_0\cdots \delta_{n-2}}$ in Equation \eqref{reid-schrei_penible1:2} yields
\[\delta_{n-1}yz_0=(\delta_0\cdots \delta_{n-2})^{-1}z_0\delta_0\cdots \delta_{n-1}y\Leftrightarrow \Delta z_0=z_0\Delta.\]

For any $l\in\intv{1,pn}$, there is a unique way to write $l$ as $jp+i$ with $j\in \intv{0,n-1}$ and $i\in \intv{1,p}$. Using this, we relabel the generators $x_{i,j}$ as $a_{jp+i}$ for $i\in \intv{1,p}$ and $j\in \intv{0,n-1}$. 

Using this relabel, we see that $\delta_j=a_{jp+1}\cdots a_{jp+p}$. In particular $\Delta=a_1\cdots a_{pn} y$. Moreover, we can rewrite the product $D_j$ as follows.
\[D_j=\begin{cases} 
y\Delta^qa_1\cdots a_{(r-1)p}&\text{if }j=n-1,\\
a_{(j+1)p+1}\cdots a_{pn}y\Delta^q a_1\cdots a_{(r+j-n)p}&\text{if }j\in \intv{n-r,n-2},\\
a_{(j+1)p+1}\cdots a_{pn}y\Delta^{q-1} a_1\cdots a_{(r+j)p}&\text{if }j\in \intv{0,n-r-1}.
\end{cases}\]
Notice that for $j\in \intv{0,n-2}$, we have.
\begin{align*}
z_{j+1}a_{(j+1)p+1}\cdots a_{pn}y &=z_{j+1}\delta_{j+1}\cdots \delta_{n-1}y\\
&=\delta_{j+1}\cdots \delta_{n-1}yz_0\\
&=a_{(j+1)p+1}\cdots a_{pn}yz_0.
\end{align*}
Using these formulae, we can rewrite the product $P_{i,j}$ as
\[P_{i,j}=\begin{cases} a_{p(n-1)+i}\cdots a_{pn}yz_0\Delta^qa_1\cdots a_{(r-1)p+i-1}&\text{if }j=n-1,\\
a_{jp+i}\cdots a_{pn}yz_0\Delta^q a_1\cdots a_{(r+j-n)p+i-1}&\text{if }j\in \intv{n-r,n-2},\\
a_{jp+i}\cdots a_{pn} yz_0\Delta^{q-1}a_1\cdots a_{(r+j)p+i-1}&\text{if }j\in \intv{0,n-r-1}.\end{cases}\]
We notice that the expressions of $P_{i,j}$ for $j=n-1$ and for $j\in \intv{n-r,n-2}$ are equal. We then have
\[P_{i,j}=\begin{cases}
a_{l}\cdots a_{pn}yz_0\Delta^q a_1\cdots a_{pr+l-pn-1}&\text{if }j\in \intv{n-r,n-1}, l=jp+i,\\
a_{l}\cdots a_{pn} yz_0\Delta^{q-1}a_1\cdots a_{pr+l-1}&\text{if }j\in \intv{0,n-r-1}, l=jp+i.\end{cases}\]

Since $i\in \intv{1,p}$, we have $j\in \intv{0,n-r-1}$ if and only if $l=jp+i\in \intv{1,pn-pr}$. Likewise, we have $j\in \intv{n-r,n-1}$ if and only if $l=jp+i\in \intv{pn-pr+1,pn}$.

Finally, we can rewrite the relations of \eqref{reid-schrei_penible1} to obtain the following presentation of $H$.
\begin{subequations}
\begin{align}
&(1) \,\, \mathrm{Generators}\!:\,  \{a_1,\ldots,a_{pn},y,z_0\};\notag\\
&(2) \,\, \mathrm{Relations}\!: \,\notag\\
& a_l^{k_i}=1 \text{ if $l=jp+i$ and $k_i<\infty$,}\\
& y^b=1 \text{ if $b<\infty$},\\
& z_0^c=1 \text { if $c<\infty$},\\
		&  a_1\cdots a_{pn}y\cdots z_0=z_0a_1\cdots a_{pn}y,\\
       &  a_{l+1}\cdots a_{pn}yz_0\Delta^q a_1\cdots a_{pr+l-pn}=a_{l}\cdots a_{pn}yz_0\Delta^q a_1\cdots a_{pr+l-pn-1}\ (*)\\& (*)\text{ for $l\in \intv{pn-pr+1,pn}$},\notag\\
       & a_{l+1}\cdots a_{pn} yz_0\Delta^{q-1}a_1\cdots a_{pr+l}=a_{l}\cdots a_{pn} yz_0\Delta^{q-1}a_1\cdots a_{pr+l-1}\ (**)\\&(**) \text{ for $l\in \intv{1,pn-pr}$}.\notag
\end{align}
\end{subequations}
where we have $a_{jp+i}=x_{i,j}=x_{p+1}^jx_ix_{p+1}^{-j}$, $y=x_{p+1}^{n}$ and $z_0=z$. Identifying $\B_*(p+1,(p+1)m;(K,bn),c)$ with its image in $J(K,bn,cm)$ using Proposition \ref{NormalClosureBraid1}, we identify $a_{jp+i}$ with $s_{p+1}^js_is_{p+1}^{-j}$, $y$ with $s_{p+1}^{n}$ and $z$ with $s_{p+2}^{m}$. As the  above presentation is precisely the defining presentation of $\B^*_*(pn,pm;d\cdot K,b,c)$, we have the desired result.
\end{proof}

\section{Generalized $J$-groups as torsion quotients of $J$-braid groups}\label{sec:gen_j-groups_as_torsion_quotients}
In the last section, we saw with Corollary \ref{cor:embedding_torsion_quotients} (embedding of torsion quotients) 
that torsion quotients of $J$-braid groups can be seen as (finite-type) generalized $J$-groups. It turns out that the converse is also true, and that the family of finite type generalized $J$-groups coincides with that of torsion quotients of $J$-braid groups.

For readability purposes, we extend the definition of reflection isomorphism by saying that a torsion quotient of $J$-braid group $W$ and a generalized $J$-group $Q$ can be reflection isomorphic if there is an isomorphism $\varphi:W\to Q$ such that $\varphi(R(W))=R(Q)$. What we are going to show is that every finite-type generalized $J$-group is reflection isomorphic to some torsion quotient of $J$-braid group. The proof is split in several intermediate results and relies ultimately on a case-by-case approach.

For the remainder of the section, we fix an integer $p$, along with elements $k_1,\ldots,k_p$, $a,b,c$ in $\N_{\geqslant 2}\cup \{\infty\}$.

We begin by giving a list of finite-type generalized $J$-groups up to reflection isomorphism.

\begin{lemma}\label{lem:finite_type_generalized_list}A generalized $J$-group has finite-type if and only if it is reflection isomorphic to a generalized $J$-group of one of the following form:
\begin{itemize}[itemsep=4pt]
\item $J\jg{k_1&\cdots&k_p&nb&mc\\1&\cdots&1&n&m}$ with $n,m\geqslant 1$,
\item $J\jg{k_1&\cdots&k_p&2a&2b&lc\\1&\cdots&1&2&2&l}$ with $l\geqslant 2$,
\item $J\jg{k_1&\cdots&k_p&2a&3b&3c\\1&\cdots&1&2&3&3}$,
\item $J\jg{k_1&\cdots&k_p&2a&3b&4c\\1&\cdots&1&2&3&4}$,
\item $J\jg{k_1&\cdots&k_p&2a&3b&5c\\1&\cdots&1&2&3&5}$.
\end{itemize}
We will denote these families of groups by $F(n,m)$, $F(2,2,l)$, $F(2,3,3)$, $F(2,3,4)$, $F(2,3,5)$ respectively..
\end{lemma}
\begin{proof}
By definition, a generalized $J$-group $J\binom K{K'}$ has finite-type if and only if $J(K)/J\binom K{K'}$ is finite. This quotient is isomorphic to $J(\overline{K'})$, where $\overline{K'}$ denotes the tuple consisting of coordinates of $K'$ different than $1$. If $\overline{K'}$ has length 2, then $J(K')$ is a finite abelian group since $K'$ contains only finite elements. If $\overline{K'}$ has length $3$ or more, then we obtain the result by Lemma \ref{lem:finite_index_J_or_coxeter}.
\end{proof}

Now that these cases have been listed, we can give a precise statement
\begin{theorem}[\textbf{Generalized $J$-groups as torsion quotients}]\label{theo:generalized_j_groups_torsion_quotients}
We have the following reflection isomorphisms:
\begin{itemize}[itemsep=4pt]
\item $J\jg{k_1&\cdots&k_p&nb&mc\\1&\cdots&1&n&m}\cong_{\mathrm{Ref}}\B_*^*(pn,pm,d\cdot(k_1,\ldots,k_p),b,c)$, where $n,m\geqslant 1$,
\item $J\jg{k_1&\cdots&k_p&2a&2b&lc\\1&\cdots&1&2&2&l}\cong_{\mathrm{Ref}}\B(2(lp+l+1),2(lp+l+1);(M,l\cdot a,l\cdot b,2\cdot c))$, where $l\geqslant 2$ and where $M:=2l\cdot (k_1,\ldots,k_p)$,
\item $J\jg{k_1&\cdots&k_p&2a&3b&3c\\1&\cdots&1&2&3&3}\cong_{\mathrm{Ref}}\B(12p+14,12p+14;(M,4\cdot b,4\cdot c,6\cdot a))$, where $M:=12\cdot (k_1,\ldots,k_p)$.
\item $J\jg{k_1&\cdots&k_p&2a&3b&4c\\1&\cdots&1&2&3&4}\cong_{\mathrm{Ref}}\B(24p+26,24p+26;(M,12\cdot a,8\cdot b,6\cdot c))$, where $M:=24\cdot (k_1,\ldots,k_p)$.
\item $J\jg{k_1&\cdots&k_p&2a&3b&5c\\1&\cdots&1&2&3&5}\cong_{\mathrm{Ref}}\B(60p+62,60p+62;(M,30\cdot a,20\cdot b,12\cdot c))$, where $M:=60\cdot(k_1,\ldots,k_p)$.
\end{itemize}
In particular, any finite-type generalized $J$-group is reflection isomorphic to a torsion quotient of $J$-braid group.
\end{theorem}

\begin{remark}\label{rem:finite-type_j-groups_generated_by_reflections}
By construction, torsion quotients of $J$-braid groups are generated by a finite number of reflections. Theorem \ref{theo:generalized_j_groups_torsion_quotients} implies that any finite-type generalized $J$-group is generated by a finite number of reflections.
\end{remark}

For the remainder of the section, we do a case-by-case study of the families described in Lemma \ref{lem:finite_type_generalized_list} to prove Theorem \ref{theo:generalized_j_groups_torsion_quotients}. Family $F(n,m)$ has already been dealt with in Corollary \ref{cor:embedding_torsion_quotients}. We separate the other cases in two subsections.

\subsection{The families $F(2,2,l)$, $F(2,3,3)$ and $F(2,3,4)$}
In fact, Corollary \ref{cor:embedding_torsion_quotients} can also help us deal with families $F(2,2,l)$, $F(2,3,3)$ and $F(2,3,4)$. We begin with an elementary lemma:

\begin{lemma}\label{lem:b**(p,p)=J}
The correspondence
\[\begin{cases} x_i\mapsto s_i&\text{for }i\in \intv{1,p},\\ y\mapsto s_{p+1},\\z\mapsto s_{p+2},\end{cases}\]
induces a reflection isomorphism $\B_*^*(p,p;K,b,c)\cong_{\mathrm{Ref}}\B(p+2,p+2;(K,b,c))\cong_{\mathrm{Ref}}J(K,b,c)$.
\end{lemma}
\begin{proof}
The euclidean division of $p$ by itself reads $p=1\times p+0$. The defining presentation of $\B_*^*(p,p;K,b,c)$ is then
\begin{subequations}
\begin{align*}
&\mathrm{Generators}\!:\,  \{x_1,\dots,x_p,y,z\};\\
&\mathrm{Relations}\!: \,\\
 &x_1^{k_1}=x_2^{k_2}=\ldots=x_p^{k_p}=y^b=z^c=1,\\
  &  x_1\cdots x_nyz=zx_1\cdots x_ny,\\
  & x_{i+1}\cdots x_nyzx_1\cdots x_{i+r}=x_i\cdots x_nyzx_1\cdots x_{i+r-1}, \, \forall 1\leqslant i \leqslant n,
\end{align*}
\end{subequations}
This is exactly the defining presentation of $\B(p+2,p+2;(K,b,c))$ and of $J(K,b,c)$.
\end{proof}

We now prove a technical lemma, which we then apply to the families $F(2,2,l)$, $F(2,3,3)$ and $F(2,3,4)$.

\begin{lemma}\label{lem:pirouette_avec_poufpoufpouf}
Let $k,d,n,m$ be positive integers. We have a reflection isomorphism
\[J\jg{k_1&\cdots&k_p&ak&bdn&cdm\\1&\cdots &1&k&dn&dm}\cong_{\mathrm{Ref}}J\jg{d\cdot(k_1&\cdots&k_p&ak)&bn&cm\\d\cdot (1&\ldots&1&k) &n&m}\]
\end{lemma}
\begin{proof}
Consider the group $J_1:=J(d\cdot(k_1,\ldots,k_p,ak),bn,cm)$, and let us denote its generators by $\sigma_1,\ldots,\sigma_{d(p+1)},\tau,\mu$. Recall that every integer $l \in \intv{1,(p+1)d}$ can be written uniquely as $l=j(p+1)+i$ with $j\in \intv{0,d-1}$ and $i\in \intv{1,p+1}$. Using this, we relabel $\sigma_{l}$ as $\sigma_{i,j}$. 

By Lemma \ref{lem:b**(p,p)=J}, the correspondence
\[\begin{cases} \sigma_{i,j}\mapsto x_{j(p+1)+i}&\text{for }j\in \intv{0,d-1},i\in \intv{1,p+1},\\ \tau\mapsto y,\\\mu\mapsto z,\end{cases}\]
induces a reflection isomorphism 
\[J_1\cong_{\mathrm{Ref}} W:=\B_*^*(d(p+1),d(p+1);d\cdot (k_1,\ldots,k_p,ak),bn,cm).\]
Now, Theorem \ref{theo:embedding_j-braid} (embedding of $J$-braid groups) gives an explicit reflection isomorphism between $W$ and the generalized $J$-group $J_2:=J\jg{k_1&\cdots&k_p&ak&bdn&cdm\\ 1&\cdots&1&1&d&d}$. Let $s_1,\ldots,s_{p+3}$ denote the generators of $J(k_1,\ldots,k_p,ak,bdn,cdm)$. Composing the two above reflection isomorphism  we obtain that the correspondence
\[\begin{cases} \sigma_{i,j}\mapsto s_{p+2}^j s_is_{p+2}^{-j}&\text{for }j\in \intv{0,d-1},i\in \intv{1,p+1},\\\tau\mapsto s_{p+2}^d,\\ \mu\mapsto s_{p+3}^d,\end{cases}\]
induces a reflection isomorphism $\varphi$ between $J_1$ and $J_2$.

Now, the generalized $J$-group $J_1':=J\jg{d\cdot(k_1&\cdots&k_p&ak)&bn&cm\\d\cdot (1&\ldots&1&k) &n&m}$ is defined as the normal closure in $J_1$ of the set
\[X:=\{\sigma_{i,j}~|~j\in \intv{0,d-1},i\in \intv{1,p}\}\cup \{\sigma_{p+1,j}^k~|~j\in \intv{0,d-1}\}\cup \{\tau^n,\mu^m\}.\]
The image of $J_1'$ under $\varphi$ is the normal closure $J_2'$ of $\varphi(X)$ in $J_2$. 

Let $\widetilde{J_2}$ denote the generalized parent $J$-group of $J_2$. The quotient $\widetilde{J_2}/J_2$ is isomorphic to $\Z/d\Z\times \Z/d\Z$, and is generated by the respective images of $s_{p+2}$ and $s_{p+3}$. Moreover, the image of the center of $\widetilde{J_2}$ in $\widetilde{J_2}/J_2$ is $(1,1)$. The quotient $Q$ of $\widetilde{J_2}/J_2$ by this image is then a cyclic group of order $d$, generated by the image of $s_{p+2}$. It is also generated by the image of $s_{p+3}$. Applying Proposition \ref{prop:conj_reflections_j-groupe} (conjugacy classes of reflections), we obtain that 
\begin{itemize} 
\item For all $i\in \intv{1,p+1}$ and $q\geqslant 1$ the set $C_{i,q}:=\{s_{p+2}^j (s_i^q)s_{p+2}^{-j}~|~j\in \intv{0,d-1}\}$ is a complete set of representatives of the conjugacy classes in $J_2$ of elements conjugate to $(s_i)^q$ in $\widetilde{J_2}$.
\item For $g\geqslant 1$, an element of $J_2$ which is conjugate to $(s_{p+2})^q$ (resp. to $(s_{p+3})^q)$) in $\widetilde{J_2}$ is also conjugate to $(s_{p+2})^q$ (resp. to $(s_{p+3})^q$) in $J_2$.
\end{itemize}
By definition of $\varphi$, we have
\[\varphi(X)=(\bigcup_{i=1}^p C_{i,1})\cup C_{p+1,k}\cup \{s_{p+2}^{dn},s_{p+3}^dm\}\]
After what is said above, the normal closure of this set in $J_2$ is the normal closure of $\{s_1,\ldots,s_p,s_{p+1}^k,s_{p+2}^{dn},s_{p+3}^{dm}\}$ in $\widetilde{J_2}$. Since $\widetilde{J_2}=J(k_1,\ldots,k_p,ak,bdn,cdm)$, the reflection isomorphism $J_1\cong_{\mathrm{Ref}}J_2$ restricts to a reflection isomorphism
\[J\jg{d\cdot(k_1&\cdots&k_p&ak)&bn&cm\\d\cdot (1&\ldots&1&k) &n&m}\cong_{\mathrm{Ref}}J\jg{k_1&\cdots&k_p&ak&bdn&cdm\\1&\cdots &1&k&dn&dm}\]
as we wanted to show.
\end{proof}

\begin{lemma}[\textbf{Family $F(2,2,l)$}]\label{lem:generalized_j_groups_torsion_quotients_case_2}
Let $M:=2l\cdot(k_1,\ldots,k_p)$. We have a reflection isomorphism
\[J\jg{k_1&\cdots&k_p&2a&2b&lc\\1&\cdots&1&2&2&l}\cong_{\mathrm{Ref}} \B_*^*(2(p+1)l,2(p+1)l;(M,l\cdot a,l\cdot b),c,c).\]
Moreover, the latter group is reflection isomorphic to $J(M,l\cdot a,l\cdot b,2\cdot c)$ and we have the second statement in Theorem \ref{theo:generalized_j_groups_torsion_quotients}.
\end{lemma}
\begin{proof}
By Proposition \ref{prop:permutation_of_parameters_generalized_j_groups} (permutation of parameters), we have a reflection isomorphism
\[J\jg{k_1&\cdots&k_p&2a&2b&lc\\1&\cdots&1&2&2&l}\cong_{\mathrm{Ref}}J\jg{k_1&\cdots&k_p&lc&2a&2b\\1&\cdots&1&l&2&2}.\]
We can apply Lemma \ref{lem:pirouette_avec_poufpoufpouf} with $d=2,n=m=1$ and $k=l$. Again by Proposition \ref{prop:permutation_of_parameters_generalized_j_groups}, we obtain reflection isomorphisms
\[J\jg{k_1&\cdots&k_p&lc&2a&2b\\1&\cdots&1&l&2&2}\cong_{\mathrm{Ref}}J\jg{2\cdot(k_1&\cdots&k_p&lc)&a&b\\2\cdot(1&\cdots&1&l)&1&1}\cong_{\mathrm{Ref}}J\jg{2\cdot(k_1&\cdots&k_p)&a&b&cl&cl\\2\cdot(1&\cdots&1)&1&1&l&l}.\]
Using Theorem \ref{theo:embedding_j-braid}, this last group is in turn reflection isomorphic to 
\[\B_*^*(2l(p+1),2l(p+1);l\cdot N,c,c).\]
where $N=(2\cdot(k_1,\ldots,k_p),a,b)$. Again by Proposition \ref{prop:permutation_of_parameters_generalized_j_groups}, this group is reflection isomorphic to 
\[\B_*^*(2(p+1)l,2(p+1)l;(M,l\cdot a,l\cdot b),c,c)\]
since $(M,l\cdot a,l\cdot b)$ is equal to $l\cdot N$ up to permutation.
The second statement is a direct application of Lemma \ref{lem:b**(p,p)=J}.
\end{proof}

\begin{lemma}[\textbf{Family $F(2,3,3)$}]\label{lem:generalized_j_groups_torsion_quotients_case_3}
Let $M:=12\cdot(k_1,\ldots,k_p)$. We have a reflection isomorphism
\[J\jg{k_1&\cdots&k_p&2a&3b&3c\\1&\cdots&1&2&3&3}\cong_{\mathrm{Ref}} \B^*_*(12(p+1),12(p+1);(M,4\cdot b,4\cdot c,4\cdot a),a,a).\]
Moreover, the latter group is reflection isomorphic to $J(M,6\cdot a,4\cdot b,4\cdot c)$ we have the third statement in Theorem \ref{theo:generalized_j_groups_torsion_quotients}.
\end{lemma}
\begin{proof}
We can apply Lemma \ref{lem:pirouette_avec_poufpoufpouf} with $d=3,n=m=1$ and $k=2$. Again by Proposition \ref{prop:permutation_of_parameters_generalized_j_groups}, we obtain reflection isomorphisms
\[J\jg{k_1&\cdots&k_p&2a&3b&3c\\1&\cdots&1&2&3&3}\cong_{\mathrm{Ref}}J\jg{3\cdot(k_1&\cdots&k_p&2a)&b&c\\3\cdot(1&\cdots&1&2)&1&1}\cong_{\mathrm{Ref}}J\jg{3\cdot(k_1&\cdots&k_p)&b&c&2a&2a&2a\\3\cdot(1&\cdots&1)&1&1&2&2&2}.\]
Using Lemma \ref{lem:generalized_j_groups_torsion_quotients_case_2}, this last group is in turn reflection isomorphic to 
\[\B_*^*(4(3p+3),4(3p+3);(N,2.a,2.a),a,a)\]
where $N=4\cdot (3\cdot (k_1,\ldots,k_p),b,c)$. Again by Proposition \ref{prop:permutation_of_parameters_generalized_j_groups}, this group is reflection isomorphic to 
\[\B^*_*(12(p+1),12(p+1);(M,4\cdot b,4\cdot c,4\cdot a),a,a)\]
since $(M,4\cdot b,4\cdot c,4\cdot a)$ is equal to $(N,2\cdot a,2\cdot a)$ up to permutation. 
The second statement is a direct application of Lemma \ref{lem:b**(p,p)=J}.
\end{proof}

\begin{lemma}[\textbf{Family $F(2,3,4)$}]\label{lem:generalized_j_groups_torsion_quotients_case_4}
Let $M:=24\cdot(k_1,\ldots,k_p)$. We have a reflection isomorphism
\[J\jg{k_1&\cdots&k_p&2a&3b&4c\\1&\cdots&1&2&3&4}\cong_{\mathrm{Ref}} \B^*_*(24(p+1),24(p+1);(M,12\cdot a,8\cdot b,4\cdot c),c,c).\]
Moreover, the latter group is reflection isomorphic to $J(M,12\cdot a,8\cdot b,6\cdot c)$ we have the fourth statement in Theorem \ref{theo:generalized_j_groups_torsion_quotients}.
\end{lemma}
\begin{proof}
Up to permutation, we have 
\[J\jg{k_1&\cdots&k_p&2a&3b&4c\\1&\cdots&1&2&3&4}\cong_{\mathrm{Ref}}J\jg{k_1&\cdots&k_p&3b&2a&4c\\1&\cdots&1&3&2&4}\]
We can apply Lemma \ref{lem:pirouette_avec_poufpoufpouf} with $d=2,n=1,m=2$ and $k=3$. Again by Proposition \ref{prop:permutation_of_parameters_generalized_j_groups}, we obtain reflection isomorphisms
\[J\jg{k_1&\cdots&k_p&3b&2a&4c\\1&\cdots&1&3&2&4}\cong_{\mathrm{Ref}}J\jg{2\cdot(k_1&\cdots&k_p&3b)&a&2c\\2\cdot(1&\cdots&1&3)&1&2}\cong_{\mathrm{Ref}}J\jg{2\cdot(k_1&\cdots&k_p)&a&2c&3b&3b\\2\cdot(1&\cdots&1)&1&2&3&3}.\]
Using Lemma \ref{lem:generalized_j_groups_torsion_quotients_case_4}, this last group is in turn reflection isomorphic to 
\[\B_*^*(12(2p+2),12(2p+2);(N,8\cdot b,4\cdot c),c,c)\]
where $N=12\cdot (2\cdot (k_1,\ldots,k_p),a)$. Again by Proposition \ref{prop:permutation_of_parameters_generalized_j_groups}, this group is reflection isomorphic to 
\[\B^*_*(24(p+1),24(p+1);(M,12\cdot a,8\cdot b,4\cdot c),c,c)\]
since $(M,12\cdot a,8\cdot b,4\cdot c)$ is equal to $(N,8\cdot b,4\cdot c)$ up to permutation. The second statement is a direct application of Lemma \ref{lem:b**(p,p)=J}.
\end{proof}

\subsection{The family $F(2,3,5)$}
At this stage, it remains to study the family $F(2,3,5)$, that is groups of the form 
\[J\jg{k_1&\cdots&k_p&2a&3b&5c\\1&\cdots&1&2&3&5}.\]
We cannot take advantage of Lemma \ref{lem:pirouette_avec_poufpoufpouf} since $2,3,5$ are pairwise coprime. We use a more direct approach, which is both more intricate and computational.

We first deal with the case where every $k_i$ is infinite. Let $p\geqslant 0$ be an integer. 

\begin{notation}Let $d\geqslant 1$ be an integer. Exceptionally in this section, we sometimes write $G(d)$ for the group $J(d\cdot \infty)$.
\end{notation}

We consider the group $G(p+3)$, and we label its generators by $x_1,\ldots,x_p,s,t,u$. We denote by $H_p$ the normal closure of $\{x_1,\ldots,x_p,s^2,t^3,u^5\}$ in $G(p+3)$. We plan to prove the following proposition

\begin{proposition}\label{prop:235_enfin}
The group $H_p$ is reflection isomorphic to $J((60p+62)\cdot \infty)$.
\end{proposition}

By definition, a presentation for the quotient $G(p+3)/H_p$ is given by
\[\langle s,t,u~|~stu=tus=ust,s^2=t^3=u^5=1\rangle.\]
By \cite[Table 2]{BMR}, this presentation is a presentation of the complex reflection group $G_{19}$, with $s,t,u$ as generating reflections. The center of $G_{19}$ is a cyclic group of order $60$ generated by $stu$. We identify the quotient $G(p+3)/H_p$ with $G_{19}$ from now on.

The case of $p=0$ can easily be dealt with by interpreting $H_0$ as a fundamental group.
\begin{lemma}\label{lem:H_0}
The group $H_0$ is reflection isomorphic to $J(62\cdot \infty)$.
\end{lemma}
\begin{proof}
The quotient $G(3)/H_0$ is isomorphic to $G_{19}$, which admits $62$ reflecting hyperplanes (=lines) when acting on $\C^2$. Let us denote by $X$ the complement in $\C^2$ of these reflecting hyperplanes. 

By \cite[Table 2]{BMR}, the defining presentation of $G(3)$ is also a presentation of the braid group $B(G_{19})$, that is the fundamental group of $X/W$. Moreover, the generators of this presentations are generators-of-the-monodromy, i.e. braided reflections in $B(G_{19})$.

The fundamental group of $X$ is denoted by $P(G_{19})$. It is the kernel of the projection map $B(G_{19})\to G_{19}$ and we have an isomorphism of short exact sequences
\[\begin{tikzcd}[ampersand replacement=\&]
	1 \& {H_0} \& {G(3)} \& {G_{19}} \& 1 \\
	1 \& {P(G_{19})} \& {B(G_{19})} \& {G_{19}} \& 1
	\arrow[from=1-1, to=1-2]
	\arrow[from=1-2, to=1-3]
	\arrow["\simeq", from=1-2, to=2-2]
	\arrow[from=1-3, to=1-4]
	\arrow["\simeq", from=1-3, to=2-3]
	\arrow[from=1-4, to=1-5]
	\arrow[equals, from=1-4, to=2-4]
	\arrow[from=2-1, to=2-2]
	\arrow[from=2-2, to=2-3]
	\arrow[from=2-3, to=2-4]
	\arrow[from=2-4, to=2-5]
\end{tikzcd}\]

We claim that the isomorphism $H_0\simeq P(G_{19})$ is the desired isomorphism.

Now, $X$ is the complement of $62$ lines in $\C^2$. By \cite{randell}, the fundamental group of $X$ is isomorphic to $J(62\cdot \infty)$, where the generators are generators-of-the-monodromy in the sense of \cite{BMR}. Moreover, by \cite[Proposition A3]{BMR}, the generators of the monodromy in $P(G_{19})$ are (up to conjugacy) equal to $s^2,t^3$ or $u^5$ in $B(G_{19})\simeq G(3)$. This gives the result.
\end{proof}

In order to complete the proof of Proposition \ref{prop:235_enfin} for $p\geqslant 1$, we first need a few intermediate lemmas.

\subsubsection{Group theoretic lemmas}
\begin{lemma}\label{lem:redresse}
Let $p\geqslant 1$ be an integer, and let $G$ be a group isomorphic to $G(p)$. Assume that $X=\{x_1,\ldots,x_p\}$ is a subset of $G$ such that
\begin{itemize}
\item The set $X$ generates $G$,
\item The product $x_1\cdots x_p$ generates $Z(G)$.
\end{itemize}
Then a group presentation of $G$ is given by
\[G=\langle x_1,\ldots,x_p~|~x_1\cdots x_p=x_2\cdots x_p x_1=\ldots=x_px_1\cdots x_{p-1}\rangle.\] 
\end{lemma}
\begin{proof}
Let $z$ denote the product $x_1\cdots x_p$. Since $z$ is central in $G$, we have
\[x_1\cdots x_p=x_2\cdots x_px_1=\ldots=x_px_1\cdots x_{p-1}\]
Let $a_1,\ldots,a_p$ be generators of the group $G(p)$. The correspondence $a_i\mapsto x_i$ for $i\in \intv{1,p}$ induces a well-defined morphism $\varphi:G(p)\to G$, which is surjective since $X$ generates $G$. In order to conclude, it remains to show that $\varphi$ is injective.

Since $G$ is (abstractly) isomorphic to $G(p)$, the quotient $\bar{G}:=G/Z(G)$ is a free group on $p-1$ generators. Since the product $a_1\cdots a_p$ generates $Z(G(p))$, and since the product $x_1\cdots x_p$ generates $Z(G)$, $\varphi$ induces a surjective morphism $\bar{\varphi}$ from $G(p)/Z(G(p))$ to $\bar{G}$. Now, since free groups of finite rank are Hopfian, and since $G(p)/Z(G(p))$ and $\bar{G}$ are both free groups on $p-1$ generators, we obtain that $\bar{\varphi}$ is an isomorphism. 

We then have $\mathrm{Ker}(\varphi)\subset Z(G(p))$, which is a cyclic group generated by $a_1\cdots a_p$. Since $\varphi(a_1\cdots a_p)=z$ is nontrivial, and since $G$ is torsion free, we finally deduce that $\varphi$ is injective. 
\end{proof}

The main purpose of Lemma \ref{lem:redresse} is to replace an abstract group isomorphism $G\simeq G(p)$ with an explicit isomorphism.

\begin{lemma}\label{lem:normal_tested_on_gens}
Let $G$ be a group and let $H\leqslant G$ be a subgroup. Let $X\subset G$ be a subset which positively generates $G$, and let $Y\subset H$ be an arbitrary generating set. The group $H$ is normal in $G$ if and only if
\[\forall x\in X,y\in Y, xyx^{-1}\in H\]
\end{lemma}
\begin{proof}
The only if part is immediate since $xyx^{-1}$ is a conjugate of an element of $H$. Conversely, let $h\in H$, and let $x\in X$. We can write $h$ as a product of elements of $Y$. The element $xhx^{-1}$ is then a product of conjugates of elements of $Y$ by $X$. Since all these elements belong to $H$ by assumption, we have $xhx^{-1}\in H$. Now, any element $g\in G$ can be written as a product $x_1\cdots x_n$ with $x_i\in X$ for $i\in \intv{1,n}$. An immediate induction then yields that $ghg^{-1}$ belongs to $H$.
\end{proof}

\begin{lemma}\label{lem:abstract}
The group $H_p$ is abstractly isomorphic to $J((60p+62)\cdot \infty)$.
\end{lemma}
\begin{proof}
Let $F$ denote the quotient $G(p+3)/Z(G(p+3))$, and let $\bbar{x_1},\ldots,\bbar{x_p},\bbar{s},\bbar{t},\bbar{u}$ denote the respective images of $x_1,\ldots,x_p,s,t,u$ in $F$. The group $F$ is a free group generated by $\bbar{x_1},\ldots,\bbar{x_p},\bbar{s},\bbar{t}$.

Let $\bbar{H_p}$ denote the image of $H_p$ in $F$. By Corollary \ref{cor:center_generalized_j_groups} (center of generalized $J$-groups), we have $Z(H_p)=Z(G(p+3))\cap H_p$, and thus $\bbar{H_p}\simeq H_p/Z(H_p)$. We have a morphism of short exact sequences
\[\begin{tikzcd}[ampersand replacement=\&]
	1 \& {H_p} \& {G(p+3)} \& {G_{19}} \& 1 \\
	1 \& {\bbar{H_p}} \& F \& {G_{19}/Z(G_{19})} \& 1
	\arrow[from=1-1, to=1-2]
	\arrow[from=1-2, to=1-3]
	\arrow[from=1-2, to=2-2]
	\arrow[from=1-3, to=1-4]
	\arrow[from=1-3, to=2-3]
	\arrow[from=1-4, to=1-5]
	\arrow[from=1-4, to=2-4]
	\arrow[from=2-1, to=2-2]
	\arrow[from=2-2, to=2-3]
	\arrow[from=2-3, to=2-4]
	\arrow[from=2-4, to=2-5]
\end{tikzcd}\]

The group $G_{19}/Z(G_{19})$ is isomorphic to the alternating group $A_5$. In particular, $\bbar{H_p}$ is a subgroup of index $60$ in a free group on $p+2$ generators. By \cite[Proposition 3.9]{reisch}, $\bbar{H_p}$ is a free group on $60(p+1)+1=60p+61$ generators.

We have a short exact sequence
\[1\to Z(H_p)\to H_p\to \bbar{H_p}\to 1\]
Which splits since $\bbar{H_p}$ is free. We obtain that $H_p$ decomposes as a direct product of $\Z$ with a free group on $60p+61$ generators, which finishes the proof since such a direct product is isomorphic to $J((60p+62)\cdot \infty)$.\end{proof}

By combining Lemma \ref{lem:abstract} and Lemma \ref{lem:redresse}, proving Proposition \ref{prop:235_enfin} amounts to finding a system of $60p+62$ generating reflections of $H_p$ whose product (in a well-chosen ordering) generates the center of $H_p$. Note that, since the center of $G_{19}$ is a cyclic group of order $60$, the center of $H_p$ is generated by $(x_1\cdots x_p stu)^{60}$.

In order to complete the proof, we distinguish the cases $p=1$ and $p>1$.

\subsubsection{The case $p=1$}\label{sec:p=1_pour_235}
This part of the proof is achieved through direct computations using {\tt GAP4}. Let us denote by $x,s,t,u$ the generators of $G(4)$. 
\begin{itemize}
\item $gxg^{-1}$ for $g$ in 
\[\begin{array}{lllll} utu,&ut,&utsu,&utsutu,&utsutuu,\\
uttu,&uu,&uutu,&uut,&uutst,\\
uutstsu,&uutsu,&uutsutuxu,&uutsutu,&uuu,\\
1,&st,&stsu,&stsutu,&stsutuu,\\
stsutxstxst,&stsut,&su,&sutuxu,&sutu,\\
sut,&sutst,&sutstsu,&sutstsuu,&sutsu,\\
sutsutuxu,&sutsutu,&sutsut,&sutsutsttu,&suu,\\
suutstu,&suutstuu,&suutstxstxst,&suutst,&suuttu,\\
tu,&t,&tst,&tsu,&tsutuxu,\\
tsutu,&tsut,&tsutst,&tsutstsu,&tsutstsuu,\\
tsutsttuxxstxst,&tsutsttu,&tsutsu,&tsutsutu,&tsutsuttu,\\
tsuttuxuu,&tsuttu,&ttuxuu,&ttu,&tttu\end{array}\]
Let us label these words by $m_1,\ldots,m_{60}$ with $m_1=utu$, $m_2=ut$ and so on. We write $x_i$ for $m_i x m_i^{-1}$.
\item $gs^2g^{-1}$ for $g$ in
\[\begin{array}{llllll}u,&utsts,&utsu,&utsutsu,&utsut,&uts,\\
uu,&uutst,&uutsus,&uutsuts,&uuts,&uuus,\\
uuut,&stsut,&st,&sus,&sutsts,&sutsus,\\
sutsut,&suts,&suus,&suutst,&suuts,&1,\\
tst,&tsu,&tsutsts,&tsutsu,&tsuts,&t\end{array}\]
Let us label these by $s_1\ldots s_{30}$.
\item $gt^3g^{-1}$ for $g$ in
\[\begin{array}{lllll}
utst,&ut,&uutsts,&uutsus,&uutsu,\\
uut,&sts,&stsu,&sus,&sutsus,\\
sutsuts,&sutsut,&su,&suutst,&suu,\\
tsus,&tsuts,&tsutsut,&tsut,&1\end{array}\]
Let us label these by $t_1,\ldots,t_{20}$
\item $gu^5g^{-1}$ for $g$ in 
\[\begin{array}{llllll}
us,&utsut,&utsu,&uutsts,&1,&stsut,\\
stsu,&sutsts,&suutst,&suuts,&t,&tsutsts
\end{array}\]
Let us label these by $u_1,\ldots,u_{12}$
\end{itemize}

The group $G(4)$ has a solution to the word problem (either considering that it is a direct product of a free group with an infinite cyclic group, or using Garside theory). We can then check directly that the product
\begin{align*}
&u_{1}~s_{1}~x_{1}~x_{2}~s_{2}~t_{1}~x_{3}~s_{3}~x_{4}~x_{5}~u_{2}~s_{4}~s_{5}~u_{3}~s_{6}~x_{6}~t_{2}~x_{7}~s_{7}~x_{8}~\\
&x_{9}~x_{10}~t_{3}~x_{11}~u_{4}~s_{8}~x_{12}~t_{4}~s_{9}~x_{13}~x_{14}~s_{10}~t_{5}~s_{11}~t_{6}~x_{15}~s_{12}~s_{13}~u_{5}~x_{16}~\\
&x_{17}~t_{7}~x_{18}~x_{19}~x_{20}~u_{6}~x_{21}~x_{22}~s_{14}~t_{8}~u_{7}~s_{15}~x_{23}~t_{9}~s_{16}~x_{24}~x_{25}~x_{26}~x_{27}~x_{28}~\\
&x_{29}~u_{8}~s_{17}~x_{30}~t_{10}~s_{18}~x_{31}~x_{32}~x_{33}~x_{34}~t_{11}~s_{19}~t_{12}~s_{20}~t_{13}~x_{35}~s_{21}~x_{36}~x_{37}~u_{9}~\\
&x_{38}~x_{39}~s_{22}~t_{14}~u_{10}~s_{23}~x_{40}~t_{15}~s_{24}~x_{41}~u_{11}~x_{42}~x_{43}~s_{25}~x_{44}~t_{16}~s_{26}~x_{45}~x_{46}~x_{47}~\\
&x_{48}~x_{49}~x_{50}~u_{12}~s_{27}~x_{51}~x_{52}~t_{17}~x_{53}~s_{28}~x_{54}~x_{55}~t_{18}~s_{29}~x_{56}~x_{57}~t_{19}~s_{30}~x_{58}~x_{59}~\\
&x_{60}~t_{20} \stepcounter{equation}\tag{\theequation}\label{eq:produ}
\end{align*}
is equal to $(xstu)^{60}$ in $G(4)$. Since we have a system of $60+62$ reflections in $H_1$ whose product generates $Z(H_1)$, it only remains to show that this system generates $H_1$:

\begin{lemma}\label{lem:gens_for_235p=1}
The group $H_1$ is generated by
\[\{x_1,\ldots,x_{60}\}\cup\{s_1\ldots s_{30}\}\cup\{t_1,\ldots,t_{20}\}\cup\{u_1,\ldots,u_{12}\}\]
\end{lemma}
\begin{proof}
Let us denote by $X$ our candidate generating set, and let $H$ be the subgroup of $G(4)$ generated by $X$. Since $X$ consists of conjugates of $x,s^2,t^3,u^5$ in $G(4)$, we have $H\subset H_1$. Moreover, we have shown that $(xstu)^{60}$ lies in $H$. Since $Z(H_1)$ is generated by $(xstu)^{60}$, we have $H=H_1$ if and only if the image of $H$ in $H_1/Z(H_1)$ is equal to $H_1/Z(H_1)$.

Now, since $Z(H_1)=Z(G(4))\cap H_1$, the group $H_1/Z(H_1)$ is identified with the image of $H_1$ in the quotient $G(4)/Z(G(4))$. Let $F$ be a free group on three letters $\chi,\sigma,\tau$. The correspondence
\[\begin{cases} x\mapsto \chi,\\s\mapsto \sigma,\\ t\mapsto \tau,\\ u\mapsto (\chi\sigma\tau)^{-1},\end{cases}\]
induces a group morphism $G(4)\to F$, and an isomorphism $G(4)/Z(G(4))\simeq F$. Under this isomorphism, the image of $H_1$ is identified with the normal closure of $\chi,\sigma^{2},\tau^3,(\chi\sigma\tau)^{-5}$, while the image of $H$ is identified with the subgroup generated by the image of $X$. 

Since the image of $H_1$ in $F$ has finite index, $H_1$ is a finitely generated free group. And we can obtain a generating set by the Reidemeister-Schreier method. Moreover, the equality problem is decidable for finitely generated subgroups of free groups \cite[Proposition 2.21]{reisch}. We use the implementation in {\tt GAP4} of both the Reidemeister-Schreier method and of the equality problem for finitely generated subgroups of free group to conclude that the images of $H$ and $H_1$ in $F$ are equal, which is sufficient to conclude that $H$ and $H_1$ are equal.
\end{proof}

\subsubsection{The case $p>1$}
Our strategy here will be to embed a copy of $H_1$ and $G(4)$ in $H_p$ and $G(p+3)$ in order to use the results of the above section. Consider the elements $x:=x_1\cdots x_p$, $s,t,u$ and let $K$ be the subgroup of $G(p+3)$ generated by $x,s,t,u$.
 
\begin{lemma}\label{lem:k_is_g4}
We have the group presentation $K=\langle x,s,t,u~|~xstu=stux=tuxs=uxst\rangle$. In particular $K\simeq G(4)$.
\end{lemma}
\begin{proof}
First, note that $z:=xstu$ generates the center of $G(p+3)$. In particular we have $Z(G(p+3))\subset G(4)$. Now, $G(p+3)$ decomposes as a direct product $F\times Z(G(p+3))$, where $F$ is the subgroup generated by $x_1,\ldots,x_p,s,t$ (which is free). Under this decomposition, $K$ becomes the direct product $\langle x,s,t\rangle \times \langle z\rangle$. The group $\langle x,s,t\rangle$ is free on $x,s,t$, and thus a presentation of $K$ is given by
\[K=\langle x,s,t,z~|~xz=zx,sz=zs,tz=zt\rangle.\]
Knowing that $u=(xst)^{-1}z$, we obtain the desired presentation.  
\end{proof}

Using this lemma, we identify the group $K$ with $G(4)$ from now on. The restriction of the projection map $G(p+3)\to G_{19}$ to $K$ is surjective since $G_{19}$ is generated by the images of $s,t,u$. The kernel of this restriction is equal to $H_p\cap K$. We can then also identify $H_p\cap K$ with $H_1$.

Let us consider the sixty words $m_1,\ldots,m_{60}$ introduced in Section \ref{sec:p=1_pour_235}, along with the set 
\[S:=\{s_1,\ldots,s_{30}\}\cup \{t_1,\ldots,t_{20}\}\cup \{u_1,\ldots,u_{12}\}.\]
By Lemma \ref{lem:gens_for_235p=1}, we know that $H_1$ is generated by $S$, along with the $m_i xm_i^{-1}$ for $i\in \intv{1,60}$. By replacing each occurrence of $x_i$ in the product \eqref{eq:produ} with the product $m_ix_1 m_i^{-1}\cdots m_ix_p m_i^{-1}$, we obtain a decomposition of $(x_1\cdots x_p stu)^{60}$ as a product of $60p+62$ reflections of $H_p$. In order to complete the proof of Proposition \ref{prop:235_enfin} in this case, it remains to show that $H_p$ is generated by $S$, along with the elements $x_{i,j}=m_jx_im_j^{-1}$ for $i\in \intv{1,p},j\in \intv{1,60}$.

\begin{proposition}\label{prop:generating_p>1}
The group $H_p$ is generated by
\[\{x_{i,j}~|~i\in \intv{1,p},j\in \intv{1,60}\}\cup S\]
\end{proposition}
\begin{proof}
Let us denote by $X$ our candidate generating set, and let $H$ be the subgroup of $G(p+3)$ generated by $X$. Since $X$ consists of conjugates of the $x_i$'s, along with conjugates of $s^2,t^3,u^5$ in $G(p+3)$, we have $H\subset H_p$. It is then sufficient to show that $H$ is normal in $G(p+3)$.

Let $z:=x_1\cdots x_pstu$. The group $G(p+3)$ is positively generated by the set $A:=\{z,z^{-1},x_1,\ldots,x_p,s,t,u\}$ (this is a standard Garside theoretic fact, see for instance \cite[Proposition I.2.4]{GrosBouquinBleu}). By Lemma \ref{lem:normal_tested_on_gens}, it is sufficient to show that any conjugate of an element of $X$ by an element of $A$ lies in $H$. This is of course true when conjugating by $z$ or $z^{-1}$. Moreover, since $m_{16}$ is the trivial word, we have $x_1,\ldots,x_p\in  H$. In particular, conjugating by $x_1,\ldots,x_p$ also leaves $H$ globally invariant. It remains to show that any conjugate of an element of $X$ by an element of $\{s,t,u\}$ lies in $H$.

Now, let $r\in \{s,t,u\}$ and $g\in S$. The element $rgr^{-1}$ lies in $H_p\cap G(4)=H_1$. Since the generating set of $H_1$ given in Lemma \ref{lem:gens_for_235p=1} is included in $H$, we have $H_1\subset H$ and thus $rgr^{-1}\in H$. It remains to consider the conjugates of the generators $x_{i,j}$ by an element of $\{s,t,u\}$.

Note that the word $m_1,\ldots,m_{60}$ form a complete system of representative of the quotient $G_{19}/Z(G_{19})$. Let $j\in \intv{1,60}$, and let $r\in \{s,t,u\}$. The image in $G_{19}/Z(G_{19})$ of $rm_j$ is represented by some word $m_k$. By Proposition \ref{prop:conj_reflections_j-groupe} (conjugacy classes of reflections), the elements $rm_j x m_{j}^{-1}r^{-1}$ and $m_k x m_k^{-1}$ are then conjugate in $H_1$, say by an element $h$. In particular, the elements $rm_j$ and $hm_k$ are equal modulo the centralizer of $x$ in $G(4)$. By Proposition \ref{prop:centralizer_reflections} (centralizer of a reflection), we then have $rm_j=hm_kx^n(xstu)^{m}$ in $G(4)$ for some $n,m\in \Z$. This implies that 
\begin{align*}
rx_{i,j}r^{-1}&=(rm_j) x_i (rm_j)^{-1}\\
&=(hm_k x^n) x_i (hm_k x^n)^{-1}\\
&=(hm_k x^nm_k^{-1}) x_{i,k} (hm_k x^nm_k^{-1})^{-1}\\
&=(h (x_{1,k}\cdots x_{p,k})^n) x_{i,k} (h (x_{1,k}\cdots x_{p,k})^n)^{-1}
\end{align*}
Lastly, since $h\in H_1\subset H$, the element $(h (x_{1,k}\cdots x_{p,k})^n)$ lies in $H$, and thus $rx_{i,j}r^{-1}$ also lies in $H$, which finishes the proof.  
\end{proof}

\subsubsection{The general case}
We can now prove that a generalized $J$-group in the family $F(2,3,5)$ is a torsion quotient of $J$-braid group:

\begin{proposition}[\textbf{Family $F(2,3,5)$}]\label{prop:generalized_j_groups_torsion_quotients_case_5}
Let $M:=60\cdot (k_1,\ldots,k_p)$. We have a reflection isomorphism
\[J\jg{k_1&\cdots&k_p&2a&3b&5c\\1&\cdots&1&2&3&5}\cong_{\mathrm{Ref}} J(M,30\cdot a,20\cdot b,12\cdot c).\]
Moreover, we have the last statement in Theorem \ref{theo:generalized_j_groups_torsion_quotients}.
\end{proposition}
\begin{proof}
Let $\widetilde{J}:=J((p+3)\cdot \infty)$, and let $J$ be the normal closure in $\widetilde{J}$ of $\{x_1,\ldots,x_p,$ $s^2,t^3,u^5\}$. 

We consider the generating set of $J$ given in Proposition \ref{prop:generating_p>1} (which is also a generating set when $p\in\{0,1\}$). Let $Q:=G_{19}/Z(G_{19})$. Notice that $Q$ is the quotient of $J(2,3,5)$ by the image of the center of $\widetilde{J}$. 

Consider the words $m_1,\ldots,m_{60}$ introduced in Section \ref{sec:p=1_pour_235}. They form a complete set of representatives of $Q$. Let $i\in \intv{1,p}$. Since the image of $x_i$ is trivial in $Q$, Proposition \ref{prop:conj_reflections_j-groupe} implies that an element of $J$ is conjugate to $x_i^{k_i}$ in $\widetilde{J}$ if and only if it is conjugate in $J$ to some $m_j(x_i^{k_i})m_j^{-1}$.  The normal closure of $x_i^{k_i}$ in $\widetilde{J}$ is then equal to the normal closure of $\{x_{i,j}^{k_i},j\in \intv{1,60}\}$ in $J$.

Similarly, the word introduced in Section \ref{sec:p=1_pour_235} defining $s_1,\ldots,s_{30}$ (resp. $t_1,\ldots,t_{20}$, $u_1,\ldots,u_{12})$ form a complete set of representatives of the cosets in $Q$ relative to the image of $s$ (resp. of $t$, of $u$). Again by Proposition \ref{prop:conj_reflections_j-groupe}, the normal closure of $s^{2a}$ (resp. of $t^{3b}$, $u^{5c}$) in $\widetilde{J}$ is the normal closure of $\{s_1^{a},\ldots,s_{30}^{a}\}$ (resp. of $\{t_1^{b},\ldots,t_{20}^{b}\}$, $\{u_1^{c},\ldots,u_{20}^{c}\}$) in $J$. 

Since $J\jg{k_1&\cdots&k_p&2a&3b&5c\\1&\cdots&1&2&3&5}$ is obtained from $J$ by quotienting by the normal closure in $\widetilde{J}$ of $\{x_1^{k_1},\ldots,x_i^{k_i},s^{2a},t^{3b},u^{5c}\}$, the isomorphism of Proposition \ref{prop:235_enfin} induces the desired reflection isomorphism. The last statement in Theorem \ref{theo:generalized_j_groups_torsion_quotients} is obtained by applying Lemma \ref{lem:b**(p,p)=J}.
\end{proof}

\section{Classification up to reflection isomorphism}\label{sec:classif}
In this section, we establish the classification of torsion quotients of $J$-braid groups up to reflection isomorphism. As a corollary, we will deduce the classification of finite-type generalized $J$-groups up to reflection isomorphism. 

\subsection{Preliminary results}\label{sec:preliminary_results}

A useful tool to establish the classification is the concept of reflecting hyperplane. Since torsion quotients of $J$-braid groups are not ``actual'' reflection groups (in the sense that we did not define a reflection representation of these groups), we cannot define reflecting hyperplanes in the usual geometric way. However, the usual definition of reflecting hyperplanes for a linear reflection groups coincides with the following more combinatorial analogue, which is introduced in \cite{GobetToric} for toric reflection groups.

In this section, we fix $W$ a torsion quotient of a $J$-braid group.

\begin{definition}[\textbf{Reflecting hyperplane}]\label{def:reflecting_hyperplane} Let $\sim$ be the relation on $R(W)$ generated by the relations $r^a\sim r^b$ for all $1\leqslant a,b<o(r)$, $r\in R(W)$ and write $[r]$ for the equivalence class of $r\in R(W)$. The set $\mathcal{H}(W)$ of \emph{reflecting hyperplanes} of $W$ is defined to by $\{[r]\}_{r\in R(W)}$.\end{definition}

The action of $W$ on $R(W)$ by conjugacy induces an action of $W$ on $\mathcal{H}(W)$. The $W$-orbit of $[r]$ will be denoted by $W.[r]$. 

There is a natural map $O:\Hh(W)\to \N_{\geqslant 2}\cup \{\infty\}$, which sends $[r]$ to $\max \{o(s)~|~s\in [r]\}$. Since the action of $W$ on $R(W)$ preserves the order of the reflections, the map $O$ is $W$-invariant.

\begin{notation}
A (finite) multiset of cardinality $n$ is defined as the orbit under the symmetric group $S_n$ of a $n$-tuple. If $(x_1,\ldots,x_n)$ is a $n$-tuple, then the associated multiset will be denoted by $\muls{x_1,\ldots,x_n}$.
\end{notation}

\begin{definition}[\textbf{Torsions}]
Let $[r_1],\ldots,[r_k]$ be a set of representatives of the $W$-orbits of $\Hh(W)$. The multiset of \emph{torsions} $T(W)$ is defined by
\[T(W):=\muls{ O([r_i])~|~i\in \intv{1,k}}.\]
Notice that $T(W)$ does not depend on the choice of the $r_i$'s.
\end{definition}

These concepts are useful as they are invariants under reflection isomorphism:

\begin{proposition}\label{prop:ref_iso_implies_torsion_equal}
Let $W_1,W_2$ be two torsion quotients of $J$-braid groups, and let $\varphi:W_1\to W_2$ be a reflection isomorphism. 
\begin{enumerate}[(a)]
\item The morphism $\varphi$ induces a bijection $\Hh(W_1)\to \Hh(W_2)$ sending $[r]$ to $[\varphi(r)]$.
\item The above bijection induces a bijection $\Hh(W_1)/W_1\to \Hh(W_2)/W_2$ sending $W_1.[r]$ to $W_2.[\varphi(r)]$.
\item The two multisets $T(W_1)$ and $T(W_2)$ are equal.
\end{enumerate}
\end{proposition}
\begin{proof}
$(a)$ By definition, $\varphi$ restricts to a bijection $R(W_1)\to R(W_2)$. For $r\in R(W_1)$ and $a,b\in \intv{1,o(r)-1}$, we have $[\varphi(r)^a]=[\varphi(r)^b]$ in $\Hh(W_2)$ by definition of $\Hh(W_2)$. Thus the map $[r]\mapsto [\varphi(r)]$ is a well-defined map $\Hh(W_1)\to \Hh(W_2)$. Considering the inverse morphism $\varphi^{-1}$ (which is also a reflection isomorphism), we obtain that the map $[r]\mapsto [\varphi(r)]$ is actually a bijection. 

$(b)$ Let $\bbar{\varphi}$ denote the bijection $\Hh(W_1)\to \Hh(W_2)$ given by point $(a)$. Let $r\in R(W_1)$ and $w\in W$. We have
\[\bbar{\varphi}(w.[r])=\bbar{\varphi}([wrw^{-1}])=[\varphi(wrw^{-1})]=\varphi(w).[\varphi(r)]=\varphi(w).\bbar{\varphi}([r]),\]
and thus $\bbar{\varphi}$ induces a well defined map $\Hh(W_1)/W_1\to \Hh(W_2)/W_2$ sending $W_1.[r]$ to $W_2.[\varphi(r)]$. Again, considering the inverse morphism $\varphi^{-1}$ gives that this map is actually a bijection.

$(c)$ Let us denote by $O_1,O_2$ the respectives maps $O$ for the groups $W_1$ and $W_2$. For $r\in R(W_1)$, we have
\[O_2(\bbar{\varphi}([r]))=O_2([\varphi(r)])=\max \{o(s)~|~s\in [\varphi(r)]\}.\]
Now, $\varphi$ induces an order preserving bijection $[r]\to [\varphi(r)]$. Thus the above maximum is also the maximum of the set $\{o(s)~|~s\in [r]\}$, that is $O_1([r])$. If $[r_1],\ldots,[r_k]$ is a set of representatives of the $W_1$ orbits of $\Hh(W_1)$, then $[\varphi(r_1)],\ldots,[\varphi(r_k)]$ is a set of representatives of the $W_2$ orbits of $\Hh(W_2)$ by point $(b)$. We then have
\[T(W_1)=\muls{O_1([r_i])~|~i\in \intv{1,k}}=\muls{O_2([\varphi(r_i)])~|~i\in \intv{1,k}}=T(W_2).\]
\end{proof}

The (multiset) cardinality of $T(W)$ is the cardinality of $\Hh(W)/W$, which is also invariant under reflection isomorphism. 

\begin{remark}
Notice that the proof of Proposition \ref{prop:ref_iso_implies_torsion_equal} is purely formal and doesn't use the particular theory of torsion quotients of $J$-braid groups. In particular, the definition and results of this section also applies to other contexts where a concept of reflection can be defined (or really, any distinguished generating set invariant under conjugacy). In particular, we can also define the reflecting hyperplanes and the torsion set of a generalized $J$-group.
\end{remark}

\subsection{Reduced $J$-groups}
In the previous section, we have shown that every finite-type generalized $J$-group is reflection isomorphic to a torsion of a $J$-braid group. Corollary \ref{cor:embedding_torsion_quotients} (embedding of torsion quotients) gives a converse to this result, but it is actually more specific. Indeed, a torsion quotient of $J$-braid group is reflection isomorphic to a particular type of finite-type generalized $J$-group, which we call reduced $J$-group. The main interest of this particular family is that we can complete its classification up to reflection isomorphism.

In this section, we fix a positive integer $p$, along with a $p$-tuple $K=(k_1,\ldots,k_p)$ in $(N_{\geqslant 2}\cup \{\infty\})^p$. We also fix another tuple $K'=(k_1',\ldots,k_p')$ in $(\N_{\geqslant 1})^p$ such that $k_i'$ divides $k_i$ for each $i$.

\begin{definition}[\textbf{Reduced $J$-group}]
A generalized $J$-group $J\binom{K}{K'}$ is said to be \emph{reduced} if $K'$ contains at most two entries distinct from $1$, and if all the elements of $K'$ are pairwise coprime.
\end{definition}

\begin{lemma}\label{lem:generalized_j_groups_are_iso_to_reduced}
A generalized $J$-group has finite-type if and only if it is isomorphic to a reduced $J$-group.
\end{lemma}
\begin{proof}
Let $J\binom{K}{K'}$ be a reduced $J$-group. By definition, $J(K')$ is a finite cyclic group, and thus $J\binom{K}{K'}$ has finite-type. Conversely, let $J\binom{K}{K'}$ be a finite type generalized $J$-group. Then $J\binom{K}{K'}$ belongs (up to reflection isomorphism) to one of the families of Lemma \ref{lem:finite_type_generalized_list}. By Lemmas \ref{lem:generalized_j_groups_torsion_quotients_case_2}, \ref{lem:generalized_j_groups_torsion_quotients_case_3}, \ref{lem:generalized_j_groups_torsion_quotients_case_4} and Proposition \ref{prop:generalized_j_groups_torsion_quotients_case_5}, four of the five families of Lemma \ref{lem:finite_type_generalized_list} are actually reflection isomorphic to generalized parent $J$-groups, which are in particular finite-type generalized $J$-group. The only case remaining is the case
\[J\jg{K\\K'}=J\jg{k_1&\cdots&k_p&nb&mc\\1&\cdots&1&n&m}.\]
Le us write $d:=n\wedge m$, $m':=\frac{m}{d}$, $n':=\frac{n}{d}$. By Corollary \ref{cor:embedding_torsion_quotients} (embedding or torsion quotients), we have
\begin{align*}
J\jg{k_1&\cdots&k_p&nb&mc\\1&\cdots&1&n&m}&\cong_{\mathrm{Ref}}\B_*^*(pdm',pdn',d\cdot(k_1,\ldots,k_p),b,c)\\
&\cong_{\mathrm{Ref}} J\jg{d\cdot(k_1&\cdots&k_p)&bn'&cm'\\d\cdot (1&\cdots&1)&n'&m'}
\end{align*}
and the latter is a reduced $J$-group.
\end{proof}

\begin{remark}\label{rem:red_j_is_j-refl}
If $p=3$, then a reduced $J$-group $J\binom{K}{K'}$ can be written as $J\jg{K&bn&cm\\1&n&m}$ with $n$ and $m$ coprime. In this case, the definition of reduced $J$-groups coincide with that of $J$-\emph{reflection group} appearing in \cite{VCRG}.
\end{remark}

\begin{lemma}\label{lem:generating_reflecting_hyperplane}
Let $r\in J\binom{K}{K'}$ be a reflection. The reflecting hyperplane $[r]$ consists of non-trivial powers of a reflection in $J\binom{K}{K'}$.
\end{lemma}

\begin{proof}
Assume that there exists $r_1,r_2\in J\binom{K}{K'}$ and $n,m\in \N_{\geqslant 1}$ such that $r_1^n=r_2^m$. 
We write $r_1=gs_i^{l_1}g^{-1}$ and $r_2=hs_j^{l_2}h^{-1}$. Up to conjugacy, we can assume that $g=1$. We get $s_i^{l_1n}=hs_j^{l_2m}h^{-1}$. In particular, Lemma \ref{lem:conj_refl_parent} tells us that $i=j$ and that $l_1n=l_2m$. 

Since $h$ centralizes $s_j^{l_2m}$, it centralizes $s_j$ by Proposition \ref{prop:centralizer_reflections} (centralizer of a reflection). We then have that $r_2=s_i^{l_2}$ is a power of $s_i$, as well as $r_1=s_i^{l_i}$. Using Lemma \ref{lem:reflexion_belonging_to_j_group}, $r_1$ and $r_2$ are powers of $s_i^{k_i'}$, which concludes the proof
\end{proof}

Recall that for a generalized $J$-group $J\binom{K}{K'}$, the quotient $Q$ of $J(K')$ by the image of $Z(J(K))$ is useful for studying the conjugacy of reflections. In the case where $J\binom{K}{K'}$ is reduced, the group $J(K')$ is cyclic and $Q$ is trivial. In particular by Proposition \ref{prop:conj_reflections_j-groupe} (conjugacy classes of reflections), we obtain

\begin{lemma}\label{lem:conj_in_reduced_j-group}
Let $J\binom{K}{K'}$ be a reduced $J$-group. 
\begin{enumerate}[(a)]
\item Two reflections $r,r'\in R(J\binom{K}{K'})$ are conjugate in $J\binom{K}{K'}$ if and only if they are conjugate in $J(K)$.
\item We have $|\Hh(J\binom{K}{K'})|=p$ and $T(J\binom{K}{K'})$ is the submultiset of $\muls{K/K'}$ consisting of elements different from 1, and where $K/K'=(k_1/k_1',\ldots,k_p/k_p')$.
\end{enumerate}
\end{lemma}
\begin{proof}
Point $(a)$ is a direct application of Proposition \ref{prop:conj_reflections_j-groupe} since the quotient $Q$ is trivial in this case. After point $(a)$, a complete system of representatives of $J\binom{K}{K'}$-orbits on $\Hh(J\binom{K}{K'})$ is the set of nontrivial elements among, $\{s_1^{k_1'},\ldots,s_p^{k_p'}\}$. Moreover, by Lemma \ref{lem:generating_reflecting_hyperplane}, any $r$ in this set generates $[r]$. Knowing that the torsion of $s_i^{k_i'}$ is $k_i/k_i'$, we have the result. 
\end{proof}

The fact that the quotient group $Q$ evoked above is trivial has another application, which is that we can compute the inner automorphism group of a reduced $J$-group. 

\begin{lemma}\label{lem:inner_autom_red_j-groups}
Assume that $p\geqslant 3$ and that $J\binom{K}{K'}$ is reduced. Then the inner automorphism group $J\binom{K}{K'}$ is isomorphic to the alternating polygonal Coxeter group $W_{K}^+$.
\end{lemma}
\begin{proof}
Since $p\geqslant 3$, we can consider the natural morphism $\pi:J(K)\to W_K^+$, whose kernel is the center of $J(K)$.  By Corollary \ref{cor:center_generalized_j_groups} (center of generalized $J$-groups), the kernel of the restriction of $\pi$ to $J\binom{K}{K'}$ is the center of $J\binom{K}{K'}$. Thus the image $W_K^+(K')$ of $J\binom{K}{K'}$ under the morphism $\pi$ is isomorphic to the inner automorphism group of $J\binom{K}{K'}$. Since $J\binom{K}{K'}$, there are at most two entries $n,m$ in $K'$ which are different from $1$. Thus a presentation for the quotient $W_K^+/W_K^+(K')$ is $\langle x,y~|~x^n=y^m=1, xy=1\rangle$. Since $n$ and $m$ are coprime, this group is trivial so that we have $W_K^+(K')=W_K^+$ and the result is shown.
\end{proof}

On top of representing every finite-type generalized $J$-group up to reflection isomorphism, the reflection isomorphism problem is rather easily solvable for reduced $J$-groups. Our solution to this problem relies on an induction argument using the following lemma 

\begin{lemma}[\textbf{Conjugacy class deletion}]\label{lem:conjugacy_class_deletion}Let $i\in \intv{1,p}$. Let $\hat{K}_{i}$ be the tuple obtained from $K$ by replacing the $i$-th entry with $k_{i}/k'_{i}$. Let alsot $t_1,\ldots,t_{p}$ denote the generators of $J(\hat{K}_{i})$. The correspondence given by $s_j\mapsto t_j$ for $j\in \intv{1,p}$ induces a surjective morphism $\pi:J(K)\to J(\hat{K}_{i})$, whose kernel is normally generated by $s_{i}^{k_{i}'}$. Moreover, if $J\binom{K}{K'}$ is reduced, then $\pi$ restricts to a surjective morphism $J\binom{K}{K'}\to J\binom{\hat{K}_{i}}{K'}$, whose kernel is the normal closure of $s_{i}^{k_{i}'}$ in $J\binom{K}{K'}$.
\end{lemma}

\begin{proof}
Note that if $k_i=k_i'$, then $\hat{K}_i=K$ and there is nothing to show. We then assume that $k_i>k_i'$. The first statement comes from the observation that adding the relation $s_{i}^{k_i'}=1$ to the defining presentation of $J(K)$ yields the defining presentation of $J(\hat{K}_{i})$ (note that if $k_i'=1$, this amounts to deleting $s_i$ entirely). Now, assume that $J\binom{K}{K'}$ is reduced. By Lemma \ref{lem:conj_in_reduced_j-group}, the conjugacy class of $s_{i}^{k_{i}'}$ in $J\binom{K}{K'}$ and in $J(K)$ coincide. Thus the normal closures of $s_{i}^{k_{i}'}$ in $J\binom{K}{K'}$ and in $J(K)$ coincide. 

The kernel of the restriction of $\pi$ to $J\binom{K}{K'}$ is then the normal closure of $s_{j}^{k_j'}$ in $J\binom{K}{K'}$. Moreover, since $J\binom{K}{K'}$ is the normal closure in $J(K)$ of $\{s_j^{k_j'}~|~j\in \intv{1,p}\}$, and since $\pi$ is surjective, the image of $J\binom{K}{K'}$ in $J(\hat{K}_{i})$ is equal to $J\binom{\hat{K}_{i}}{K'}$.
\end{proof}

Notice that in the above lemma, the resulting group $J\binom{\hat{K}_i}{K'}$ is again a reduced $J$-group (up to deleting a possible $\binom{1}{1}$ column).

\begin{proposition}[\textbf{Classification of reduced $J$-groups}]\label{prop:classif_red_j-group}~\\
Let $W_1:=J\binom{K}{K'}$ and $W_2:=J\binom{L}{L'}$ be two reduced $J$-groups. 
\begin{enumerate}[(a)]
\item The group $W_1$ is abelian if and only if $p\leqslant 2$, in which case $W_1\cong_{\mathrm{Ref}} J\jg{k_1/k_1'\\1}$ or $W_1\cong_{\mathrm{Ref}}J\jg{k_1/k_1'&k_2/k_2'\\1&1}$ depending on the length of $K$.
\item If $W_1$ and $W_2$ are nonabelian, then we have $W_1\cong_{\mathrm{Ref}}W_2$ if and only if $p=q$ and if there is $\sigma\in S_p$ such that $\sigma(K)=L$ and $\sigma(K')=L'$.
\end{enumerate}
\end{proposition}
\begin{proof}
$(a)$ The first statement is already known by Corollary \ref{cor:center_generalized_j_groups} since the generalized parent $J$-group is known to be abelian if and only if $p\leqslant 2$. If $p=2$, then $J\binom{K}{K'}=J\jg{k_1&k_2\\k_1'&k_2'}\cong_{\mathrm{Ref}}J\jg{k_1/k_1'&k_2/k_2'\\1&1}$. If $p=1$, then $J\binom{K}{K'}=J\jg{k_1\\k_1'}\cong_{\mathrm{Ref}}J\jg{k_1/k_1'\\1}$.

$(b)$ Let $W_1:=J\binom{K}{K'}$ and $W_2:=J\binom{L}{L'}$ be two reduced $J$-groups, and let $\varphi:W_1\to W_2$ be a reflection isomorphism. By Lemma \ref{lem:inner_autom_red_j-groups}, we have $W_L^+\simeq W_K^+$. By Proposition \ref{ClassificationPolygonalCox}, we have $K=L$ up to permutation, and in particular the lengths of $K$ and $L$ coincide. Let us denote it by $p$. Since we assume $W_1$ and $W_2$ to be nonabelian, we can assume that $p\geqslant 3$. We now proceed by induction on $p$.

The case $p=3$ is already handled in \cite[Theorem 1.6]{VCRG}. Indeed, reduced $J$-groups in this case coincide with the $J$-reflection groups introduced in \emph{loc. cit} (see Remark \ref{rem:red_j_is_j-refl}).

Assume now that $p\geqslant 4$. Up to permuting $K$ and $K'$, we can assume that $k_1'=1$. Considering the bijection $\Hh(W_1)/W_1\to \Hh(W_2)/W_2$ induced by $\varphi$ after Proposition \ref{prop:ref_iso_implies_torsion_equal}, we can consider the unique $i\in \intv{1,p}$ such that $W_2.[\varphi(s_1)]=W_2.[t_i^{\ell_i'}]$. By Lemma \ref{lem:conjugacy_class_deletion}, we have
\[J\jg{\hat{K}_1\\K'}\cong_{\mathrm{Ref}}J\jg{\hat{L}_{i}\\L'}.\]

Since $k_1'=1$, the actual defining matrix $\binom{\kappa}{\kappa'}$ of $J\jg{\hat{K}_1\\K'}$ is obtained from $\jg{\hat{K}_1\\K'}$ by removing one column $\binom{1}{1}$. We then have $|\kappa|=p-1$. As we explained above, the same must hold for the defining matrix $\binom{\lambda}{\lambda'}$ of $J\jg{\hat{L}_i\\L'}$. Which is possible if and only if $\ell'_i=1$ (and $|\lambda|=p-1)$. The induction hypothesis gives that there is some $\sigma\in S_{p-1}$ such that $\sigma(\kappa)=\lambda$ and $\sigma(\kappa')=\sigma(\lambda')$. Since the order of $\varphi(s_1)$ is equal to the order of $s_1$, we have $k_1/k_1'=\ell_i/\ell_i'$. The equality $k_1'=\ell_i'=1$ then gives that the columns $\binom{k_1}{1}$ and $\binom{\ell_i}{1}$ are equal, which finishes the proof.
\end{proof}

\subsection{Classification of torsion quotients of $J$-braid groups}
We are now able to complete the classification of torsion quotients of $J$-braid group up to reflection isomorphism, by taking advantage of Corollary \ref{cor:embedding_torsion_quotients} (embedding of torsion quotients) and of the classification of reduced $J$-reflection groups.

In this section, we fix four positive integers $n,m,p,q$. We let $d:=n\wedge m$ (resp. $d':=p\wedge q$) denote the gcd of $n$ and $m$ (resp. of $p$ and $q$). We also fix two tuples $K:=(k_1,\ldots,k_d)$, $L:=(\ell_1,\ldots,\ell_{d'})$ of elements in $\N_{\geqslant 2} \cup \{\infty\}$. Lastly, we fix $b,c,\beta,\gamma\in \N_{\geqslant 2}\cup \infty$.

We will use the results of Section \ref{sec:preliminary_results}, in particular, the following lemma will be useful:

\begin{lemma}\label{lem:torsion_multiset_of_torsion_quotients}
Let $W$ be a torsion quotient of $J$-braid group.
\begin{itemize}
\item If $W=\B_*^*(n,m;K,b,c)$, then $|\Hh(W)/W|=d+2$ and $T(W)=\muls{K,b,c}$.
\item If $W=\B^*(n,m;K,c)$, then $|\Hh(W)/W|=d+1$ and $T(W)=\muls{K,c}$.
\item If $W=\B_*(n,m;K,b)$, then $|\Hh(W)/W|=d+1$ and $T(W)=\muls{K,b}$.
\item If $W=\B(n,m;K)$, then $|\Hh(W)/W|=d$ and $T(W)=\muls{K}$.
\end{itemize}
\end{lemma}
\begin{proof}
We only prove the first point, the other points are obtained similarly. After Lemma \ref{lem:conjugacy_classes_of_braid_reflections} (conjugacy classes of braid reflections), a complete system of representatives of $W$ orbits on $\Hh(W)$ is $\{[x_1],\ldots,[x_d],[y],[z]\}$. Moreover, any $r\in \{x_1,\ldots,x_d,y,z\}$ generates $[r]$. Knowing that the torsion of $x_1,\ldots,x_d,y,z$ is given by $k_1,\ldots,k_d,b,c$, we have the result. 
\end{proof}

We already have seen several instances of reflection isomorphisms between torsion quotients of $J$-braid groups in Corollary \ref{cor:permutation_parameters_torsion_quotient} (permuting torsion coefficients) and Corollary \ref{cor:swap_parameters_torsion_quotient} (swap of parameters in torsion quotients). Let us list a few more such isomorphisms.

\begin{lemma}\label{lem:iso_exceptionnel}
We have the following reflection isomorphisms
\begin{enumerate}[(a)]
\item $\B^*(n,m;K,c)\cong_{\mathrm{Ref}}\B(n+\frac{n}{m},m+1;(K,c))$ if $m$ divides $n$.
\item $\B_*^*(n,m;K,b,c)\cong_{\mathrm{Ref}}\B_*(n+\frac{n}{m},m+1;(K,c),b)$ if $m$ divides $n$.
\item $\B_*(n,m;K,b)\cong_{\mathrm{Ref}}\B(n+1,m+\frac{m}{n};(K,b))$ if $n$ divides $m$.
\item $\B^*_*(n,m;K,b,c)\cong_{\mathrm{Ref}}\B^*(n+1,m+\frac{m}{n};(K,b),c)$ if $n$ divides $m$.
\item $\B(2,2;(k_1,k_2))\cong_{\mathrm{Ref}}\B^*(1,m;(k_1),k_2)\cong_{\mathrm{Ref}}\B_*(n,1;(k_1),k_2)$.
\item $\B(1,n;(k_1))\cong_{\mathrm{Ref}}\B(1,m;(k_1))$.
\end{enumerate}
\end{lemma}
\begin{proof}
$(a)$ is a consequence of $(b)$. Assume that $m$ divides $n$, in other words that we have $m'=1$, $d=m$ and $n'=\frac nm$. By Corollary \ref{cor:embedding_torsion_quotients} and Corollary \ref{cor:permutation_parameters_torsion_quotient}, we have 
\begin{align*}
    \B^*_*(n,m;K,b,c)&\cong_{\mathrm{Ref}} J\jg{K&bn'&cm'\\d\cdot 1&n'&m'}\\
    &\cong_{\mathrm{Ref}}J\jg{K&bn'&c\\1&n'&1}\\
    &\cong_{\mathrm{Ref}}J\jg{K&c&bn'\\1&1&n'}\\
    &\cong_{\mathrm{Ref}}\B_*(n+\frac{n}{m},m+1;(K,c),b)
\end{align*}
as claimed. $(c)$ and $(d)$ are respective consequences of $(a)$ and $(b)$ using Corollary \ref{cor:swap_parameters_torsion_quotient}.

$(e)$ Again by Corollary \ref{cor:embedding_torsion_quotients} and Corollary \ref{cor:permutation_parameters_torsion_quotient}, we have 
\begin{align*}
    \B^*(1,m;(k_1),k_2)&\cong_{\mathrm{Ref}}J\jg{k_1&k_2m'&1\\1&m'&1}\\
    &\cong_{\mathrm{Ref}}J\jg{k_1&k_2m'\\1&m'}\\
    &\cong_{\mathrm{Ref}}J\jg{k_1&k_2\\1&1}\\
    &\cong\B(2,2;(k_1,k_2)).
\end{align*}
The second isomorphism is obtained using Corollary \ref{cor:swap_parameters_torsion_quotient}.

$(f)$ Both of the defining presentations are $\langle x~|~x^{k_1}=1 \text{ if $k_1<\infty$}\rangle$.
\end{proof}

Since reflection isomorphisms are in particular group isomorphisms, two reflection isomorphic torsion quotients of $J$-braid groups are either both abelian or both nonabelian. The following lemma completely classifies abelian torsion quotients of $J$-braid groups. 

\begin{lemma}[\textbf{Determination of abelian torsion quotients}]\label{lem:classif_of_abelian_torsion_quotients}\hfill
\begin{enumerate}[(a)]
\item The group $\B_*^*(n,m;K,b,c)$ is not abelian.
\item The group $\B^*(n,m;K,c)$ is abelian if and only if $n=1$.
\item The group $\B_*(n,m;K,b)$ is abelian if and only if $m=1$.
\item The group $\B(n,m;K)$ is abelian if and only if $1\in \intv{n,m}$, or if $n=m=2$.
\end{enumerate}
\end{lemma}
\begin{proof}
$(a)$ Let $W:=\B_*^*(n,m;K,b,c)$. By Corollary \ref{cor:embedding_torsion_quotients}, we have
\[W\cong_{\mathrm{Ref}}J\jg{K&bn'&cm'\\d\cdot 1&n'&m'}.\]
Since $b,c\geqslant 2$, the defining matrix of $W$ as a reduced $J$-group has $d+2>2$ nontrivial columns. The group $W$ is then nonabelian by Corollary \ref{cor:center_generalized_j_groups} (center of generalized $J$-group). 

$(b)$ Let $W:=\B^*(n,m;K,c)$. If $n=1$, then we saw in the proof of Theorem \ref{TheoremBraidIsoIntroCircularFR} (isomorphism type of $J$-braid group) that $\B^*(1,m)$ is free abelian generated by $x_1$ and $z$. 

Assume now that $n'\neq 1$. By Corollary \ref{cor:embedding_torsion_quotients}, we have
\[W\cong_{\mathrm{Ref}}J\jg{K&n'&cm'\\d\cdot 1&n'&m'}.\]
If $n'\neq 1$, then the defining matrix of $W$ as a reduced $J$-group has $d+2>2$ nontrivial columns. If $n'=1$, then $d=n\geqslant 2$ by assumption, and the defining matrix of $W$ as a reduced $J$-group then has $d+1>2$ nontrivial columns. In both cases, the group $W$ is nonabelian by Corollary \ref{cor:center_generalized_j_groups}.

$(c)$ The group $\B_*(n,m;K,b)$ is reflection isomorphic to $\B^*(m;n;K,c)$ by Corollary \ref{cor:swap_parameters_torsion_quotient}. The results then comes from point $(b)$.

$(d)$ Let $W:=\B(n,m;K)$. By Corollary \ref{cor:swap_parameters_torsion_quotient}, we can assume that $n\geqslant m$. If $1\in \{n,m\}$ or if $n=m=2$, then $\B(n,m)$ is abelian by Theorem \ref{TheoremBraidIsoIntroCircularFR}.

Assume now that we are not in one of these cases, i.e. $n\geqslant 2$ and $m\geqslant 3$. By Corollary \ref{cor:embedding_torsion_quotients}, we have
\[W\cong_{\mathrm{Ref}}J\jg{K&n'&m'\\d\cdot 1&n'&m'}.\]
If $n'\neq 1$, then $m'\neq 1$ and the defining matrix of $W$ as a reduced $J$-group has $d+2>2$ nontrivial columns. If $n'=1$, then $d=n\geqslant 2$ by assumption, and the defining matrix of $W$ as a reduced $J$-group then has $d+1>2$ nontrivial columns. In both cases, the group $W$ is nonabelian by Corollary \ref{cor:center_generalized_j_groups}.
\end{proof}

\begin{theorem}\label{theo:classif_tors_quot_j-braid}
The reflection isomorphism relation on the set of torsion quotients of $J$-braid groups is generated by the isomorphisms appearing in Corollary \ref{cor:permutation_parameters_torsion_quotient}, Corollary \ref{cor:swap_parameters_torsion_quotient} and Lemma \ref{lem:iso_exceptionnel}. In other words, it is generated by the following isomorphisms: 
\begin{itemize}[itemsep=3pt]
\item Permutation of the tuple $K$.,
\item $\B_*^*(n,m;K,b,c)\cong_{\mathrm{Ref}}\B_*^*(n,m;K,c,b)$,
\item $\B^*(n,m;K,c)\cong_{\mathrm{Ref}}\B_*(m,n;K,c)$,
\item $\B(n,m;K)\cong_{\mathrm{Ref}}\B(m,n;K)$, 
\item $\B^*_{(*)}(n,m;K(,b),c)\cong_{\mathrm{Ref}} \B_{(*)}(n+\frac{n}{m},m+1;(K,c)(,b))$ if $m$ divides $n$,
\item $\B(2,2;(k_1,k_2))\cong_{\mathrm{Ref}} \B^*(1,m;(k_1),k_2)\cong_{\mathrm{Ref}}\B_*(n,1;(k_1),k_2)$,
\item $\B(1,n;(k_1))\cong_{\mathrm{Ref}}\B(1,m;(k_1))$.
\end{itemize}
\end{theorem}

In order to prove Theorem \ref{theo:classif_tors_quot_j-braid}, we have to consider two reflection isomorphic torsion quotient of $J$-braid groups $W_1$ and $W_2$, and to show that they are related by a sequence of reflection isomorphism appearing in either Corollary \ref{cor:permutation_parameters_torsion_quotient}, Corollary \ref{cor:swap_parameters_torsion_quotient} or Lemma \ref{lem:iso_exceptionnel}.

Notice that by Corollary \ref{cor:swap_parameters_torsion_quotient}, we do not need to consider the case where $W_1=\B_*(n,m;K,b)$ (resp. where $W_2=\B_*(p,q;K,\beta)$).

\noindent\underline{Case 0 : $W_1$ and $W_2$ are abelian}
By Lemma \ref{lem:iso_exceptionnel} and Lemma \ref{lem:classif_of_abelian_torsion_quotients}, It is sufficient to consider the cases of groups of the form $\B(2,2;(k_1,k_2))$ or $\B(1,1,(k_1))$. These two groups are not reflection isomorphic as they do not have the same number of conjugacy classes of reflecting hyperplanes.

For the remainder of the proof, we assume that neither $W_1$ nor $W_2$ is abelian.

\noindent\underline{Case 1 : $W_1=\B(n,m;K),W_2=\B(p,q;L)$.} \\First, we have $d=|\Hh(W_1)/W_1|=|\Hh(W_2)/W_2|=d'$ by  Proposition \ref{prop:ref_iso_implies_torsion_equal} and Lemma \ref{lem:torsion_multiset_of_torsion_quotients}. Moreover, these two results also give that $T(W_1)=\muls{K}=\muls{L}=T(W_2)$. Applying Corollary \ref{cor:permutation_parameters_torsion_quotient}, we can assume that $L=K$. Applying Corollary \ref{cor:swap_parameters_torsion_quotient}, we can assume that $n\leqslant m$ and $p\leqslant q$. 

Now, by Corollary \ref{cor:embedding_torsion_quotients}, we have
\[W_1\cong_{\mathrm{Ref}}J\jg{K&n'&m'\\d\cdot 1&n'&m'}\text{~and~}W_2\cong_{\mathrm{Ref}}J\jg{K&p'&q'\\d\cdot 1&p'&q'}.\]
In order to properly write $W_1,W_2$ as reduced $J$-groups, we may need to remove columns for instance if $n'=1$ or if $m'=1$. However, after these reductions, applying Proposition \ref{prop:classif_red_j-group} yields that the number of elements in $(n',m')$ which are equal to $1$ is equal to the number of elements in $(p',q')$ which are equal to $1$. Moreover, since $n'\leqslant m'$ and $p'\leqslant q'$, we deduce that in each case we have $n'=p'$ and $m'=q'$. Since $d=d'$, we deduce that $n=p$ and $m=q$. In other words $W_1=W_2$.

\noindent\underline{Case 2 : $W_1=\B^*(n,m;K,c),W_2=\B(p,q;L)$.}\\First, we have $d+1=|\Hh(W_1)/W_1|=|\Hh(W_2)/W_2|=d'$ by  Proposition \ref{prop:ref_iso_implies_torsion_equal} and Lemma \ref{lem:torsion_multiset_of_torsion_quotients}. Moreover, these two results also give that $T(W_1)=\muls{K,c}=\muls{L}=T(W_2)$.

Now, by Corollary \ref{cor:embedding_torsion_quotients}, we have
\[W_1\cong_{\mathrm{Ref}}J\jg{K&n'&cm'\\d\cdot 1&n'&m'}\text{~and~}W_2\cong_{\mathrm{Ref}}J\jg{K&p'&q'\\d\cdot 1&p'&q'}.\]
If $m'=1$ (i.e. if $m$ divides $n$), then $W_1\cong_{\mathrm{Ref}}\B(n+\frac{n}{m},m+1;(K,c))$ by Lemma \ref{lem:iso_exceptionnel} and we are in Case 1. 

If $m'>1$, then since $c\geqslant 2$, the column $\binom{cm'}{m'}$ then has two distinct coefficient, and the bottom coefficients is different from 1. Since no column in the defining matrix of $W_2$ satisfies this condition (even if $p'=1$ or $q'=1$), then Proposition \ref{prop:classif_red_j-group} implies that $W_1$ and $W_2$ cannot be reflection isomorphic.

\noindent\underline{Case 3 : $W_1=\B^*_*(n,m;K,b,c),W_2=\B(p,q;L)$.}\\
If $n'=1$, then $W_1\cong_{\mathrm{Ref}}\B^*(n+1,m+\frac{m}{n};(K,b),c)$ by Lemma \ref{lem:iso_exceptionnel} and we are in Case 2. Similarly, if $m'=1$, then $W_1\cong_{\mathrm{Ref}}\B_*(n+\frac{n}{m},m+1;(K,c),b)$ by Lemma \ref{lem:iso_exceptionnel}. By Corollary \ref{cor:swap_parameters_torsion_quotient} we are then also in Case 2.

Now, by Corollary \ref{cor:embedding_torsion_quotients}, we have
\[W_1\cong_{\mathrm{Ref}}J\jg{K&bn'&cm'\\d\cdot 1&n'&m'}\text{~and~}W_2\cong_{\mathrm{Ref}}J\jg{L&p'&q'\\d'\cdot 1&p'&q'}.\]
Since $b,c\geqslant 2$, the defining matrix of $W_1$ contains two columns whose coefficients are distinct while the bottom coefficient is different from 1. Since the defining matrix of $W_1$ contains no such column (even in $p'=1$ or $q'=1$), Proposition \ref{prop:classif_red_j-group} implies that $W_1$ and $W_2$ cannot be reflection isomorphic.

\noindent\underline{Case 4 : $W_1=\B^*(n,m;K,c),W_2=\B^*(p,q;L,\gamma)$.}\\
First, we have $d+1=|\Hh(W_1)/W_1|=|\Hh(W_2)/W_2|=d'+1$ by  Proposition \ref{prop:ref_iso_implies_torsion_equal} and Lemma \ref{lem:torsion_multiset_of_torsion_quotients}. Moreover, these two results also give that $T(W_1)=\muls{K,c}=\muls{L,\gamma}=T(W_2)$.

Now, by Corollary \ref{cor:embedding_torsion_quotients}, we have
\[W_1\cong_{\mathrm{Ref}}J\jg{K&n'&cm'\\d\cdot 1&n'&m'}\text{~and~}W_2\cong_{\mathrm{Ref}}J\jg{L&p'&\gamma q'\\d\cdot 1&p'&q'}.\]
If $m'=1$, then $W_1\cong_{\mathrm{Ref}}\B(n+\frac{n}{m},m+1;(K,c))$ by Lemma \ref{lem:iso_exceptionnel} and we are in Case 2. Similarly, if $q'=1$, then $W_2\cong_{\mathrm{Ref}}\B(p+\frac{p}{q},q+1;(L,\gamma))$ by Lemma \ref{lem:iso_exceptionnel} and we are also in Case 2. We can then assume that $m',q'>1$.

Since $c\geqslant 2$ (resp. $\gamma\geqslant 2$), the column $\binom{cm'}{m'}$ (resp. $\binom{\gamma q'}{q'})$ is the only column in the defining matrix of $W_1$ (resp. of $W_2$) whose coefficients are distinct while the bottom coefficient is different from $1$. By Proposition \ref{prop:classif_red_j-group}, we then have $m'=q'$ and $c=\gamma$. Moreover, since $d=d'$, we have $m=q$. 

If $n'=1$, then the defining matrix of $W_1$ contains $d+1$ columns. By Proposition \ref{prop:classif_red_j-group}, we must then have $p'=1$. If $n'>1$, then $\binom{n'}{n'}$ is the only column in the defining matrix of $W_1$ whose coefficients are equal. Since the only possible column in the defining matrix of $W_2$ whose coefficients are equal is $\binom{p'}{p'}$, Proposition \ref{prop:classif_red_j-group} implies that $n'=p'>1$. In each case, we have $p'=n'$ and $p=n$ since $d=d'$. In other words $W_1=W_2$.

\noindent\underline{Case 5 : $W_1=\B^*_*(n,m;K,b,c),W_2=\B^*(p,q;L,\gamma)$.}\\If $n'=1$, then $W_1\cong_{\mathrm{Ref}}\B^*(n+1,m+\frac{m}{n};(K,b),c)$ by Lemma \ref{lem:iso_exceptionnel} and we are in Case 4. If $m'=1$, then $W_1\cong_{\mathrm{Ref}}\B_*(n+\frac{n}{m},m+1;(K,c),b)$ by Lemma \ref{lem:iso_exceptionnel}. By Corollary \ref{cor:swap_parameters_torsion_quotient} we are then also in Case 4. If $q'=1$, then $W_2\cong_{\mathrm{Ref}}\B(p+\frac{p}{q},q+1;(L,\gamma))$ by Lemma \ref{lem:iso_exceptionnel} and we are in Case 3. We can then assume that $n',m',q'>1$.

Now, by Corollary \ref{cor:embedding_torsion_quotients}, we have
\[W_1\cong_{\mathrm{Ref}}J\jg{K&bn'&cm'\\d\cdot 1&n'&m'}\text{~and~}W_2\cong_{\mathrm{Ref}}J\jg{L&p'&\gamma q'\\d'\cdot 1&p'&q'}.\]
Since $b,c\geqslant 2$, the defining matrix of $W_1$ contains two columns whose coefficients are distinct while the bottom coefficient is different from 1. Since $\gamma\geqslant 2$, the defining matrix of $W_1$ contains one such column (even in $p'=1$), Proposition \ref{prop:classif_red_j-group} implies that $W_1$ and $W_2$ cannot be reflection isomorphic.

\noindent\underline{Case 6 : $W_1=\B^*_*(n,m;K,b,c),W_2=\B^*_*(p,q;L,\beta,\gamma)$.}\\First, we have $d+2=|\Hh(W_1)/W_1|=|\Hh(W_2)/W_2|=d'+2$ by Proposition \ref{prop:ref_iso_implies_torsion_equal} and Lemma \ref{lem:torsion_multiset_of_torsion_quotients}. Moreover, these two results also give that $T(W_1)=\muls{K,b,c}=\muls{L,\beta,\gamma}=T(W_2)$. Applying Corollary \ref{cor:swap_parameters_torsion_quotient}, we can assume that $n\leqslant m$ and $p\leqslant q$. 

Now, by Corollary \ref{cor:embedding_torsion_quotients}, we have
\[W_1\cong_{\mathrm{Ref}}J\jg{K&bn'&cm'\\d\cdot 1&n'&m'}\text{~and~}W_2\cong_{\mathrm{Ref}}J\jg{L&\beta p'&\gamma q'\\d\cdot 1&p'&q'}.\]
If $n'=1$, then $W_1\cong_{\mathrm{Ref}}\B^*(n+1,m+\frac{m}{n};(K,b),c)$ by Lemma \ref{lem:iso_exceptionnel} and we are in Case 5. If $m'=1$, then $n'\leqslant m'$ is also equal to 1 and we are also in Case 4. If $p'=1$ or if $q'=1$, then exchanging the roles of $W_1$ and $W_2$ also brings us back to Case 4. We can then assume that $n',m',p',q'>1$.

Since $b,c\geqslant 2$ (resp. $\beta,\gamma\geqslant 2$), the columns $\binom{bn'}{n'}$,$\binom{cm'}{m'}$ (resp. $\binom{\beta p'}{p'}$,$\binom{\gamma q'}{q'}$) are the only columns in the defining matrix of $W_1$ (resp. of $W_2)$ whose coefficients are distincts while the bottom coefficient is different from $1$. By Proposition \ref{prop:classif_red_j-group}, and since $n'\leqslant m'$, $p'\leqslant q'$, we then have $n'=p'$, $m'=q'$. Since $d=d'$, we deduce that $n=p$ and $m=q$. Moreover, we also have $bn'= \beta p'$ and $cm'=\gamma q'$, which implies $b=\beta$ and $c=\gamma$. In other words $W_1=W_2$.

Since Cases 0 to 6 cover all possibilities, this finishes the proof of Theorem \ref{theo:classif_tors_quot_j-braid}.

\vspace{1ex}We finish this section with a consequence of Theorem \ref{theo:classif_tors_quot_j-braid} on the classification of $J$-braid groups.

Considering the isomorphisms given in Corollary \ref{cor:permutation_parameters_torsion_quotient}, Corollary \ref{cor:swap_parameters_torsion_quotient} and Lemma \ref{lem:iso_exceptionnel} when all the torsions are infinite, we obtain the following result:

\begin{corollary}\label{cor:classif_j-braid_group}
The reflection isomorphism relation on the set of $J$-braid groups is generated by the following relations
\begin{enumerate}
    \item $\B_*^*(n,m)\cong_{\mathrm{Ref}}\B_*^*(m,n)$; $\B^*(n,m)\cong_{\mathrm{Ref}}\B_*(n,m)$; $\B(n,m)\cong_{\mathrm{Ref}}\B(m,n)$,
    \item $\B^*(n,m)\cong_{\mathrm{Ref}} \B(n+\frac{n}{m},m+1)$ if $m$ divides $n$,
    \item $\B_*^*(n,m)\cong_{\mathrm{Ref}} \B_*(n+\frac{n}{m},m+1)$ if $m$ divides $n$,
    \item $\B_*(n,m)\cong_{\mathrm{Ref}} \B(n+1,m+\frac{m}{n})$ if $n$ divides $m$,
    \item $\B_*^*(n,m)\cong_{\mathrm{Ref}} \B^*(n+1,m+\frac{m}{n})$ if $n$ divides $m$,
    \item $\B(2,2)\cong_{\mathrm{Ref}}\B^*(1,m)\cong_{\mathrm{Ref}}\B_*(n,1)$,
    \item $\B(1,n)\cong_{\mathrm{Ref}}\B(1,m)$.
\end{enumerate}
\end{corollary}

Since each reflection isomorphism given in Corollary \ref{cor:permutation_parameters_torsion_quotient}, Corollary \ref{cor:swap_parameters_torsion_quotient} and Lemma \ref{lem:iso_exceptionnel} lifts at the level of $J$-braid group, we also have the following corollary:

\begin{corollary}\label{cor:iso_quot_tors_implies_iso_j-braid}
Let $\B,\B'$ be two $J$-braid groups, and let $W,W'$ be respective torsion quotients of $\B$ and $\B'$. If $W$ and $W'$ are reflection isomorphic, then so are $\B$ and $\B'$.
\end{corollary}

\begin{remark}
Specializing Corollary \ref{cor:iso_quot_tors_implies_iso_j-braid} to $p=3$ solves \cite[Conjecture 3.2.3]{TheseIgor}.
\end{remark}

\section{Seifert links}\label{sec:seifert_links}
\subsection{Reminders on Seifert links}
So far our point of view on $J$-braid groups and their torsion quotients has been purely combinatorial and given in terms of generators and relations. A topological interpretation of $J$-braid groups is given in \cite{Necklace} as link groups of particular links called \emph{torus necklaces}, which generalize the classical torus knots. Torus necklaces make up the majority of a family of links called \emph{Seifert links}. We give some reminders on Seifert links here, following the exposition of \cite{Necklace}.

In this section, we fix two positive integers $n,m$.

\begin{figure}[ht]
    \centering
    \begin{subfigure}{0.2\textwidth} 
        \centering
        \includegraphics[width=\textwidth]{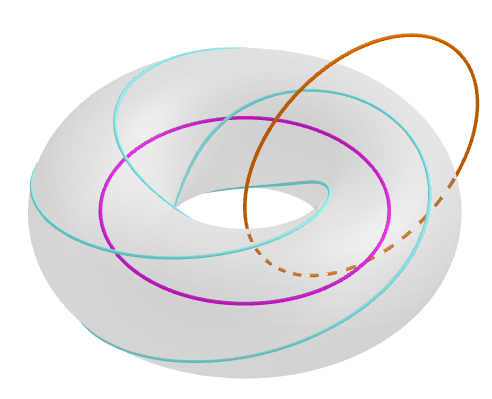}
        \caption{$\mathrm L_*^*(3,4)$}
 
    \end{subfigure}
        \begin{subfigure}{0.2\textwidth} 
        \centering
        \includegraphics[width=\textwidth]{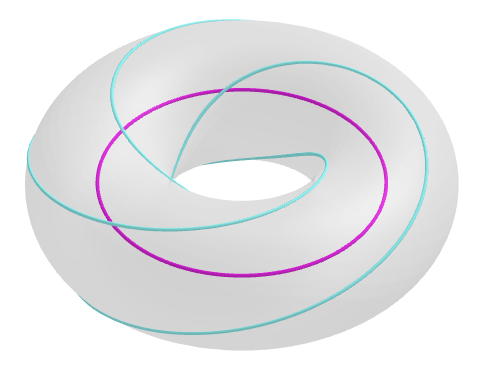}
        \caption{$\mathrm L_*(3,4)$}

    \end{subfigure}
        \begin{subfigure}{0.2\textwidth} 
        \centering
        \includegraphics[width=\textwidth]{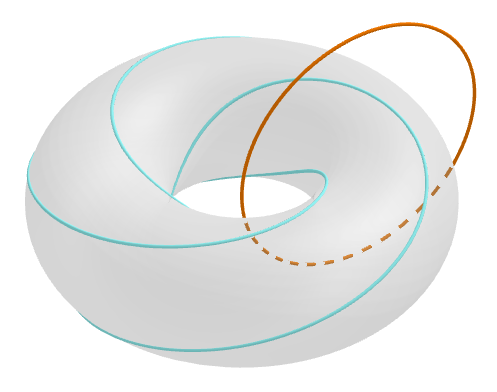}
        \caption{$\mathrm L^*(3,4)$}

    \end{subfigure}
     \begin{subfigure}{0.2\textwidth} 
        \centering
        \includegraphics[width=\textwidth]{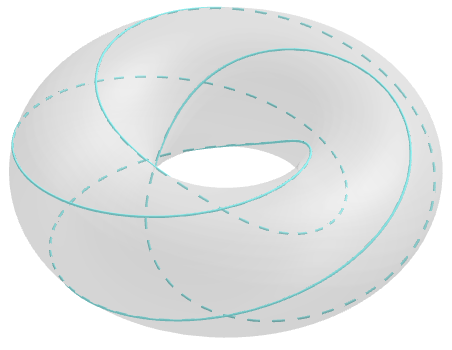}
        \caption{$\mathrm L(3,4)$}

    \end{subfigure}
\caption{The torus necklaces $\mathrm L_*^*(3,4)$, $\mathrm L_*(3,4)$, $\mathrm L^*(3,4)$ and $\mathrm L(3,4)$}
\label{L23IntroFR}
\end{figure}

We start by defining torus necklaces following (for a precise definition as the closure of particular braids, see \cite[Definition 4.2]{Necklace}).
\begin{itemize}
    \item The \emph{torus link} $L(n,m)$ is the closure of the braid $(\sigma_1\cdots \sigma_{n-1})^m$, where $\sigma_1,\ldots,\sigma_{n-1}$ denote the Artin generators of the braid group on $n$ strands. (see Figure \ref{L23IntroFR}).
    \item The link $L_*(n,m)$ is the disjoint union of $L(n,m)$ with a circle going around the internal heart of the torus (see Figure \ref{L23IntroFR}).
    \item The link $L^*(n,m)$ is the disjoint union of $L(n,m)$ with a circle going around the external heart of the torus (see Figure \ref{L23IntroFR}).
    \item The link $L^*_*(n,m)$ is the disjoint union of $L(n,m)$ with two circles, one around the internal heart of the torus, and one around its external heart (see Figure \ref{L23IntroFR}).
\end{itemize}

\begin{figure}[ht]
        \centering
        \includegraphics[width=0.4\textwidth]{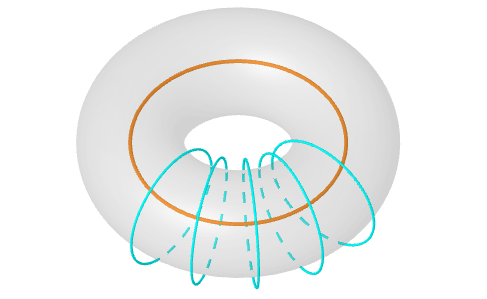}
        \caption{The keychain $\mathrm K(5)$}
 \label{TrousseauExampleIntro}
\end{figure}

On top of torus necklaces, we are interested in keychain links: The \emph{keychain} $K(k)$ is the disjoint union of $k$ circles around the internal heart of the torus (see Figure \ref{TrousseauExampleIntro}).

The study of these links is justified by the following result which classifies all links whose associated group admits a nontrivial center:

\begin{theorem}\cite{BZTorus} and \cite[Theorem 1]{SeifertLinks}
\begin{enumerate}
    \item The knot group of a knot $K$ has a nontrivial center if and only if $K$ is isotopic to a torus knot.
    \item The link group of a link $L$ has a nontrivial center if and only if $L$ is isotopic to a torus necklace or to a keychain.
\end{enumerate}
\end{theorem}

The links which appear in the second statement of the above theorem are called \emph{Seifert links} since they are precisely the links whose complement in $\mathbb{S}^3$ are Seifert fibered spaces. Not only are Seifert links interesting because the center of their link groups is nontrivial, they are also central objects in the study of the so-called JSJ decomposition of link complements (see \cite{Budney} for a discussion of the JSJ decomposition in the context of link). The relationship between Seifert links and $J$-braid groups appeared in \cite{Necklace} by the second author:

\begin{theorem}\cite[Theorem 4.7]{Necklace}\label{theo:j-braid_groups_are_link_groups}
\begin{enumerate}
    \item The link group of $L_*^*(n,m)$ is isomorphic to $\B_*^*(n,m)$.
    \item The link group of $L_*(n,m)$ is isomorphic to $\B_*(n,m)$.
    \item The link group of $L^*(n,m)$ is isomorphic to $\B^*(n,m)$.
    \item The link group of $L(n,m)$ is isomorphic to $\B(n,m)$.
\end{enumerate}
Moreover, under these isomorphisms, meridians of the link group correspond to braid reflections of the $J$-braid group.
\end{theorem}

The statement in \cite[Theorem 4.7]{Necklace} assumes that $n,m$ are coprime, but as pointed out in \cite[Section 5.1]{Necklace}, the isomorphisms of \cite[Theorem 4.7]{Necklace} does not require any coprimality assumption. Moreover, examining the proof shows that the assumptions that $m\geqslant 2$ is not required for the second statement to be true. Likewise for the third and fourth statements.


On top of the torus necklaces groups, let us also mention that the group of the keychain link $K(k)$ is a direct product $F_k\times \Z$ of a free group on $k$ letters with the group of integers. Under this isomorphism, the generators of $F_k\times \Z$ correspond to meridians in the link group.

A classification of Seifert links up to isotopy is given by the following result: 

\begin{proposition}[\textbf{Isotopy classes of Seifert links}]\cite[Proposition 3.5]{Budney}\label{prop:budney_isotopy}\\
The equivalence relation $\sim$ of unoriented isotopy on Seifert links is generated by the relations:
\begin{enumerate}[(a)]
\item $L_*^*(n,m)\sim L_*^*(m,n))$; $L_*(n,m)\sim L^*(m,n)$; $L(n,m)\sim L(m,n)$,
\item $L_*(n,m)\sim L(n+\frac nm,m+1)$ if $m|n$,
\item $L_*^*(n,m)\sim L^*(n+\frac nm,m+1)$ if $m|n$,
\item $L(2,2)\sim L_*(n,1)\sim L^*(1,m)$,
\item $L(1,n)\sim L(1,m)$,
\item $L(2,2)\sim K(1)$,
\item $L(1,1)\sim K(0)$.
\end{enumerate}
\end{proposition}

Combining this result with Corollary \ref{cor:classif_j-braid_group}, we can show that the "topological isomorphism type" of a link group completely characterizes the isotopy type of a Seifert link:

\begin{theorem}\label{theo:isotopy=reflection_iso_link_group}
Let $L,L'$ be two Seifert links, and let $\Ll,\Ll'$ be their respective link groups. The following are equivalent
\begin{enumerate}[(i)]
\item The links $L,L'$ are isotopic,
\item There is an isomorphism $\Ll\to \Ll'$ which maps meridians in $\Ll$ bijectively to meridians in $\Ll'$.
\end{enumerate}
\end{theorem}
\begin{proof}
$(i)\Rightarrow (ii)$ is true in general for all links, so that we only have to prove the converse statement.

First, we assume that $L$ and $L'$ are torus necklaces. In this case $\Ll$ (resp. $\Ll'$) is isomorphic to a $J$-braid group $\B$ (resp. $\B'$) in a way that maps meridians to generating (braid) reflections. Thus statement $(ii)$ implies that $\B$ and $\B'$ are reflection isomorphic.  We can then apply Corollary \ref{cor:classif_j-braid_group}, which implies that $\B$ and $\B'$ are related by a sequence of isomorphism appearing in Corollary \ref{cor:classif_j-braid_group}. For each of these isomorphisms, the associated links $L$ and $L'$ are isotopic (compare the isomorpisms of Corollary \ref{cor:classif_j-braid_group} with the isotopies of Proposition \ref{prop:budney_isotopy}). We then have that $L$ and $L'$ are isotopic. 

Now, assume that $L'$ is a keychain link say $K(k)$ with $k\geqslant 0$. Since the link group of $K(k)$ is $\Ll'=F_k\times \Z$ (where the generators are meridians), a generator of the $\Z$ part is a meridian. 

Assume that $L$ is a torus necklace. Let $\B$ be a $J$-braid group isomorphic to $\Ll$ in a way that maps meridians to generating reflections. The isomorphism given by $(ii)$ then sends a generator of $Z(\Ll')$ to a generating reflection $r$ of $\B$. Since $\B$ is reflection isomorphic to a generalized $J$-reflection group, Proposition \ref{prop:centralizer_reflections} (centralizer of a reflection) implies that $\B=C_{\B}(r)$ is abelian. The group $\Ll=F_k\times \Z$ is then also abelian, and $k\leqslant 1$. In this case, the link $L'$ is isotopic to a torus necklace by Proposition \ref{prop:budney_isotopy} and we are back in the first case.

Lastly, assume that $L=K(k')$ is another keychain. We have $\Ll\simeq \Ll'$ if and only if $k=k'$, in which case $L$ and $L'$ are isotopic.
\end{proof}

\subsection{Torsion quotients of Seifert links}
Since for every $J$-braid group $\B$ there is a Seifert link group $\mathcal L$ and an isomorphism $\B\to\mathcal{L}$ sending braid reflections to meridians, it is very natural to consider torsion quotients of Seifert link groups as well.

\begin{definition}[\textbf{Torsion quotient of a link}]\label{def:torsion_quotient_link}
Let $L$ be a link and let $L_1,\dots,L_p$ be its knot components. A torsion quotient of $L$ is a group of the form
\begin{equation}\label{eq:pres_torsion_quotient_link}
    \pi_1(\mathbb{S}^3\backslash L)/\llangle m_1^{k_1},\dots,m_p^{k_p}\rrangle,
\end{equation}
where for each $i\in \llbracket1,p\rrbracket$, $m_i$ is a meridian of $L_i$ and $k_i\in\N_{\geqslant 2}\cup\{\infty\}$. The conjugates of non-trivial powers of the images of $m_1,\dots,m_p$ in this group are still called meridians.\\
Two such torsion quotients are said to be topologically isomorphic if there exists an isomorphism which induces a bijection between the meridians.
\end{definition}

By Corollary \ref{cor:embedding_torsion_quotients} (embedding of torsion quotients) and Theorem \ref{theo:generalized_j_groups_torsion_quotients} (generalized $J$-groups as torsion quotients), the family of torsion quotients of $J$-braid groups coincides with the family of finite-type generalized $J$-group up to reflection isomorphism. By Theorem \ref{theo:j-braid_groups_are_link_groups}, we obtain that the following families of groups coincide up to isomorphism preserving the generators:
\begin{itemize}
    \item Finite-type generalized $J$-groups
    \item Torsion quotients of $J$-braid groups
    \item Torsion quotients of torus necklaces
\end{itemize}

Since we know that the finite groups in the first family are precisely the complex reflection groups of rank 2, we obtain the following result:

\begin{corollary}[\textbf{Finite torsion quotients of Seifert links}]\label{cor:finite_torsion_quotients_seifert}
The family of finite torsion quotients of Seifert links precisely coincides with the family of rank two complex reflection groups.
\end{corollary}

\begin{proof}
By the above discussion, the finite torsion quotients of torus necklaces are the finite generalized $J$-groups, which are the rank two complex reflection groups by Corollary \ref{cor:finite_generalized_j_groups_are_crgs}.\\
It remains to determine the finite torsion quotients of keychain links, which is immediate: if $k\leqslant 2$, every finite quotient torsion of $K(k)$ is a direct product of two cyclic group (which is a rank two complex reflection group) and if $k\geqslant 3$, every quotient torsion of $K(k)$ contains a free product of two nontrivial groups, hence it is infinite.
\end{proof}

It turns out that torsion quotients are extremely rigid, at least for Seifert links, in that they detect isotopy of links:

\begin{theorem}[\textbf{Classification of torsion quotients of Seifert links}]\label{theo:class_torsion_quotient_seifert}~\\
Let $L,L'$ be two Seifert links, and let $W,W'$ be respective torsion quotients of $L$ and $L'$. If $W$ and $W'$ are topologically isomorphic, then $L$ and $L'$ are isotopic.
\end{theorem}
\begin{proof}
Assume that $W$ and $W'$ are topologically isomorphic torsion quotients of $L$ and $L'$. 

If $L$ and $L'$ are both torus necklaces, then the result is a combination of Corollary \ref{cor:iso_quot_tors_implies_iso_j-braid} with Theorem \ref{theo:isotopy=reflection_iso_link_group}. 

Now, assume that $L$ is a torus necklace and that $L'=K(k)$ is a keychain. There is a meridian in $L'$ which is central. The image in $W$ of this meridian is then a central reflection. Proposition \ref{prop:centralizer_reflections} then implies that $W=C_{W}(r)$ is abelian. Since $L'$ is a direct product of $\Z$ with a free product of $k$ nontrivial cyclic groups, we obtain that $k\leqslant 1$. In each cases, $L'$ is isotopic to a torus necklace by Proposition \ref{prop:budney_isotopy} and we are back in the first case.

Lastly, assume that both $L=K(k)$ and $L'=K(k')$ are keychains. We write
\[W:=\Z/a_1\Z*\cdots*\Z/a_k\Z\times \Z/p\Z\text{~and~}W':=\Z/b_1\Z*\cdots*\Z/b_{k'}\Z\times \Z/q\Z.\]
We know that $W$ (resp. $W'$) is abelian if and only if $k\leqslant 1$ (resp. $k'\leqslant 1$). Since $W$ and $W'$ are in particular isomorphic, we have $k\leqslant 1$ if and only if $k'\leqslant 1$. If $k=0$, then $W=\Z/p\Z$ admits 1 reflecting hyperplane in the sense of Definition \ref{def:reflecting_hyperplane} (replacing reflections with meridians). If $k=1$, then $W\simeq \Z/a_1\Z\times \Z/p\Z$ admits two reflecting hyperplanes. We then have $k=k'$ if $k=0$ or if $k=1$, and $L,L'$ are isotopic.

Assume now that $k,k'>1$. The center of $W$ (resp. of $W'$) is the direct factor $\Z/p\Z$ (resp. $\Z/q\Z$). The groups $W/Z(W)$ and $W'/Z(W')$ are the two free products of $k$ and $k'$ nontrivial cyclic groups. Since the minimal number of generators of a free product of $n$ nontrivial cyclic groups is $n$ by Grushko's Theorem, we obtain that $k=k'$ and that $L$ and $L'$ are isotopic.
\end{proof}

\printbibliography
\end{document}